\documentclass[a4paper,10pt]{article}
\usepackage{amsmath,amssymb,amsthm}
\usepackage[english]{babel}

\newtheorem{theo}{Theorem}[section]
\newtheorem{defi}[theo]{Definition}
\newtheorem{lemm}[theo]{Lemma}
\newtheorem{prop}[theo]{Proposition}
\newtheorem{coro}[theo]{Corollary}

\newcommand{\Na}{\mathbb N}                   

\newcommand{\Ra}{\mathbb R}                   

\newcommand{\Ca}{\mathbb C}                   
\newcommand{\re}[1]{\mbox{Re} \ #1}    
\newcommand{\im}[1]{\mbox{Im} \ #1}

\newcommand{\scal}[1]{\langle #1 \rangle}

\newcommand{\csubset}{\Subset}

\newcommand{\Schn}[1]{{\bf S_{#1}}}

\newcommand{\finpreuve}{\hfill $\Box$}

\newcommand{\name}{$\underline{\qquad \qquad}$}
\newcommand{\etan}{$ \& $ }

\newcommand{\refe}[1]{\ref{#1}}
\newcommand{\reff}[1]{(\ref{#1})}

\newcommand{\e}{\epsilon}

\begin{document}



\author{{\sc J.M. Bouclet} \\
Universit\'e de Lille 1 \\
 UMR CNRS 8524, \\
59655 Villeneuve d'Ascq \\
 Jean-Marc.Bouclet@agat.univ-lille1.fr}
\title{ \bf \sc Generalized scattering phases for
asymptotically hyperbolic  manifolds}

\maketitle

\begin{abstract} We prove two asymptotic expansions of the 
generalized scattering phases. These phases are generalizations
of the Birman-Krein spectral shift function associated to pairs 
of perturbations of the Laplacians of asymptotically hyperbolic manifolds.
The first expansion, of `heat type', holds for all `long range' metric 
perturbations of the Laplacian, whereas the second one is shown under a non 
trapping condition.

\end{abstract}

\setcounter{section}{0}

\section{Introduction and results} \label{introduction}
\setcounter{equation}{0}

\subsection{Introduction}
In this paper,  we define and study some properties of the
generalized scattering phases associated to a pair of self-adjoint
elliptic differential operators $(P_0,P_1)$ on an asymptotically
hyperbolic manifold $X$ of dimension $n$. The manifolds and operators that we 
consider are of the same type as those considered in \cite{MaMe1,FrHi0} or in 
many other related papers, especially \cite{JoSa}. Typically, 
$ P_0 $ can be the Laplacian of the quotient of the hyperbolic space $ {\mathbb
H}^n $ by some discrete group of isometries, with infinite volume,
and $ P_1 $  a non compactly
supported perturbation of $ P_0 $. The precise definition of the operators is
given in subsection $ \refe{notation} $.
Let us first explain the terminology.

The generalized scattering phase of order $q \in \Na$ is the
distribution $ \xi_{\bf q} $ defined by
\begin{eqnarray}
 \int_{\Ra} f (\lambda) \ \mbox{d} \xi_{\bf q}(\lambda) = \
\mbox{tr} \left( f (P_1) - \sum_{j=0}^{q-1} \frac{1}{j !}
\frac{d^j}{d \e^j} f (P_{\e})_{| \e = 0} \right), \qquad f \in
C_0^{\infty} (\Ra) \label{traq}
\end{eqnarray}
provided, of course, that the right hand side makes sense. Here  $
P_{\e} = P_{0} + \e (P_1 - P_0) $ and  the left
hand side stands for $ \scal{\xi_{\bf q}^{\prime},f} $ if
$\scal{.,.}$ is the duality between test functions and
distributions. This equality actually defines  the derivative of
$\xi_{\bf q}$ and we choose for $ \xi_{\bf q}$ the unique primitive of
$ \xi_{\bf q}^{\prime}$ vanishing near $-
\infty$, which is possible if $P_0$ and $P_1$ are semi-bounded
from below.

 The distribution $\xi_{\bf 1}$ is well
known: it is  the spectral shift function introduced by Birman and
Kre$\breve{\i}$n \cite{BiKr1,Yafa1}. We recall that
$ \xi_{\bf 1} $ is, in general, a measurable
function which is defined if
$ (P_1 + i)^{-N} - (P_0+i)^{-N} \in \Schn{1} $
 for some $N$ large enough. Here $\Schn{1}$ is the set of trace class operators.
 
 The point of considering
$ \xi_{\bf q} $ for $q \geq 2$ is that we can relax this trace class condition
 and consider operators such that $ (P_1 +
i)^{-N} - (P_0+i)^{-N} \in \Schn{q} $, the Schatten class of order
$q$ (for instance the Hilbert-Schmidt class if $q=2$). We
furthermore emphasize that, if $ q \geq 2$, $ \xi_{\bf q} $ is a priori a
distribution  and showing that it is a (smooth or continuous
or even measurable) function on the continuous spectrum is not trivial.
The distributions $ \xi_{\bf q} $ have been introduced by Koplienko
\cite{Kopl1} for bounded or $1$-dimensional Schr\"odinger
operators and studied in the higher dimensional case by the author
(see \cite{Bouc2,Bouc4} for more details). We also quote
a similar approach considered in \cite{HiPo1}.  In \cite{Bouc2,Bouc4},
$\xi_{\bf q}$ was called {\it spectral distribution} by analogy
with  {\it spectral function}, however the name {\it scattering phase} is
natural as well since we have shown that
\begin{eqnarray}
 \xi_{\bf q}^{\prime}(\lambda) = \lim_{\delta \downarrow 0}
 \frac{d}{d \lambda} \
\mbox{arg} \ \mbox{det}_{\bf q} \left( (P_1 - \lambda - i\delta)
(P_0 - \lambda - i \delta)^{-1} \right) \label{detq} .
\end{eqnarray}
The precise meaning of formula $ \reff{detq} $, which is well known for $ q = 1 $ if
$(P_1 - P_0)(P_0 - z)^{-1}$ is trace class
\cite{BiKr1,Yafa1}, is explained in
\cite{Bouc4}. We just specify that $\mbox{det}_{\bf q}$ coincides
in many cases with  the standard Fredholm determinant
$\mbox{Det}_{\bf q}$ defined for  perturbations of identity by
elements of $ \Schn{q} $ (see \cite{GoKr1,Yafa1}) and that we need such
an extension of $\mbox{Det}_{\bf q}$ since $ (P_1 - P_0)(P_0
-z)^{-1}  $ is not necessarily compact in general. We also refer to the recent paper 
\cite{BJP1} where similar determinants with $ q = 1 $ are studied for hyperbolic 
surfaces.

\smallskip

We now recall the definition of an asymptotically hyperbolic
manifold. A complete non compact
Riemannian manifold $(X,G)$, without boundary, is asymptotically hyperbolic if
it is  isometric, outside a compact set, to
\begin{eqnarray}
 \left( (R , + \infty) \times Y , d r^2 + e^{2r} g (e^{-r}) \right) .
\label{definitionmetrique}
\end{eqnarray}
Here $ Y $ is a connected compact manifold without boundary, and
 $g (x)$ is
a family of metrics on $Y$ depending smoothly on $x \in
[0,e^{-R})$. More precisely, if $ S^2 T^* Y $ is the vector bundle of bilinear
symmetric forms on $ Y $ and $ \Gamma (S^2 T^* Y  ) $  is the space of its 
smooth sections, we assume that
$ g \in C^{\infty}( [0,e^{-R}) , \Gamma (S^2 T^* Y  )) $ and of course
that $ g (x)|_{T_p Y \times T_p Y} $ is positive definite
for each $ p \in Y $ and $ x \in [0,e^{-R}) $.  
We furthermore equip  $Y$ with the metric 
$g (0) $.  We could actually consider
manifolds with finitely
many such ends, i.e. with $ (X, G) $ isometric, outside a compact set, 
to a finite union of  manifolds like $ \reff{definitionmetrique} $ but we restrict our
attention to the one end case for notational convenience. 
As explained in \cite{JoSa,Melr0}, such a
manifold $(X,G)$ can be obtained from a compact manifold with
boundary $(Z,\widetilde{G})$, with a boundary defining function
$x$ such that $\partial Z = Y = \{ x = 0 \}$, by setting $X =  \{
x > 0 \}$, the interior of $Z$, equipped with the metric
$\widetilde{G}/x^{2}$. 
These manifolds are also called conformally compact
manifolds and the most basic example is the hyperbolic space $
{\mathbb H}_n $.

\medskip

The results of this paper are two asymptotic formulas for $ \xi_{\bf q} $
which are similar to the heat expansion
and Weyl formula for the eigenvalues counting function on a compact manifold.
More precisely, these asymptotics are of the same type as those obtained
 in Euclidean scattering by
\cite{PePo1,MaRa1,Robe2,Robe4,Chri0,VaWa1,KoPu1} for $ q = 1$ and
 \cite{Bouc2,Bouc4} for $ q \geq 2 $. 
 See also \cite{Chri0,Chri1,Carr1,FrHi1} and \cite{GuZw2} in more geometric frameworks. 
The scattering phases are natural and basic tools of scattering theory in view of
their relation with time delay \cite{Robe4}, relative scattering determinants 
\cite{Bouc2,Bouc4,PeZw2} or resonances. More specifically, their
asymptotic behavior
 is of interest for several reasons such as  relative index theory \cite{BMS1}, 
 trace formulas \cite{Coli1,Guil1,Bouc2,BeRy1} 
or Breit-Wigner formula \cite{GMR1,PeZw1,BrPe1}.

The most popular scattering phase is Birman-Kre$\breve{\i}$n's function 
$ \xi_{\bf 1} $ but its use leads to restrictions on the pairs of operators as already
mentionned. For instance in \cite{BJP1}, the authors are able to define the determinant 
$ \mbox{det}_{\bf 1} $ for a pair of operators which are, up to a unitary
transform, the Laplacians associated with two metrics $G_0$ and $G_1$ as above
for which $ g_1 (x) - g_0(x) = {\mathcal O}( x^2)$. This last condition 
implies that 
$ (P_0 + i)^{-N} - (P_1+ i)^{-N} $ is of trace class since they work in dimension
$ 2 $, but in higher dimension, the difference 
$ (P_0 + i)^{-N} - (P_1+ i)^{-N} $ is not trace class, in general,
under the sole condition $ g_1 (x) - g_0 (x) = {\mathcal O}(x^2) $. Our theorem
$ \refe{theo1} $ combined with the method of \cite{Bouc4} proves directly
the existence of $\mbox{det}_{\bf q}((P_1-z)(P_0-z)^{-1})$ for  $q \geq n$ for 
any $n$ under the weaker condition that $ g_1 (x) - g_0 (x) = {\mathcal O}(x)$.

For the proof of theorem $ \refe{theo2} $, we adapt the method that we used in \cite{Bouc4},
namely a refined analysis of Isozaki-Kitada's construction \cite{IsKi1}. We recall the
principle of this method, which has only been used in the Euclidean context so far, in 
subsection $ \refe{IsozakiKitadaformel} $. We devote section
$\refe{classical}$ to the  relevant estimates on geodesics in the
hyperbolic framework and the explicit construction is given in subsection $ \refe{main1}$.

\subsection{Notations} \label{notation}
For any  $C^q$ function $T_{\e}$ of the real variable $\e \in
[0,1]$, scalar or vector valued, we set
\begin{eqnarray}
 \left[  T_{\e} \right]_{\bf q} \  = \  T_{1} - \sum_{j=0}^{q-1}
\frac{1}{j!} \frac{d^j}{d \e^j} T_{\e | \e= 0}
 \ = \  \frac{1}{(q-1)!} \int_0^1 (1-\e)^{q-1} \frac{d^q}{d \e^q} T_{\e} \
\mbox{d} \e .
\end{eqnarray}
We will also use the notation
$$ \left\{ T_{\e}  \right\}_{\bf q} \ = \
\frac{1}{(q-1)!} \int_0^1 (1-\e)^{q-1} \frac{d^{q-1}}{d \e^{q-1}}
T_{\e} \ \mbox{d} \e .  $$ Note that, if ${\mathcal T}_{\e}$ is a
primitive of $T_{\e}$, we have $[{\mathcal T}_{\e}]_{\bf q}=\{
T_{\e} \}_{\bf q} $.

We will have to consider distributions smooth with respect to a
parameter. We shall say that a family of distributions $u_{\tau}
\in {\mathcal S}^{\prime}(\Ra)$ is $C^j$ with respect to $\tau \in
J$, an interval of $\Ra$, if for all $f \in {\mathcal S}(\Ra)$ the
function $ \scal{u_{\tau},f}$ is $C^j$ on $J$. In particular, $
\partial_{\tau}^j u_{\tau}$ and $ \int_J u_{\tau}  \mbox{d} \tau
$ are defined by
$$ \scal{ \partial_{\tau}^j u_{\tau} ,  f } = \partial_{\tau}^j \scal{u_{\tau},f}, \qquad
\left\langle \int_{J} u_{\tau} \ \mbox{d} \tau , f \right\rangle =
\int_J \scal{u_{\tau},f} \ \mbox{d} \tau .
$$
Another useful distributional notation is the following. If $H$ is
a self-adjoint operator on  a separable Hilbert space ${\mathcal H
}$ and if $T$ is an operator acting on (some subspace) of
${\mathcal H}$ such that the map $f \mapsto \mbox{tr} \left( f(H)T
\right)$, with $f \in {\mathcal S}(\Ra)$, defines a distribution,
that is if the trace is well defined and depends continuously on
$ f \in {\mathcal S}(\Ra)$, then we shall write this distribution
$$ \mbox{tr} \left( \frac{\partial E}{\partial \mu} T \right) $$
if $E (\mu) = E (- \infty,\mu)$ is the spectral projection of $H$
on $(- \infty , \mu)$.

We will also use extensively Schatten classes $ \Schn{q} $ of real
order $ q \geq 1 $. We simply recall that, by definition,  a
bounded operator $ A $ on a separable Hilbert space $ {\mathcal H}
$ belongs to $ \Schn{q} = \Schn{q}({\mathcal H}) $ if  $ |A|^q =
(A^*A)^{q/2} $ is trace class and that the norm $ ||.||_{\bf q} $
on $ \Schn{q} $ is defined by $ ||A||_{\bf q}^q = \mbox{tr}
(|A|^q) $. We will need the following H\"older type estimates
\begin{eqnarray}
|| A B ||_{\bf q} \leq  || A ||_{\bf q_1} ||B||_{\bf q_2}, \qquad
q^{-1} = q_1^{-1} + q_2^{-1} \label{HolderSchatten}
\end{eqnarray}
for all $A \in \Schn{q_1} $ and $  B \in \Schn{q_2} $. This
estimate still holds if $ A $ (resp. $ B $) is bounded with $q_1 =
\infty $ (resp. $q_2 = \infty$), which is consistent with the
notation $ ||.||_{\infty} $ for the operator norm on  $ {\mathcal
H }$. For a more general presentation of Schatten classes, we
refer to \cite{GoKr1} and \cite{Yafa1}.

On the manifold $X $, we will mainly work near infinity. If $ y =
( y_1, \cdots , y_{n-1} ) $ are coordinates on $ Y $ defined on $
U_{n-1} \subset Y $, such that $ y : U_{n-1} \rightarrow \Omega
\subset \Ra^{n-1} $ is a diffeomorphism, then $ r ,y_1 , \cdots ,
y_{n-1} $ are coordinates on an subset $ U_n \subset X $
diffeomorphic to $ (R,\infty) \times \Omega $. We call $
(R,\infty) \times \Omega $ a  chart at infinity. We can write $ Y
$ as the finite union of open sets $ U_{n-1} $ and thus we get an
atlas of a neighborhood of infinity on $ X $ given by a finite
number of charts at infinity. If necessary, we can choose each $
\Omega $ to be convex, since we fix the coordinates once for
all in each chart. We assume that, in each chart at infinity, the
volume density giving the $L^2$ structure on $X$ can be written
\begin{eqnarray}
 \mbox{d}vol = a (r,y)| d r \wedge d y_1 \wedge \cdots \wedge d
y_{n-1} | \label{formevolume}
\end{eqnarray}
 with $ a $ smooth and bounded on $  (R,\infty) \times
\Omega $. In particular, this implies that the pullback on $ X $
of any function $ u \in L^2 ((R,\infty) \times \Omega ,
\mbox{d}r\mbox{d}y) $, supported in the chart, 
belongs to $ L^2 (X) = L^2 (X,\mbox{d}vol) $.

We now describe the operators $P_0$ and $ P_1 $. Let us start with
an example.  The expression of the Laplace Beltrami operator
associated to $ \reff{definitionmetrique} $ in a chart at infinity
is
$$ \Delta^G = - \partial_r^2 - e^{-2r} \Delta^{g(e^{-r})} - \left( n-1 - 
e^{-r} \frac{\partial_x \mbox{det}  g(x)}{2 \mbox{det} g(x)}|_{x=e^{-r}} \right) \partial_{r}  $$
where $ \Delta^{g(x)} $ is the expression of the (negative) Laplacian
of $ Y $ associated to the metric $ g (x) $. Furthermore $ \mbox{d}vol(G)= e^{(n-1)r} \sqrt{\mbox{det} \ g(e^{-r},y)} |dr \wedge d y_1 \wedge
\cdots \wedge d y_{n-1}| $ is  the volume form induced
by $G$, 
hence 
\begin{eqnarray}
e^{(n-1)r/2} \Delta^G e^{-(n-1)r/2} = - \partial_r^2 - 
e^{-2r} \Delta^{g(e^{-r})} + c_n^2 + e^{-r}
 \frac{\partial_x \mbox{det} g(x)}{2 \mbox{det}  g(x)}|_{x=e^{-r}} (\partial_r - c_n) , 
 \label{operateur}
\end{eqnarray}
with $ c_n = (n-1)/2 $,
is selfadjoint w.r.t the density $ \sqrt{\mbox{det} \ g(e^{-r},y)} |dr \wedge d y_1 \wedge
\cdots \wedge d y_{n-1}|  $ which is of the form $ \reff{formevolume} $. Guided by
 $ \reff{operateur} $, we shall consider operators $P_0,P_1$ which are both 
second order elliptic differential operators,
symmetric w.r.t to $ \mbox{d}vol $, whose expressions in each chart at
infinity are
\begin{eqnarray}
  P_j = - \partial_r^2 - g_j(e^{-r},y,e^{-r} \partial_y) 
  + e^{-  r}  \sum_{|\alpha|+l \leq 1} v_j^{\alpha,l}(e^{-r},y) 
 (e^{-r} \partial_y)^{\alpha}\partial_r^l , \qquad j =0,1. \label{ordreundiff}
\end{eqnarray}
Here $  g_j(x,y,\eta)$ is the principal symbol
of $\Delta^{g
_j(x)}$, $ j= 0,1$, that is the expression of the metric $g_j$ on the fibers 
of $T^*Y$ with coordinates $\eta_1,\cdots,\eta_{n-1}$ dual of  
$ y_1 , \cdots , y_{n-1} $, and the functions $v^{\alpha,l}_j(x,y)$ are 
smooth and bounded on $[0,e^{-R})\times \Omega$. The fact that $ P_1 $ is a
perturbation of $ P_0 $ is reflected by the assumption that 
\begin{eqnarray}
 g_1|_{x=0} = g_0|_{x=0} . \label{recollement}
\end{eqnarray}
We will use extensively
 the principal symbol of $ P_0 + \e (P_1 - P_0) $, denoted by $ p_{\e} $, which
 has the form 
$$ p_{\e} (r,y,\rho,\eta) = \rho^2 + e^{-2r} g_{\e}(r,y,\eta) $$
where $ g_{\e}(r,y,\eta)  = g_0(e^{-r},y,\eta) + \e (g_1 -g_0)(e^{-r},y,\eta) $ 
and $ \rho $ is dual variable of $r$. Note that
 $g_1(x,y,\eta) - g_{0}(x,y,\eta)
= {\mathcal O}(x)|\eta|^2$ by $ \reff{recollement} $.
\subsection{Results}
In the next two theorems,  $P_0$ and $P_1$ are two operators as described above
such that $ \reff{recollement} $ holds. We recall that $ P_{\e} = P_0 + \e 
(P_1 - P_0) . $
\begin{theo}   \label{theo1} i) For all $q \geq n$ and 
$ f \in {\mathcal S}(\Ra)$, $[f(P_{\e})]_{\bf q}$ is trace class and 
there exists a unique $\xi_{\bf
q} \in {\mathcal S}^{\prime}(\Ra) $ which vanishes below 
$\inf (\sigma(P_0) \cup \sigma(P_1) )$ such that
$$ \int_{\Ra} f (\lambda) \ \emph{d} \xi_{\bf q}(\lambda) = \
\emph{tr} \left[ f (P_{\e}) \right]_{\bf q}. $$
ii) The Laplace transform of $\xi_{\bf q}^{\prime}$ has a complete asymptotic
expansion as $t \downarrow 0+$, namely
\begin{eqnarray}     
 \emph{tr} \left[ e^{-t P_{\e} } \right]_{\bf q} \sim t^{-n/2} 
\sum_{k \geq 0} a_k t^k , \qquad 
                       \mbox{with} \  \
 a_0 =     \Gamma \left(\frac{n}{2}+1 \right) (2 \pi)^{-n} \omega_n
  \int_X [\emph{d}vol_{\e}]_{\bf q}   .
\label{heatexpansion}     
\end{eqnarray}     
Here $\emph{d}vol_{\e}$ is volume density obtained naturally from $p_{\e}$, that
is in local coordinates
$$ \emph{d}vol_{\e} =  e^{(n-1)r} \emph{det} \left( \partial_{\eta}^2
g_{\e}(r,y,\eta) /2 \right)^{- 1/2} |dr \wedge d y_1 
\wedge \cdots \wedge d y_{n-1}| . $$ 
The other coefficients $ a_1 , a_2 , \cdots $ 
can be expressed as integrals of functions of the symbols of $P_{\e}$.
\end{theo}

\noindent {\it Example.} If $(X,G_0)$ and $(X,G_1)$ are of the form $
\reff{definitionmetrique} $, outside a compact set, associated respectively
to $g_0$ and $g_1$ satisfying $
\reff{recollement} $, then the operators
$ P_0 = e^{(n-1)r/2} \Delta^{G_0} e^{-(n-1)r/2} - c_n^2 $ and
$ P_1 = U^{1/2} e^{(n-1) r/2} \Delta^{G_1} e^{-(n-1)r/2} U^{-1/2} - c_n^2 $, 
with 
$U = \mbox{d}vol(G_1) / \mbox{d}vol(G_0) $,
satisfy the assumptions of the theorem.

\smallskip

The leading term of $ \reff{heatexpansion} $ involves the 
regularized volume term $ a_0 $. Other kinds of regularization have been
considered for similar purposes, for instance the $0-$volume used 
in \cite{GuZw2} or those used in \cite{Chri0,GrZw1}.

If we knew that $\xi_{\bf q}$ was a monotone function, this result
would yield immediately an equivalent for $ \xi_{\bf q} (\lambda) $  
 as $\lambda \uparrow
\infty $ by Karamata's Tauberian theorem. Unfortunately
 we don't know that it is a function neither that it is monotone (it is not
 the case in general). Nevertheless, we have the following result.
\begin{theo} \label{theo2}
i) $\xi_{\bf q}$ is a continuous function on $(0,\infty)$. \newline
ii) Assume that $G$ is of the form $ \reff{definitionmetrique} $ with $g=g_1$.
Assume furthermore that $ G $ is non trapping (see below) and that the principal symbols
of $P_0$ and $P_1$ coincide outside the region $\{ r > r_0 \}$, for some
$ r_0 $ large enough, then we have the complete asymptotic expansion
\begin{eqnarray}
 \xi_{\bf q} (\lambda) \sim \lambda^{n/2} \sum_{k \geq 0} b_k \lambda^{-k},
 \qquad \lambda \uparrow + \infty . 
 \label{expansionphase}
\end{eqnarray}
The coefficients $ b_k $ can be deduced from $ \reff{heatexpansion} $  
and in particular $ b_0 = (2 \pi)^{-n} \omega_n 
\int_X [\emph{d}vol_{\e}]_{\bf q} $.
\end{theo}

Note that the condition on the principal symbol means that the terms of order $
2$ of $ P_1 - P_0 $ are supported near infinity. For instance $P_1$ can be a
perturbation of $P_0$ with $P_0 = - \partial_r^2 - e^{-2r} \Delta^{g_0(0)} $
near infinity (the `product case') or with $P_0$ associated to a metric with
constant curvature near infinity. The latter can be of special interest in 
view of the recent results of \cite{CuVo1}.

We recall that $G$ is a {\it non trapping metric}, 
if for any compact subset $K $ of  $ T^* X \setminus 0 $ 
there exists $T_K \geq 0$ such that $\phi^t (K) \cap K  $ is empty for all
$|t| \geq T_K$, if $\phi^t$ is the geodesic flow on $ T^*X$. The results
already obtained on $ \Ra^n $ let us hope that the non trapping condition
in theorem $ \refe{theo2} $ can be relaxed to get a Weyl formula, i.e.
an equivalent for $ \xi_{\bf q}$. However this is an open question.

We shall use methods of semi-classical analysis to prove these theorems
and we will consider
\begin{eqnarray}
 H_{\e} = h^2 P_{\e} = H_{0} + \e V, \qquad h \in (0,1] .
 \label{operateurssemiclassiques}
\end{eqnarray}
We will use the notations $ E_{\e}(\mu) $ for the spectral projection
of $ H_{\e} $ on $ (-\infty,\mu) $ and
$$ R_{\e}(z) = (H_{\e}  - z)^{-1} , \qquad U_{\e}(t) = e^{- i t H_{\e} / h} . $$
In theorem $ \refe{theo1} $, we will choose $ h =t^{1/2} $ and in theorem
$ \refe{theo2} $, we shall set $  h = \lambda^{-1/2} $ and consider the
rescaled scattering phases $ \xi_{\bf q} (\mu , h) $ associated to $H_0,H_1$,
 near the energy $\mu = 1$,
since one has clearly $ \xi_{\bf q} (\mu h^{-2}) = \xi_{\bf q} ( \mu , h ) $. 

The non trapping condition will be used to show that
\begin{eqnarray}
 \int_{\Ra} \ \left| \left| \scal{r}^{-M} f (H_{\e}) U_{\e}(t) \scal{r}^{-M} 
 \right| \right|_{\infty} \mbox{d}t \ \leq \ C_{M,f,h_0,w}
\label{propagationestimate}
\end{eqnarray}
for some $M > 0$, all $f \in C^{\infty}_0$ supported close to $ 1 $ and uniformly w.r.t.
$h \in (0,h_0]$ and $\e \in [0,1]$. Such estimates are essentially
well known under the non trapping condition. They have been proved
in the Euclidean case by Robert-Tamura \cite{RoTa1} and simplified by
 G\'erard-Martinez \cite{GeMa0} (see also \cite{Robe1}), using the theory of 
 Mourre \cite{Mour1} with a {\it conjugate operator} 
 (see appendix $\refe{appendixuf}$)
defined as a suitable perturbation of generator of dilations $ r .
hD_r + h D_r.r $.
 As explained by Hislop and Froese in
\cite{FrHi0}, the generator of dilations does not fit the
hyperbolic framework and they built explicitly another conjugate
operator which makes Mourre theory applicable. In  appendix $
\refe{appendixuf} $, we sketch the proof leading to $
\reff{propagationestimate} $ by combination of the ideas of
\cite{GeMa0,Robe1} and \cite{FrHi0}, which is necessary since 
$ \reff{propagationestimate} $ is 
 only  known for fixed $h$ in the asymptotically hyperbolic case \cite{FrHi0,DBHS1}.

Note finally that if $ \reff{propagationestimate} $ can be improved so that
$  t^p \scal{r}^{-M} f (H_{\e}) U_{\e}(t) \scal{r}^{-M} $ has a polynomial bound 
in $h^{-1}$, then our method would show that 
$ \xi_{\bf q}^{(p)} $ has a complete expansion, obtained by differentiating
$ \reff{expansionphase} $.

\medskip

\noindent {\bf Acknowledgments:} I want to thank Peter Hislop for
helpful discussions as well as Gilles Carron and Didier Robert for
their interest and useful remarks.
\section{The basic tools of the proof}
\setcounter{equation}{0}
The purpose of this section, which is of pedagogic nature,
is to describe the main tools of the proof of theorem $ \refe{theo2} $.
The formulas that we are going to display hold for a much wider class of
operators than those defined on asymptotically hyperbolic manifolds and we want
 to separate the general ideas, sketched in this section,
from the specific analysis of the hyperbolic context given in sections
 $\refe{classical}$  and  $\refe{pseudo}$.

\subsection{Representation formula of $\xi_{\bf q}$} \label{representation}
Let us first assume that $ V $, defined by $ \reff{operateurssemiclassiques} $,
is compactly supported so that the Birman-Krein
spectral shift function $\xi_{1,\e} (\mu,h)$ is well defined for the
pair $H_0,H_{\e}$ (this is essentially standard but is anyway a consequence
of lemma $\refe{aprioriSchatten0}$).
 Then we can use the well known Birman-Solomyak formula 
\cite{BiSo1}
\begin{eqnarray}
\mbox{tr} \left( f (H_{\e}) - f (H_0) \right) = \int_0^1 \mbox{tr}
\left( f^{\prime} (H_{s \e} ) \e V \right) \ \mbox{d}s . 
\label{Birman-Solomyak}
\end{eqnarray}
The left hand side of $ \reff{Birman-Solomyak} $ is $ -
\scal{\xi_{1,\e},f^{\prime}} $ so, with the notations of
subsection $ \refe{notation} $, we get directly
\begin{eqnarray}
 \xi_{1,\e}(\mu,h) & = & - \int_0^1  \mbox{tr} \left( \frac{\partial
E_{s \e}}{\partial \mu} \e V \right) \ \mbox{d}s  =
- \int_0^{\e} \mbox{tr} \left( \frac{\partial
E_{s }}{\partial \mu}  V \right) \ \mbox{d}s ,
 \label{start}
\end{eqnarray}
since both sides of $ \reff{start} $ have the same derivative (in
the distributions sense w.r.t. $\mu$) and vanish for $\mu \ll 0 $.
 For the general case, i.e. $V$ non compactly supported, we consider the family
\begin{eqnarray}
 V^{(\kappa)} =  \theta (\kappa r)  V  \theta(\kappa r) , \qquad  \kappa \in ] 0 ,1 ] \label{perturbationlimite}
\end{eqnarray}
  with
$\theta \in C_0^{\infty}$, $ \theta = 1$ near $0$, and thus we can define
 the associated family of spectral shift
functions $\xi_{1,\e}^{(\kappa)}$. We obtain the following
representation formula for $\xi_{\bf q}$ :
\begin{lemm} \label{lemmeformexi} In the distributions sense on $ \Ra_{\mu}$, we have
\begin{eqnarray}
 \xi_{\bf q} (\mu,h) =
 -  \lim_{\kappa \downarrow 0} \left\{  \emph{tr}
\left(  \frac{\partial E_{  \e}^{(\kappa)}}{\partial \mu}
V^{(\kappa)} \right) \right\}_{\bf q}
\label{formulelimite}
\end{eqnarray}
Here $E_{\e}^{(\kappa)}$ is the spectral resolution associated
with $ H_{ \e }^{(\kappa)} = H_0 + \e V^{(\kappa)} $.
\end{lemm}
The proof of this lemma (given in section $\refe{pseudo}$)
is not very hard and formally obvious from $
\reff{start} $ and the definition of $ \xi_{\bf q} $.  Formula $
\reff{formulelimite} $ leads obviously to the study of distributions of the
form
\begin{eqnarray}
 \mbox{tr} \left(   \frac{\partial E_{ \e}^{(\kappa)}}{\partial \mu} f
(H_{\e}^{(\kappa)})  V^{(\kappa)}   \right) \label{formezero}
\end{eqnarray}
 with $ f \in C_0^{\infty} $, $f \equiv 1$ close to $1$. 
 Note that such a function $f$ can be added for free since we only consider
  $\mu $ close to $ 1 $.   
In order to simplify our notations, we drop the index $\kappa$ in
the sequel but the reader must keep in mind that we work with a
perturbation of the form $ \reff{perturbationlimite} $.

As usual,
we consider the (semi-classical) Fourier transform of
$\reff{formezero} $ which is $ \mbox{tr} \left( U_{\e}(t) f
(H_{\e}) V  \right) $. It is thus natural to consider
distributions of the more general form
\begin{eqnarray}
  \mbox{tr} \left( U_{\e}(t) f
(H_{\e}) K_{\e}  \right) \label{densite}
\end{eqnarray}
  with $K_{\e}$ trace class, for instance
$K_{\e} = \tilde{f}(H_{\e}) V$ with $\tilde{f} \in C_0^{\infty}$,
$\tilde{f}f=f$. Recall that we consider perturbations of the form
$ \reff{perturbationlimite} $, so that $\tilde{f(H_{\e}^{(\kappa)})V^{(\kappa)}} $ is trace class for
each $\kappa > 0$, but for $\kappa = 0$ we only have $K_{\e} \in
\Schn{q}$ in which case the trace makes sense only once
$\{.\}_{\bf q}$ has been taken. 

We will have  to consider the $\e$ derivatives
of $ f (H_{\e}) U_{\e}(t) $, that is why we quote quote the formula
\begin{eqnarray}
  \partial_{\e}
U_{\e} (t) = \frac{i}{h} \int_{0}^{t} U_{\e}(t-s) V U_{\e}(s) \ \mbox{d}s
\label{deriveeprop}
\end{eqnarray}
which holds, for instance, in the strong sense on the domain of the
operators $H_{\e}$ (which is independent of $ \e $, see proposition 
$ \refe{domaine}$). Since we want to use $\reff{propagationestimate}$, it is
important to keep a spectral cutoff in front of each $U_{\e}(t)$. To that end,
we use a very simple trick. Let us
introduce the notations
\begin{eqnarray}
 U^{f}_{\e}(t) = f (H_{\e}) U_{\e}(t), \qquad
S_{\e}^f = f (H_{\e}) . \label{notationpropagation}
\end{eqnarray}
 Then for any small neighborhood $\tilde{I}$
of $\mbox{supp} \ f$,
we can choose $\tilde{f}$ smooth, supported in $\tilde{
I}$, such that $\tilde{f}f=f$ and we obtain
$$
 \partial_{\e}  U_{\e}^f (t)  =  \partial_{\e} \left( S_{\e}^{f}
U_{\e} (t) S_{\e}^{\tilde{f}} \right)
 =  \left( \partial_{\e}  S_{\e}^{f} \right)
U_{\e} (t) S_{\e}^{\tilde{f}} +  S_{\e}^{f} \left( \partial_{\e}
U_{\e} (t) \right) S_{\e}^{\tilde{f}}  + S_{\e}^{f} U_{\e} (t)
\big( \partial_{\e}  S_{\e}^{\tilde{f}} \big) .
$$
and using $ \reff{deriveeprop} $, we get the  formula
\begin{eqnarray}
 \partial_{\e}  U_{\e}^f (t) =   \partial_{\e}  S_{\e}^{f}
\ U_{\e}^{\tilde{f}}(t) + \frac{i}{h}
 \int_{0}^{t} U^{f}_{\e}(t-s) V U^{\tilde{f}}_{\e}(s) \ \mbox{d}s +
 U_{\e}^{f}(t) \partial_{\e} S_{\e}^{\tilde{f}} .
\end{eqnarray}
This can be obviously iterated
and proves the following more general result
\begin{lemm} \label{Leibnitz0} For all  $k \geq 1$,
$ \partial_{\e}^k U_{\e}^f (t)  $ is a linear combination
with universal coefficients of
\begin{eqnarray}
 h^{-j}   \int_{ F_t^j }
\partial_{\e}^{l_0} S_{\e}^{\tilde{f}_0 } U_{\e}^{f_0}(t_0) 
  \partial_{\e}^{l_1} S_{\e}^{\tilde{f}_1 } V
\cdots
  V  \partial_{\e}^{l_{j-1}} S_{\e}^{\tilde{f}_{j-1} } 
U_{\e}^{f_{j}} (t_j)
\partial_{\e}^{l_{j}} S_{\e}^{\tilde{f}_{ j  } } \
    \emph{d}^{j}  (t_0,\cdots,t_j)  , \label{integrale} \\
    \mbox{with}  \qquad  1 \leq j \leq k, \ \
    l_0 + \cdots + l_{j+1} = k - j, \ \ f_0 , \cdots ,
f_{j},\tilde{f}_0, \cdots,\tilde{f}_{j+1} \in C_{0}^{\infty}(\tilde{I})
                       \nonumber
\end{eqnarray}
and of
$ \partial_{\e}^l  S_{\e}^{\tilde{f} }
U_{\e}^{f} (t) \partial_{\e}^{k-l} S_{\e}^{\tilde{f} } $, for
$ 0 \leq l \leq k $. In $ \reff{integrale} $, we have used the notations
$$ F_t^j = \{(t_0, \cdots , t_j ) \in [0,t]^{j+1} \ |  \  t_0 + \cdots + t_j = t \} $$
and $ \emph{d}^{j}(t_0, \cdots ,t_j)$ for  the
($j$-dimensional) Lebesgue measure on the hyperplane $t_0 + \cdots + t_j = t$.
\end{lemm}

In the applications, we will get estimates on $ \reff{integrale} $
using $\reff{propagationestimate}$ combined with
the following easy estimate
\begin{eqnarray}
\int_{\Ra} \left| \int_{F_t^j} \psi_0 (t_0) \psi_1 (t_1) \cdots \psi_j (t_j) \
\mbox{d}^j (t_0, \cdots , t_j) \right| \ \mbox{dt} \leq \int_{\Ra} \left| \psi_0 \right|
 * \left| \psi_1 \right| *
\cdots * \left| \psi_j \right|(t) \ \mbox{d}t \label{Convolution}
\end{eqnarray}
valid for all integrable functions $\psi_0 , \cdots, \psi_j$. 
\subsection{Two microlocal tools}
The operators $K_{\e}$ in $\reff{densite}$ will essentially be
pseudo-differential
operators. Sometimes we will need to shift the support of their
symbols by the classical Hamilton flow $\phi_{\e}^t$. To that end,
we will use the fact that for any $t_0$ we have
\begin{eqnarray}
  \mbox{tr} \left( U_{\e}(t) f (H_{\e}) K_{\e}  \right) =
\mbox{tr} \left( U_{\e}(t) f (H_{\e}) U_{\e}(t_0) K_{\e}
U_{\e}(-t_0)  \right) . \label{shifttrick}
\end{eqnarray}
 This remark was used by Robert in \cite{Robe2} and
follows trivially by centrality of the trace. In general, we
cannot obtain an explicit formula for $ U_{\e}(s) K_{\e} U_{\e}( -
s) $ and rather get an approximation. This implies that we have to
study an error term and this is why we display the following
explicit formulas. If $K_{\e}^s$ is $C^1$ w.r.t. $s$, satisfying
$K_{\e}^0 = K_{\e}$, we have
\begin{eqnarray}
U_{\e}(t_0) K_{\e} U_{\e}(-t_0) = K_{\e}^{t_0} + \frac{i}{h}
\int_0^{t_0} U_{\e}(t_0 -s) \left( h i \frac{ \partial  }{
\partial s } K_{\e}^s - \left[ H_{\e} , K_{\e}^s \right] \right)
U_{\e} (s - t_0) \ \mbox{d} s  \label{Egorov}          
\end{eqnarray}
and this leads to the following exact formula
\begin{eqnarray}
 \mbox{tr} \left( U_{\e}(t) f
(H_{\e}) K_{\e}  \right) =  \mbox{tr} \left( U_{\e}(t) f (H_{\e})
 K_{\e}^{t_0}  \right) +  \mbox{tr} \left(
 U_{\e} (t) f (H_{\e} ) \ \frac{i}{h} \!  \int_0^{t_0}  h i
\partial_s  K_{\e}^s - \left[ H_{\e} , K_{\e}^s \right]  \
\mbox{d} s
 \right) \label{evolutiontrick}
\end{eqnarray}
since the terms $ U_{\e} (t_0 -s) $ and $ U_{\e} (s-t_0) $ cancel
out by centrality. Thus, if we are able to find such a $K_{\e}^s$
with $
 h i
\partial_s  K_{\e}^s - \left[ H_{\e} , K_{\e}^s \right]
 $ small in a certain sense, we see that the study of $ \mbox{tr} \left( U_{\e}(t) f
(H_{\e}) K_{\e}  \right) $ reduces to the one of $  \mbox{tr}
\left( U_{\e}(t) f (H_{\e})
 K_{\e}^{t_0}  \right) $, up to a remainder which is given
 explicitly by $ \reff{evolutiontrick} $. The method leading to
 the calculation of $ K_{\e}^s $ is the usual one given by
the Egorov theorem. 
We refer to \cite{Robe1} for a proof of this theorem. The main point is that
the symbol (in each chart) 
of $ K_{\e}^s $ has an explicit expression in term of the symbol
of $ K_{\e} $ and of the Hamiltonian flow of $ p_{\e} $. 

If $K_{\e}$ is a pseudo-differential operator with a symbol
supported in a suitable region of $T^*X$ (and this can be achieved
by replacing $ K_{\e} $ by $K_{\e}^{t_0}$ thanks to the above
trick), we shall see that it can be factorized as
\begin{eqnarray}
  K_{\e} = A_{\e} B_{\e}^*, \qquad  \mbox{with} \ \  A_{\e},
B_{\e} : L^2(\Ra^n) \rightarrow L^2 (X)  .
\label{basefactorisation}
\end{eqnarray}
In practice such a factorization will  only be obtained
approximately, i.e. $ K_{\e} - A_{\e} B_{\e}^* $ {\it negligible}
(in a sense defined rigorously in section $ \refe{pseudo} $).
 The operator $A_{\e}$
will be used to {\it intertwin} $U_{\e}(t)$ with a {\it free
dynamic} given by $U (t) = \exp (- i t P / h)$, i.e. make $
U_{\e}(t) A_{\e} - A_{\e} U (t) $ small, in some sense, with $P = p
(hD)$  differential operator with constant coefficients
on $\Ra^n$. The explicit formula for $ U_{\e}(t) A_{\e} - A_{\e} U
(t)  $ is easily seen to be
\begin{eqnarray}
  U_{\e}(t) A_{\e} - A_{\e} U (t) = \frac{1}{i h} \int_0^t
U_{\e}(t-s) \left( H_{\e} A_{\e} - A_{\e} P \right) U (s) \
\mbox{d}s  \label{entrelacement}
\end{eqnarray}
and this shows that, if $ \reff{basefactorisation} $ holds, then we
have
\begin{eqnarray}
 \mbox{tr} \left( U_{\e}(t) f
(H_{\e}) K_{\e}  \right) -  \mbox{tr} \left( f (H_{\e}) A_{\e} U
(t) B_{\e}^* \right)  =  \qquad \qquad \qquad \qquad \qquad \nonumber \\
\frac{1}{i h} \int_0^t \mbox{tr} \left( f (H_{\e}) U_{\e}(t-s)
\left( H_{\e} A_{\e} - A_{\e} P \right) U (s) B_{\e}^* \right) \
\mbox{d}s . \label{Isozakitrick}
\end{eqnarray}
Thus if we are able to show that  $  \left( H_{\e} A_{\e} -
A_{\e} P \right) U (s) B_{\e}^*   $ is small (see lemma
 $ \refe{clefnegligeable}$), we see that we are left with the study of 
 $ \mbox{tr} \left( f (H_{\e}) A_{\e} U (t) B_{\e}^* \right) $. This
trace is easy to study since we shall have explicit expressions
for the operators $A_{\e},B_{\e}^*$ and $U(t)$. The construction
of the operators $A_{\e}$ and $B_{\e}$ will follow the scheme of
Isozaki-Kitada's method explained in the next subsection.

\subsection{The method of Isozaki-Kitada} \label{IsozakiKitadaformel}
In this part, we recall the principle of the construction of the
operators $A_{\e},B_{\e}$ by the method of Isozaki-Kitada introduced in 
\cite{IsKi1}. This
method has only been used on $\Ra^n$ for Euclidean scattering 
\cite{GeMa1,Robe2,Robe4,DeGe,Bouc2,Bouc4} 
but it turns out that it can be used in 
our framework as well, with
some changes on which we shall put emphasize in section $
\refe{main} $.

We first recall the algebraic formulas which enters into the game
for a general differential operator of order $2$ (on a general
manifold $X$), which we still denote $H_{\e} = h^2 P_{\e}(x,D_x)$ with
$$ P_{\e}(x,D_x) = p_{\e}(x,D_x) + p_{\e}^{(1)}(x,D_x) + p_{\e}^{(2)}(x) $$  in
coordinates $x = (x_1,\cdots, x_n)$. Here $p_{\e}(x,\xi)$ is the
principal symbol, i.e. homogeneous of degree $2$ w.r.t. to the
dual coordinates $\xi = (\xi_1,\cdots,\xi_1)$, and $p_{\e}^{(j)}$
are homogeneous  of degree $2-j$ for $j=1,2$.
We look for $A_{\e},B_{\e}$ defined as operators of the form $J (
\varphi_{\e} , a_{\e})$ and $J ( \varphi_{\e}, b_{\e}) $ where
\begin{eqnarray}
 J(\varphi_{\e},a_{\e}) u (x) = (2 \pi h)^{-n} \int \! \! \int
e^{ \frac{i}{h} ( \varphi_{\e}(x,\xi) - x^{\prime}.\xi )}
a_{\e}(x,\xi,h) u (x^{\prime}) \ \mbox{d} x^{\prime} \ \mbox{d}
\xi . \label{Fourierintegral}
\end{eqnarray}
Note that $ \reff{Fourierintegral} $ defines actually an operator
from $L^2(\Ra^n)$ into itself (under suitable conditions on
$\varphi_{\e}$ and $a_{\e}$); however, in the applications,
$a_{\e}$ (and $b_{\e}$) will be supported into a region of
$\Ra^n_x \times \Ra^n_{\xi}$ whose projection onto $\Ra_x^n$ is
included into a coordinate chart of $X$. Thus, up to an invertible operator, 
we can consider $ \reff{Fourierintegral} $ as an
operator from $L^2(\Ra^n)$ to $L^2(X)$.
\newline
Following $ \reff{entrelacement} $, we see that we have to study $
H_{\e} J (\varphi_{\e},a_{\e}) - J (\varphi_{\e},a_{\e}) P $.
Since $P=p(hD)$ is a Fourier multiplier, we have obviously $ J
(\varphi_{\e},a_{\e})P  = J (\varphi_{\e},a_{\e} p )   $. On the
other hand, we see easily that
\begin{eqnarray}
H_{\e} \left( e^{\frac{i}{h} \varphi_{\e}} a_{\e} \right) =
e^{\frac{i}{h} \varphi_{\e}} \left(  p_{\e} (x, \partial_x
\varphi_{\e}) a_{\e} + h i^{-1} L_{\e}(x,\partial_x) a_{\e} + h^2
P_{\e}(x,D_{x}) a_{\e} \right) .
\end{eqnarray}
Here, $L_{\e}(x,\partial_x)$ is  a differential operator of order
$1$ defined as
$$ L_{\e} (x,\partial_x) = w_{\e} (x,\xi). \partial_x + c_{\e}(x,\xi) $$
where we have set
$$ w_{\e} (x,\xi) = (\partial_{\xi} p_{\e})(x,\partial_x \varphi_{\e})  , \qquad c_{\e}(x,\xi) =
p_{\e}(x,\partial_x) \varphi_{\e} + i p_{\e}^{(1)}(x,\partial_x
\varphi_{\e}) .
$$
All this shows that if we look for $a_{\e} = a_{\e}^{(0)} + h
a_{\e}^{(1)} + \cdots + h^N a_{\e}^{(N)}$, then
\begin{eqnarray}
 H_{\e} J
(\varphi_{\e},a_{\e}) - J (\varphi_{\e},a_{\e}) P = J
(\varphi_{\e}, \widetilde{a}_{\e} ) \label{entrelacagepetit}
\end{eqnarray}
 with $ \widetilde{a}_{\e} =
\widetilde{a}_{\e}^{(0)} + h \widetilde{a}_{\e}^{(1)} + \cdots +
h^{(N+2)} \widetilde{a}_{\e}^{(N+2)} $. For $ 0 \leq j \leq N + 2
$, the functions $ \widetilde{a}_{\e}^{(j)} $ are given by
\begin{eqnarray}
\widetilde{a}_{\e}^{(j)} = \left( p_{\e} (x, \partial_x
\varphi_{\e}) - p (\xi) \right) a_{\e}^{(j)} - i
L_{\e}(x,\partial_x) a_{\e}^{(j-1)} + P_{\e} (x,D_x)
a_{\e}^{(j-2)},
\end{eqnarray}
where we use the convention that $ a_{\e}^{(k)} = 0 $ for $k < 0$
or $k > N$. 

Since we want to make $ \reff{entrelacagepetit} $ small, we look for 
$ \varphi_{\e} (x,\xi) $ such that
\begin{eqnarray}
 p_{\e} (x, \partial_x
\varphi_{\e}) = p (\xi) \label{HamiltonJacobi}
\end{eqnarray}
which is usually called the {\it Hamilton-Jacobi equation}. We
also need to find $a_{\e}^{(k)}$ such that
\begin{eqnarray}
  L_{\e}(x,\partial_x) a_{\e}^{(0)}& = & 0 ,
 \label{transport0}  \\
   L_{\e}(x,\partial_x) a_{\e}^{(k)} & = & - i P_{\e}(x,D_x) a_{\e}^{(k-1)}  , \qquad k \geq 1 ,  \label{transportj}
\end{eqnarray}
which are the {\it transport equations}. The resolution of $
\reff{HamiltonJacobi} $, $ \reff{transport0} $ and $
\reff{transportj} $ rely upon estimates on classical trajectories.
The technical part leading to such estimates in the asymptotically
hyperbolic case is the purpose of section $ \refe{classical} $.
Here we recall the general method. Assume that we can find  $
S_{\e} (t,x,\xi) $ defined on $ [0,\infty) \times \Gamma $, for
some open set $\Gamma \subset \Ra^{2n}$, such that
$$ \partial_t S_{\e} = p_{\e}(x, \partial_x S_{\e}) , \qquad S_{\e}(0,x,\xi) =
x. \xi . $$
The existence of $S_{\e}$ will follow from suitable estimates on
$ \phi_{\e}^t $, the
Hamilton flow of $p_{\e}$. In practice, $\Gamma$ is such that
\begin{eqnarray}
 \partial_{\xi} S_{\e}(t,x,\xi) \rightarrow \infty,
\qquad \mbox{as} \ \ t \uparrow \infty . \label{unbounded}
\end{eqnarray}
Now, using the fact that $S_{\e}$ is a {\it generating function} of the flow, i.e.
\begin{eqnarray}
 \phi_{\e}^t (x, \partial_x S_{\e}) = (\partial_{\xi}S_{\e},\xi) ,
 \label{generatingflow}
\end{eqnarray}
 the invariance of $p_{\e}$ by
the flow and $ \reff{unbounded} $ imply that
\begin{eqnarray}
  \lim_{t \uparrow \infty} p_{\e}(x, \partial_x S_{\e}) = p (\xi)
  \label{HJapprox}
\end{eqnarray}
 provided
$$ p (\xi) = \lim_{x \rightarrow \infty} p_{\e}(x,\xi) . $$
This shows that we have to built $ \varphi_{\e} $  such that
$\partial_{x}
\varphi_{\e} = \lim_{t \uparrow \infty} \partial_x S_{\e}$, or equivalently
such that
\begin{eqnarray}
 \partial_x \varphi_{\e}(x,\xi )
& = & \xi + \int_0^{\infty} \partial_t \partial_x S_{\e}(t,x,\xi)
 \ \mbox{d}t . \nonumber
\end{eqnarray}
Of course, if such a function $ \varphi_{\e} $
exists it is not unique. A
possible construction is the following
\begin{eqnarray}
 \varphi_{\e}(x,\xi )
& = & x.\xi + \int_0^{\infty} \partial_t \left( S_{\e}(t,x,\xi) -
\widetilde{S}_{\e} (t,\xi) \right)
 \ \mbox{d}t \label{definitionphase}
\end{eqnarray}
with $\widetilde{S}_{\e}(t,\xi)$ independent of $x$ and such that the integral
of the right hand side converges. In practice, $\widetilde{S}_{\e}$ will be easy
to find. In that case, by construction, $
\reff{HJapprox} $ holds and we get $ \reff{HamiltonJacobi} $.

Note that, for Euclidean scattering on $\Ra^n$, we have
$p(\xi)=|\xi|^2$. In the asymptotically hyperbolic situation, with
the coordinates $x=(r,y)$ and $\xi=(\rho,\eta)$ we will consider
$p(\xi)= \rho^2$.

For the resolution of the transport equations $ \reff{transport0}
$, $ \reff{transportj} $ we study the solution $\check{x}_{\e}^t
(x,\xi)$ of
$$ \frac{d \check{x}_{\e}^t }{d t}  = w_{\e}(\check{x}_{\e}^t,\xi), \qquad \check{x}_{\e}^0 (x,\xi) = x . $$
Indeed, any solution $a(x,\xi)$ of $ L_{\e}(x, \partial_x) a = 0 $
satisfies
$$ a (\check{x}_{\e}^t,\xi) = a (x,\xi) \exp \left( - \int_0^t c_{\e}(\check{x}_{\e}^s,\xi) \ \mbox{d}s \right) . $$
In the applications, we shall prove that $ \check{x}_{\e}^t $ is defined
on $[0,\infty) \times \Gamma$ and satisfies (in a sense to be made
precise)
\begin{eqnarray}
 \check{x}_{\e}^t \sim x + t \partial_{\xi} p (\xi), \qquad t \uparrow
\infty . \label{spacetime}
\end{eqnarray}
 This is actually well known in the Euclidean case (see for instance
\cite{IsKi1,GeMa1,Robe4,DeGe}).
Hence, if we look for $a_{\e}^{(0)}$ such that $  a_{\e}^{(0)} (x,\xi)
\rightarrow 1 $ as $ (x,\xi) \rightarrow \infty $ in $ \Gamma $,
 we obtain
\begin{eqnarray}
  a_{\e}^{(0)} (x,\xi) = \exp \left(  \int_0^{\infty} c_{\e}(\check{x}_{\e}^t,\xi) \
\mbox{d} t \right) . \label{abrzerob}
\end{eqnarray}
 Similarly, by the
method of variation of constants, we see that any solution of $
L_{\e} a = \iota $ satisfies
$$ a (\check{x}_{\e}^t,\xi) = \left( a ( x,\xi ) + \int_0^t \iota
(\check{x}_{\e}^s,\xi) \exp \left( \int_0^s
c_{\e}(\check{x}_{\e}^u,\xi) \ \mbox{d} u \right) \ \mbox{d} s
\right) \exp \left( - \int_0^t c_{\e}(\check{x}_{\e}^s ,\xi) \
\mbox{d}s \right) .
$$
Thus if we look for a solution $a_{\e}^{(k)}$ of $
\reff{transportj} $ such that $a_{\e}^{(k)}(x,\xi) \rightarrow 0$ as 
$ (x,\xi) \rightarrow  \infty $ in $ \Gamma $, we get
\begin{eqnarray}
a_{\e}^{(k)} (x,\xi) = - \int_0^{\infty}  i P_{\e}(x,D_x)
a_{\e}^{(k-1)} (\check{x}_{\e}^t,\xi) \ \exp \left( \int_0^t
c_{\e}(\check{x}_{\e}^s,\xi) \ \mbox{d} s \right) \ \mbox{d} t .
\label{abrjb}
\end{eqnarray}
 Here again, the
convergence of the integrals is justified by the appropriate
estimates on the functions $c_{\e}$ and  $P_{\e}(x,D_x)
a^{(k-1)}_{\e}$ which need to be shown. Such estimates are well known in the
Euclidean case and will follow from section $ \refe{classical} $
for the hyperbolic one.

The formulas $ \reff{abrzerob} $ and $ \reff{abrjb} $ define
$a_{\e}^{(k)}$, for $k \geq 0$, in $ \Gamma $. In section $ \refe{main} $, we
will explain how to define them globally, i.e. how to cut them
off outside a suitable area.

 Recall that we want to consider a factorization of the form $
J (\varphi_{\e},a_{\e}) J (\varphi_{\e},b_{\e})^* $ whose Schwartz
kernel is
$$  {\mathcal K}_{\e} (x,x^{\prime}) = (2 \pi h)^{-n} \int e^{\frac{i}{h} ( \varphi_{\e}(x,\xi) -
\varphi_{\e} (x^{\prime},\xi) ) }  a_{\e}(x,\xi,h)
\overline{b_{\e}(x^{\prime},\xi,h)} \ \mbox{d} \xi . $$ This
kernel is the one of a pseudo-differential operator since {\it
Kuranishi's trick}, namely
$$ \varphi_{\e}(x,\xi) - \varphi_{\e} (x^{\prime},\xi) = (x-x^{\prime}). \theta_{\e} (x,x^{\prime},\xi)  $$
which is of course obtained by Taylor's formula, allows to write
$$  {\mathcal K}_{\e} (x,x^{\prime}) = (2 \pi h)^{-n} \int e^{\frac{i}{h} (x - x^{\prime} ).\theta }  \vartheta_{\e}
(x,x^{\prime},\theta,h)  \ \mbox{d} \theta
$$
provided the following map is a diffeomorphism for each
$x,x^{\prime}$ belonging to the projections of the supports of $ a_{\e}$ and $
b_{\e}$:
$$ \xi \mapsto \theta_{\e} (x,x^{\prime},\xi) = \int_0^1 \partial_x \varphi_{\e} (x^{\prime}+ t (x-x^{\prime}),\xi)
\ \mbox{d}t . $$ Note that, in view of $ \reff{definitionphase} $,
we see that if $\partial_x \partial_t S_{\e}(t,x,\xi)$ is small
(which will indeed be the case), then the map $\xi \mapsto
\theta_{\e}$ is close to the identity and easily seen to be a
diffeomorphism. In this case we have
$$ \vartheta_{\e}
(x,x^{\prime},\theta,h) =
a_{\e}(x,\theta_{\e}^{-1}(x,x^{\prime},\theta),h)
\overline{b_{\e}(x^{\prime},\theta_{\e}^{-1}(x,x^{\prime},\theta),h)}
 \left| \mbox{det} \ \partial_{\theta} \theta_{\e}^{-1} (x,x^{\prime},\theta)
 \right| .
$$
Using a general elementary property of pseudo-differential operators,
we have
$$ {\mathcal K}_{\e} (x,x^{\prime}) \sim (2 \pi h)^{-n} \sum_{j \geq 0} h^j \int e^{\frac{i}{h}(x-x^{\prime}).\theta}
\sum_{|\alpha| = j}  \frac{1}{\alpha !
}\partial_{x^{\prime}}^{\alpha} D_{\theta}^{\alpha} \vartheta_{\e}
(x,x^{\prime},\theta,h)_{|x^{\prime}=x}
 \mbox{d} \theta . $$
By identification of the powers of $h$, this allows to find $ b_{\e}  = b_{\e}^{(0)} + \cdots + h^N b_{\e}^{(N)}$
 such that, modulo $h^N$, the right hand side is the expansion of
 the Schwartz kernel of $K_{\e}$, with the notation of $\reff{basefactorisation}$.
 For instance, if the principal symbol of $K_{\e}$ is $ \sigma (x,\theta)$,
  we must have
 $$ a_{\e}^{(0)}(x,\theta_{\e}^{-1}(x,x,\theta))
\overline{b_{\e}^{(0)}(x,\theta_{\e}^{-1}(x,x,\theta))}
 \left| \mbox{det} \ \partial_{\theta} \theta_{\e}^{-1} (x,x,\theta)
 \right| = \sigma (x,\theta)  $$
 which implies that
\begin{eqnarray}
 \overline{ b_{\e}^{(0)} (x,\xi) } = \sigma (x, \theta_{\e}(x,x,\xi))  a^{(0)}(x,\xi)^{-1}
 \left| \mbox{det} \ \partial_{\xi} \theta_{\e} (x,x,\xi)
 \right| . \label{formebzero}
\end{eqnarray}
 More generally, one can get explicit expressions for
 $b_{\e}^{(1)},\cdots,b_{\e}^{(N)}$ and  the important remark is
 that they are linear combinations of products of (derivatives of)
 $a_{\e}^{(0)}(x,\xi)^{-1}$,
 $ a_{\e}^{(k)}(x,\xi) $ for $k \geq 1$,
 $  \theta_{\e} (x,x^{\prime},\xi) $ (evaluated at $x^{\prime}=x$)
 and the symbols of $K_{\e}$ evaluated at $ \left( x ,
 \theta_{\e}(x,x,\xi) \right) $.
In practice, $ \theta_{\e}(x,x,\xi) -\xi $ will be small in the
region that we will consider and $ \Gamma $ will be a neighborhood of 
$ \mbox{supp} \sigma $, so that $ b_{\e}^{(0)}, b_{\e}^{(1)}, \cdots $ will only
depend on  $ \sigma $ and on the values of $ a_{\e}^{(0)}, a_{\e}^{(1)}, \cdots
$ and $ \varphi_{\e}$ in $ \Gamma $. In particular $ \reff{abrzerob}$ will
ensure that $ a_{\e}^{(0)} $ doesn't vanish.

\section{Some estimates on the geodesics} \label{classical}
\setcounter{equation}{0}
 In this technical section, we prove long
time estimates on classical trajectories, in suitable areas of
$T^* X$, needed to justify Isozaki-Kitada's method in the
asymptotically hyperbolic case.
\subsection{Results of this section}
We consider the function defined for $r\in \Ra  $, $ \rho \in \Ra
$ and $y,\eta \in \Ra^{n-1}$ by
$$ p_{\e} (r,y,\rho,\eta) = \rho^2 + e^{-2r} g_{\e}(r,y,\eta) $$
where the function $g_{\e} (r,y,\eta)$ defines a metric on
$\Ra^{n-1}$, i.e. $g_{\e}$ is a smooth function which is an
homogeneous polynomial of degree $2$ w.r.t. $\eta$ and such that
for some $C_0 > 0$
\begin{eqnarray}
C_0^{-1} |\eta|^2 \leq g_{\e}(r,y,\eta) \leq C_0 |\eta|^2,
 \label{elliptic0}
\end{eqnarray}
for all $ r \in \Ra  $,  $ y,\eta \in \Ra^{n-1}$ and $\e \in
[0,1]$. We assume
 furthermore that $g_{\e}$ is of the following form
\begin{eqnarray}
 g_{\e} (r,y,\eta) = g (y,\eta) + \e e^{-r} \tilde{g} \left(
e^{-r} ,y,\eta \right) \label{tildemetrique}
\end{eqnarray}
 where $g$ and $\tilde{g}$ satisfy the
following estimates
\begin{eqnarray}
|\partial_{x}^{k}\partial_{y}^{\alpha} \partial_{\eta}^{\beta}
\tilde{g}(x,y,\eta)| + |\partial_{y}^{\alpha}
\partial_{\eta}^{\beta} g(y,\eta)| \leq
C_{k,\alpha,\beta} |\eta|^{2-|\beta|}, \qquad x \in [0,1] , \ y,
\eta \in \Ra^{n-1}. \label{bound}
\end{eqnarray}
Note that in particular, we have 
$ g_{\e}(r,y,\eta) - g (y,\eta) = {\mathcal O}(e^{-r})|\eta|^2  $.
 All these estimates are satisfied
by the principal symbol of $P_{\e}$ in any chart at infinity
$(R,\infty) \times \Omega$. Thus we can assume that this principal
symbol is the restriction of some function $p_{\e}$ satisfying the
above estimates. We are going to prove estimates on trajectories
of some vector fields  and these trajectories will turn out to lye
inside $(R,\infty)\times \Omega \times \Ra^n$. This means that the
results of the present section can be obviously considered as
local results in $T^*X$.

Our first goal is the study of the Hamiltonian flow $\phi_{\e}^t $
of the function $p_{\e}$, that is the solution of 
$$ \dot{\phi}^t_{\e} = {\mathbf H}_{\e} \left( \phi^t_{\e} \right) $$
where the notation $\dot{ }$ stands for $d/dt$ throughout the
section and ${\mathbf H}_{\e}={\mathbf H}_{\e}(r,y,\rho,\eta)$ is
the Hamiltonian vector field of $p_{\e}$, that is
\begin{eqnarray}
{\mathbf H}_{\e} = \left(  \begin{matrix}
\partial_{\rho} p_{\e} \\
\partial_{\eta} p_{\e} \\ - \partial_r p_{\e} \\ - \partial_y p_{\e}
\end{matrix} \right)
= \left(  \begin{matrix}
2 \rho \\
e^{-2r} \partial_{\eta} g_{\e}(r,y,\eta) \\  e^{-2r} \left( 2 g_{\e}(r,y,\eta) - \partial_r g_{\e}(r,y,\eta) \right)  \\
- e^{-2r} \partial_y g_{\e} (r,y,\eta)
\end{matrix} \right)
 .
\label{flow}
 \end{eqnarray}
 We denote the components of the flow by
$(r^t_{\e},y^t_{\e},\rho^t_{\e},\eta^t_{\e})$. They are functions
of the initial condition $(r,y,\rho,\eta)$. We remark that
$\phi_{\e}^t$ is defined for $t \in \Ra$ for if this was wrong,
then $|\phi_{\e}^t|$ should blow up in finite time. This cannot
happen since the conservation of energy, $p_{\e} = p_{\e} \circ
\phi_{\e}^t $, implies easily that $\dot{\phi}_{\e}^{t}$ is
bounded, hence $\phi_{\e}^t - \phi_{\e}^0 = {\mathcal O}(t)$ can
not blow up in finite time.

\smallskip

Let us introduce some notations. We consider  energy intervals of
width $2w$ defined by
$$ I (w) = (1-w,1+w) .$$
In the end, $w$ will be chosen small enough but for the time being
we only assume that $0<w<1/2$. We also introduce the outgoing
(resp. incoming) parameters $\sigma_{+} > 0$ and $\delta_{+} > 0$
(resp. $\sigma_{-} > 0$ and $\delta_{-} > 0$) defined by
$$ \sigma_{\pm} = \left(1 \mp w \right)^{1/2} - \delta_{\pm} . $$
Of course, this makes sense provided $0 < \delta_{\pm} < \left( 1
\mp w \right)^{1/2}$.
 The last two parameters we shall need are a positive real number
 $R > 0$ and an arbitrary open subset $ \Omega \subset \Ra^{n-1}$. We can
 now define the {\it outgoing area} $\Upsilon^{+}(R,\sigma_+,w,\Omega)$ and the {\it incoming area}
 $\Upsilon^{-}(R,\sigma_{-},w,\Omega)$ by
\begin{eqnarray}
 \Upsilon^{\pm} (R,\sigma_{\pm},w,\Omega):= \left\{ (r,y,\rho,\eta) \in \Ra^{2n} \ | \
p_0 (r,y,\rho,\eta) \in I (w), \ \pm \rho > - \sigma_{\pm} , \ r >
R, \ y \in \Omega \right\} . \label{presquesortante}
\end{eqnarray}
 The point of considering such
areas is that we have a splitting
\begin{eqnarray}
 \Upsilon^{+} (R,\sigma_{+},w,\Omega) \cup \Upsilon^{-}
(R,\sigma_{-},w,\Omega) = \left\{ (r,y,\rho,\eta) \ | \
 p_0 (r,y,\rho,\eta) \in I(w) , \ r > R , \ y \in \Omega   \right\}
 . \label{split}
\end{eqnarray}
The first result is the following.
\begin{prop} \label{aprioritraj} Let $0<w<1/2$,  $\sigma_{\pm} > 0$ and $\Omega$ as above.
There exists $R$ large enough and $C > 0$, depending only on $C_0$
and a finite number of constants $C_{k,\alpha,\beta}$ in
$\reff{bound}$, such that the following estimates hold
\begin{eqnarray}
r_{\e}^t & \geq & r  \pm t - C , \label{rt}  \\
\dot{\rho}_{\e}^t & \geq & 0  \label{rhot} \\
|y^t_{\e}-y| & \leq & C e^{-r} \label{yt}
\end{eqnarray}
for all $\e \in [0,1]$, $(r,y,\rho,\eta) \in
\Upsilon^{\pm}(R,\sigma_{\pm},w,\Omega)$ and $\pm t \geq 0$.
\end{prop}
The meaning of these statements is that the properties are true on
$\Upsilon^{+}$ (resp. $\Upsilon^{-}$) for $t \geq 0$ (resp. $t
\leq 0$). This result says that the geodesics starting in outgoing
(resp. incoming) areas stay in the neighborhood of infinity in $X$
for $t \geq 0$ (resp. $t \leq 0$).  More precise estimates on the
flow will be given in the theorem below. Note moreover that $
\reff{yt} $ shows that $y^t_{\e}$ lye in any arbitrary small
neighborhood of $\Omega$. It explains why our results can be
localized in charts at infinity on $X$.

\medskip

\noindent {\bf Remark ${\bf 1}$.} This proposition, or rather its
proof, shows in particular that the flow of $p_{\e}$ in $
\Upsilon^{\pm} (R,\sigma_{\pm},w,\Omega) $ depends only on the
values of $p_{\e}$ in the region $ \{ r \geq R - C \} $. This
implies that if we replace $\tilde{g}$ by $\tilde{g}^{(\kappa)} $
depending on a parameter $\kappa \in [0,1]$, in $
\reff{tildemetrique} $, such that the constants $C_0$ and
$C_{k,\alpha,\beta}$ can be chosen uniformly w.r.t. $\kappa$ for
$r$ large enough, then  the above proposition is true uniformly
w.r.t. $\kappa$.
 In
particular, this holds if one considers the principal symbol of $
\reff{perturbationlimite} $ which involves
\begin{eqnarray}
\tilde{g}^{(\kappa)}(e^{-r},y,\eta) = \theta (\kappa r)^2
\tilde{g} (e^{-r},y,\eta) \label{metriqueparametre} .
\end{eqnarray}

Using this remark, we now claim that  {\it all the results listed
in this subsection hold true uniformly w.r.t. $ \kappa \in [0,1] $
if one considers} $ \reff{metriqueparametre} $.

\smallskip

\begin{theo} \label{Precised} Let $w,\sigma_{\pm}$ and $\Omega$ be as in
proposition $\refe{aprioritraj}$. There exists $R$ large enough
such that for all $\gamma$ defined by $\partial^{\gamma} =
\partial_{r}^k
\partial_{y}^{\alpha} \partial_{\rho}^l \partial_{\eta}^{\beta}
\partial_{\e}^j$,
one can find $C_{\gamma} \geq 0$ satisfying
\begin{eqnarray}
\left| \partial^{\gamma} \left( r_{\e}^t - r - 2 t \rho \right)
 \right| & \leq &
C_{\gamma} \scal{t} e^{-(j + 1)r} |\eta|   \label{rt2} \\
\left| \partial^{\gamma} \left( y_{\e}^t - y \right) \right| &
\leq &
C_{\gamma} e^{-(j+1) r}  |\eta|   \label{yt2} \\
\left| \partial^{\gamma} \left(  \rho_{\e}^t - \rho \right)
\right| & \leq &
C_{\gamma} e^{-(j+1) r}  |\eta|  \label{rhot2} \\
\left| \partial^{\gamma} \left(  \eta_{\e}^t - \eta \right)
\right| & \leq & C_{\gamma} e^{-(j+1) r} |\eta|  \label{etat2}
\end{eqnarray}
for all $\e \in [0,1]$ provided the following condition holds
\begin{eqnarray}
\pm t \geq 0 \qquad \mbox{and} \qquad (r,y,\rho,\eta) \in
\Upsilon^{\pm} (R,\sigma_{\pm},w,\Omega) .
\end{eqnarray}
\end{theo}

\begin{coro}  \label{classicalscat} For $R$ large enough and for all
$(r,y,\rho,\eta) \in \Upsilon^{\pm}(R,\sigma_{\pm},w,\Omega)$ the
following limits exist for all $\e \in [0,1]$
$$ y_{\e}^{\pm}= \lim_{t \rightarrow \pm \infty} y_{\e}^t, \qquad
\eta_{\e}^{\pm}= \lim_{t \rightarrow \pm \infty} \eta_{\e}^t
\qquad \lim_{t \rightarrow \pm \infty} \rho_{\e}^t = \pm
p_{\e}^{1/2}
$$
where $p_{\e}=p_{\e}(r,y,\rho,\eta)$.
\end{coro}
This corollary follows very easily from the motion equations,
proposition $\refe{aprioritraj}$ and theorem $\refe{Precised}$
since, in particular, $\dot{y}^t_{\e}$, $\dot{\rho}^t_{\e}$ and
$\dot{\eta}^t_{\e}$ are ${\mathcal O}(e^{-2t})$.

\medskip

 Another important consequence of the theorem
is the following. Since $\nabla_{\rho,\eta} \left(
\rho^t_{\e},\eta^t_{\e} \right)$ is close to the identity matrix
if $e^{-r}\eta$ is small enough, we can hope that $(\rho,\eta)
\mapsto (\rho^t_{\e},\eta^t_{\e})$ is a diffeomorphism under
suitable conditions. This is why we introduce
$$ \Gamma^{\pm} (R,\varepsilon,w,\Omega) = \{ (r,y,\rho,\eta) \ | \
p_0 (r,y,\rho,\eta) \in I (w), \ r > R , \pm \rho > 0 ,  \ y \in
\Omega , \  e^{-2r} g (y,\eta) < \varepsilon  \}  . $$

Before studying the above diffeomorphism, let us note that any
area $\Gamma^{\pm}$ can be reached in finite time from
$\Upsilon^{\pm}$. Precisely we have
\begin{lemm} \label{Lemme} Let $0 < w   < w^{\prime} < 1/2$ and $\sigma_{\pm} > 0$.
 Assume that $\Omega$ is bounded and that $\Omega \Subset \Omega^{\prime}$.  
 Then, there exists $ R > 0 $ large enough such that for any
 $\varepsilon
> 0$ and $R^{\prime} > R $, there exists  $T > 0$ such that
$$  \phi_{\e}^t \left(  \Upsilon^{\pm}(R,\sigma_{\pm},w,\Omega) \right)
\subset \Gamma^{\pm} (R^{\prime},\varepsilon, w^{\prime},
\Omega^{\prime} ) , \qquad \forall \ \pm t \geq T
 , \ \ \e \in [0,1]. $$
\end{lemm}
This lemma follows easily from the previous results, the main tool
being the fact  that $e^{-r^t_{\e}} \eta^t_{\e} \rightarrow 0$ as
$\pm t \rightarrow + \infty$.
 Now let us  consider the  map
$$ \Phi_{\e}^t : (r,y,\rho,\eta) \mapsto (r,y,\rho^t_{\e},\eta^t_{\e}) . $$
Then we have the following result.
\begin{prop} \label{eiko1} Let $\Omega_0,\Omega_1$ be bounded, connected and such that $\Omega_1 \Subset \Omega_0$.
There exists $R_0 , R_1 > 0$ large enough and
$w_0,w_1,\varepsilon_0 , \varepsilon_1> 0$ small enough such that
for all $ \pm t \geq 0$
\begin{eqnarray}
  \Gamma^{ \pm } (R_1,\varepsilon_1,w_1,\Omega_1)
 \subset  \Phi_{\e}^t \left( \Gamma^{+} (R_0,\varepsilon_0,w_0,\Omega_0) \right)
 .\label{Part0}
\end{eqnarray}
Furthermore, $\Phi^t_{\e}$ is a diffeomorphism from $ \Gamma^{\pm}
(R_0,\varepsilon_0,w_0,\Omega_0) $ onto its range for all $\pm t
\geq 0$. If we denote the inverse map by $ (r,y,\rho,\eta) \mapsto
(r,y,\tilde{\rho}^t_{\e},\tilde{\eta}^t_{\e})$ we have
\begin{eqnarray}
\left| \partial^{\gamma} \! \left( \tilde{\rho}_{\e}^t - \rho
\right) \right|
 + \left| \partial^{\gamma} \! \left( \tilde{\eta}_{\e}^t - \eta \right) \right|
 & \leq &
C_{\gamma} e^{-(j+1)r}|\eta| ,  \qquad \mbox{if} \ \
\partial^{\gamma} =
\partial_r^k \partial_y^{\alpha} \partial_{\rho}^{l}
\partial_{\eta}^{\beta} \partial_{\e}^{j} \label{Diffe}
\end{eqnarray}
 with $C_{\gamma}$ independent of $\e \in [0,1]$, $\pm t \geq 0$
 and $(r,y,\rho,\eta) \in \Gamma^{\pm}(R_1,\varepsilon_1,w_1,\Omega_1)$.
\end{prop}
The main application of this proposition is the resolution of the
following eikonal equation
\begin{eqnarray}
 \partial_t S_{\e}^{\pm} = p_{\e} \left( r,y,\partial_{r}
S_{\e}^{\pm} , \partial_{y} S_{\e}^{\pm}
 \right) , \qquad  S_{\e}^{\pm}(0,r,y,\rho,\eta) = r \rho + y.\eta
 \label{eikonal}
\end{eqnarray}
 with $S_{\e}^{\pm} = S_{\e}^{\pm} (t,r,y,\rho,\eta,\e)$ defined
 on $\Gamma^{\pm} (R_1,\varepsilon_1,w_1,\Omega_1) $ for any
 $\pm t \geq 0$. We can solve $\reff{eikonal}$ since any
 $S$ which satisfy
\begin{eqnarray}
 S (t,r,y,\rho^t_{\e},\eta^t_{\e},\e) = r\rho + y.\eta + t
p_{\e} - \int_0^t \left(r \partial_r p_{\e} + y.\partial_y p_{\e}
\right) \circ \phi^s_{\e} \ \mbox{d}s \label{eikonalexplicit}
\end{eqnarray}
solves $\reff{eikonal}$. Thus the composition of the right hand
side of $\reff{eikonalexplicit}$ with the inverse of $\Phi^t_{\e}$
is a solution to $\reff{eikonal}$. We shall show the following
\begin{prop} \label{eiko2} We have a solution $S_{\e}^{\pm}$
of $\reff{eikonal}$ on
$\Gamma^{\pm}(R_1,\varepsilon_1,w_1,\Omega_1)$ for $\pm t \geq 0$.
 Furthermore, if $R_1$ is large enough and $\varepsilon_1$ small enough,
 we have
 \begin{eqnarray}
 \left| \partial^{\gamma} \left( S_{\e}^{\pm} -r\rho - y.\eta \mp t \rho^2 \right)
 \right| \leq C_{\gamma} e^{-(j+1)r}|\eta|, \qquad
 \partial^{\gamma} = \partial_r^k \partial_y^{\alpha} \partial_{\rho}^{l}
\partial_{\eta}^{\beta} \partial_{\e}^{j} \label{approxeiko}
 \end{eqnarray}
for $\e \in [0,1]$, $\pm t \geq 0$ and $(r,y,\rho,\eta) \in
\Gamma^{\pm}(R_1,\varepsilon_1,w_1,\Omega_1)$.
\end{prop}

\medskip

\noindent {\bf Remark ${\bf 2}$.} We add this short remark in
order to explain in what sense these constructions are continuous
w.r.t. $\kappa$, if we consider $ \reff{metriqueparametre} $. Let
us recall again that all the coming proofs in the next subsections
work uniformly w.r.t. the parameter $\kappa \in [0,1]$ since they
rely on the estimates $ \reff{elliptic0} $ and $ \reff{bound} $
which are uniform w.r.t. $\kappa$, when working with $
\reff{metriqueparametre} $. We will omit this parameter not to
burden the notations but we hope that the proofs are explicit
enough to make this uniformity clear. Then the continuity w.r.t.
$\kappa \in [0,1]$, and actually the smoothness, is obvious since
in the case of the flow (theorem $\refe{Precised}$) it follows
from the standard result of smoothness of O.D.E. with respect to
parameters and in the case of the diffeomorphism (propositions
$\refe{eiko1}$ and $\refe{eiko2}$) it is due to the implicit
functions theorem.

\medskip

We now give the proofs of  these results. We shall only consider
the outgoing situation since the incoming one can be treated
similarly.
\subsection{Proof of proposition ${\bf \refe{aprioritraj}}$ and theorem ${\bf \refe{Precised}}$}
Most of the results that we are going to prove rely on the fact
that
\begin{eqnarray}
 - \partial_r p_{\e} = 2 ( p_{\e} - \rho^2 ) +  {\mathcal
O}(e^{-r}) (p_{\e} - \rho^2) . \label{reason}
\end{eqnarray}
 The first step is the following.
\begin{lemm} \label{time1} One can choose $\bar{R}$ large enough
so that for all  $ (r,y,\rho,\eta) \in \Upsilon^+
(\bar{R},\sigma_+,w,\Omega) $ and all  $\e \in [0,1]$ the
following condition holds:
$$ r_{\e}^{t_1} \geq \bar{R} \ \ \mbox{and} \ \ \rho_{\e}^{t_1} > 0
 \ \ \mbox{for some} \ t_1 \geq 0  \Rightarrow \ r_{\e}^t \geq \bar{R} \ \
\mbox{and} \ \ \rho_{\e}^t \geq \rho_{\e}^{t_1} , \ \ \forall t
\geq t_1.
 $$
\end{lemm}

\noindent {\it Proof.} First we remark that there exists $\bar{R}
> 0$ and $c > 0$ such that
\begin{eqnarray}
- \partial_r p_{\e}
\geq c g (y,e^{-r}\eta) , \qquad r \geq \bar{R}, \ y,\eta \in
\Ra^{n-1} , \ \e \in [0,1]. \label{technic0}
\end{eqnarray}
This follows easily from $\reff{elliptic0}$ and $\reff{bound}$. It
is the first example of application of $\reff{reason}$. Then we
consider the set
$$ {\mathcal T} = \left\{ t \geq t_1 \ | \ \rho^s_{\e} \geq \rho^{t_1}_{\e}
 \ \mbox{and} \ r^s_{\e} \geq \bar{R}, \ \
\forall s \in [t_1,t] \right\}. $$ It is clear that $t_1 \in
{\mathcal T}$ thus ${\mathcal T}$ is a non empty interval. We
shall prove that $T:=\sup {\mathcal T}$ is $+\infty$. We argue by
contradiction and assume that $T$ is finite. By continuity, it is
clear that $T \in {\mathcal T}$, and that
 there exists $T^{\prime}>T$ such that $\rho^s_{\e}
> \rho^{t_1}_{\e}/2$ for all $s \in [t_1,T^{\prime}]$. Now using the fact
$\dot{r}^t_{\e} = 2 \rho^t_{\e}$ we get
$$ r_{\e}^s \geq r_{\e}^{t_1} + (s-t_1) \rho^{t_1}_{\e} \geq \bar{R}, \qquad
t_1 \leq s \leq T^{\prime} . $$ By $\reff{technic0}$, this implies
that $\dot{\rho}^{\e}_s \geq 0 $ on $[t_1,T^{\prime}]$. In
particular $\rho_s^{\e} \geq \rho_{t_1}^{\e}$  for $s$ in a larger
interval than $[t_1,T]$ which is a contradiction. \finpreuve

\medskip

In the next lemma, we show that $\rho^t_{\e}$ reaches a positive
value at some time $t_1$. More precisely we show that
$\rho_{\e}^{t_1} \geq 1/2$ for some positive time $t_1$ which is
uniform with respect to the initial conditions in the outgoing
area.
\begin{lemm} \label{time0} There exists $\widetilde{R}(w,\delta_+) > 0 $  and
$t_1 (w,\delta_+) > 0 $ such that for all initial condition
$(r,y,\rho,\eta) \in \Upsilon^+ (\widetilde{R},\sigma_+,w,\Omega)
$ and all $\e \in [0,1]$
$$ r^{t_1}_{\e} \geq \bar{R} \qquad \mbox{and} \qquad
\rho^{t_1}_{\e} \geq \frac{1}{2} . $$
\end{lemm}

\noindent {\it  Proof.} We first note that the result is clear if
$r \geq \bar{R}$ and $\rho \geq 1/2$ by the previous lemma. Thus
we assume that $ \delta_+ - (1 - w)^{1/2} < \rho < 1/2$. We also
choose $w< \tilde{w} < 1/2$ such that $\tilde{w}-w < \delta_{+}$.
By possibly increasing $\bar{R}$ we may assume that
$\reff{technic0}$ holds and that
\begin{eqnarray}
p_{0}(r,y,\rho,\eta) \in I (w) \  \ \mbox{and} \ \ r \geq \bar{R}
\Rightarrow  p_{\e} (r,y,\rho,\eta) \in I (\tilde{w}) \qquad
\forall \ \e \in [0,1]
\end{eqnarray}
since $p_0-p_{\e}= {\mathcal O}( e^{-r} ) $ on $p_0^{-1}(I(w))$.
We may assume moreover that $\bar{R}$ is large enough so that
\begin{eqnarray}
\left| e^{-2r} \frac{\partial}{\partial r}  e^{-r} \tilde{g}
\left( e^{-r},y,\eta \right) \right| \leq \frac{\delta_+}{2}
\qquad  \mbox{if} \ r \geq \bar{R} \ \ \mbox{and} \ \ p_{\e}
(r,y,\rho,\eta) \in I (\tilde{w})  \ \ \mbox{for all} \ \e \in
[0,1]. \label{time4}
\end{eqnarray}
Now we choose $\widetilde{R}$ such that
$$ \widetilde{R} > \bar{R} + 4 \frac{\sigma_+}{c_+ \delta_{+}},
\qquad \mbox{with} \ c_{+} = 2^{1/2} -1 -\frac{1}{2} .  $$  We
start by proving   that for any $(r,y,\rho,\eta) \in \Upsilon^+
(\widetilde{R},\sigma_+,w,\Omega)$ we have
\begin{eqnarray}
r_{\e}^t \geq \bar{R} \ \ \mbox{and} \ \ \rho^t_{\e} \geq \rho
\qquad \mbox{ for all } \ \ t \geq 0, \ \e \in [0,1].
\label{time2}
\end{eqnarray}
To that end, we proceed as in the previous lemma. We consider the
 set
$$ {\mathcal T} = \{ t \geq 0 \ | \ r^{\e}_{s} \geq \bar{R} \ \ \mbox{and} \ \
\rho^{\e}_{s} \geq \rho, \ \ \forall \ s \in [0,t]  \}  $$ which
is non empty interval since $0 \in {\mathcal T}$ and we show that
$T:=\sup {\mathcal T}$ is $+\infty$. Assume that this is wrong,
then $T$ belongs to ${\mathcal T}$. Note that we may assume that
$\rho^t_{\e} \leq 0$ on $[0,T]$ since otherwise we can find $t_1
\in [0,T]$ such that $\rho^{t_1}_{\e} > 0$ and we get a
contradiction using lemma $\refe{time1}$. Then  $\reff{time4}$ and
the conservation of energy yields
\begin{eqnarray}
\dot{\rho}_{s}^{\e} \geq 2
p_{\e}^{1/2}(p_{\e}^{1/2}+\rho_{s}^{\e}) - \frac{\delta_+}{2} \geq
c_+ \delta_+
\end{eqnarray}
where the second inequality follows from the fact that
$p_{\e}^{1/2}
> (1/2)^{1/2}$ and $$ p_{\e}^{1/2}+\rho \geq  (1-\tilde{w})^{1/2}
- (1-w)^{1/2} + \delta_{+} \geq (1-1/2^{1/2}) \delta_+ . $$  As a
consequence we get
 $ \rho^T_{\e} - \rho \geq T c_+ \delta_+ $ and this implies easily that
$$ T < \frac{\sigma_{+}}{c_+ \delta_+}  $$
since $\rho^{T}_{\e} \leq 0 $ and $\rho > \delta_+ - (1-w)^{1/2}$.
On the other hand, the first equation of motion yields
$$ r_{\e}^T \geq r - 2 T E_{+}(w)^{1/2} > \widetilde{R} - 4 T > \bar{R} . $$
Thus $\rho^{T}_{\e} > \rho$ and $r_{\e}^T > \bar{R}$. This shows
that
 there exists
$T^{\prime} > T$ such that $\rho_{s}^{\e}\geq \rho_0$ and
$r_{s}^{\e}\geq \bar{R}$ on $[0,T^{\prime}]$ which is a
contradiction and completes the proof of $\reff{time2}$.

We shall now prove the existence of $t_1$ in the same spirit.
Recall that we can assume that $\rho<1/2$. Then there exists $t_0$
depending on $\e$ and the initial conditions such that
$\rho^{t_0}_{\e} = 1/2$ otherwise $\rho^t_{\e}<1/2$ for all $t
\geq 0 $ and
\begin{eqnarray}
 \dot{\rho}^t_{\e} \geq 2 (p_{\e}^{1/2} + \rho^t_{\e})(
p_{\e}^{1/2} - \rho^t_{\e}) - \frac{\delta_+}{2} \geq 2
(2^{1/2}-1/2)(1 - 2^{-1/2}) \delta_{+} - \frac{\delta_+}{2} >
\frac{\delta_+}{2} \label{lowerbtime}
\end{eqnarray}
since we can use $\reff{time4}$ thanks to $\reff{time2}$. This
implies that $\rho^t_{\e} \geq \rho + t \delta_+ /2 \rightarrow +
\infty$ as $t \rightarrow \infty$ which is forbidden by the energy
conservation. This is a contradiction thus $\rho_{\e}^t$ reaches
the value $1/2$ at some positive time. Moreover
$\reff{lowerbtime}$ holds as long as $\rho_{\e}^t$ is lower than
$1/2$, thus if $t_0$ is the smallest positive time for which
$\rho_{\e}^{t_0}=1/2$ we have $1/2 - \rho \geq t_0 \delta_+ / 2$.
This yields
$$ t_0 < \frac{2}{\delta_+} (\sigma_+ + 1/2)  . $$
Using lemma $\reff{time1}$ it is clear that $\rho_{\e}^t \geq 1/2$
for $t \geq (2 \sigma_+ +1 ) / \delta_+$ and this completes the
proof. \finpreuve

\medskip

\noindent {\it Proof of proposition $ \refe{aprioritraj}$.}
We choose  $R>\widetilde{R}$ with $\widetilde{R}$ as in lemma
$\reff{time0}$. Then for $t \geq t_1$, with $t_1$ as in lemma
$\refe{time0}$, we have
\begin{eqnarray}
r_{\e}^t  \geq  (t-t_1) + r^{t_1}_{\e}  .
\end{eqnarray}
On the other hand, for $t \in [0,t_1]$ we have
$$ r^{t}_{\e} - r \geq - 2 c_1 t $$ for any
 $c_1 > p_{\e}(r,y,\rho,\eta)$ (note that such a
$ c_1 $ can be chosen uniformly with respect to $\e$ and the
initial conditions in $ \Upsilon^{+}(R,\sigma_+,w,\Omega) $).
These two estimates yields $\reff{rt}$. The proof of $\reff{rhot}$
follows directly from $\reff{reason}$ since one knows that
$e^{-r^t_{\e}}$ is small uniformly for $t \geq 0$. Finally, the
conservation of energy implies that $e^{-r^t_{\e}}\eta^t_{\e}$ is
bounded and thus $\dot{y}^t_{\e} = {\mathcal O}(e^{-r-t})$, by
$\reff{rt}$. This implies $ \reff{yt} $ and completes the proof.
\finpreuve

\bigskip

The rest of this subsection is now devoted to the proof of theorem
$\refe{Precised}$. We start with the following lemma

\begin{lemm} \label{Lemme2}   If $ R $ is large
enough, there exists $C > 0$ such that
$$ |\eta^t_{\e}-\eta| \leq C  e^{-r}|\eta|  $$
for all $(r,y,\rho,\eta) \in \Upsilon^{+}
(R,\sigma_{+},w,\Omega)$, $ \e \in [0,1] $  and $t \geq 0$.
\end{lemm}

\noindent {\it Proof.} Choosing $R$ large enough and using
proposition $ \refe{aprioritraj}
 $ shows the existence of $c > 0$ such that
$$|(\partial_y g_{\e})(r^t_{\e},y_{\e}^t,\eta^t_{\e}) | \leq c g_{\e}(r^t_{\e},y_{\e}^t,\eta^t_{\e}) $$
since $r^t_{\e}$ is large. The last motion equation and the
conservation of energy then show that $\dot{\eta}^t_{\e}$ is
bounded and thus $\eta^t_{\e} - \eta = {\mathcal O}(t)$. Putting
this new estimate into the last motion equation and using the fact
that $e^{- 2 r^t_{\e}} = {\mathcal O}(e^{-2r-2t})$ show that
$\eta^t_{\e}-\eta$ is bounded. Since $e^{-r}\eta$ is bounded on
$\Upsilon^+$, the lemma is proved provided $e^{-r}\eta$ is away
from any neighborhood of $0$. Thus we assume now that $e^{-r}\eta$
is small enough so that $\rho \geq 0$. Since $\dot{\rho}^t_{\e}
\geq 0$ and $\rho \geq 0$ we have
$$ \rho^2 \leq (\rho^t_{\e})^2 \leq p_{\e} $$
and then, using the fact that
$e^{-2r_{\e}^t}g_{\e}(r_{\e}^t,y_{\e}^t,\eta_{\e}^t) = p_{\e} -
(\rho_{\e}^t)^2 $ we obtain
$$ |\dot{\eta}^t_{\e}| \leq c \left( p_{\e} - \rho^2 \right) \leq c C_0 (e^{-r}|\eta|)^2 . $$
This shows that $\eta^t_{\e}-\eta = {\mathcal O}(t e^{-r}|\eta|)$
and putting this estimate into the motion equation of
$\eta^t_{\e}$ as before we obtain $\eta^{t}_{\e} - \eta =
{\mathcal O}(e^{-r}|\eta|)$ thanks to the exponential decay in
time of $e^{- r^t_{\e}}$. \finpreuve

\medskip

We explain the strategy of the proof of theorem $\refe{Precised}$
with the case $\partial^{\gamma} =
\partial_{\e}$.  By standard results on ordinary differential equations we
know that $\phi^t_{\e}$ is smooth with respect to $\e$. In
particular, if we consider the matrix
\begin{eqnarray}
 M_{\e}(r,y,\rho,\eta) = \left(
\begin{matrix}
\partial^2_{r \rho}  p_{\e} & \partial^2_{y \rho} p_{\e} &
\partial^2_{\rho \rho} p_{\e}  & \partial^2_{\eta \rho} p_{\e} \\

\partial^2_{r \eta} p_{\e} &  \partial^2_{y \eta} p_{\e} &
\partial^2_{\rho \eta} p_{\e}  & \partial^2_{\eta \eta} p_{\e} \\

- \partial^2_{r r} p_{\e} & -\partial^2_{y r} p_{\e} &
- \partial^2_{\rho r} p_{\e}  & - \partial^2_{ \eta r} p_{\e} \\

- \partial^2_{r y} p_{\e}  & - \partial^2_{y y} p_{\e} &
- \partial^2_{\rho y} p_{\e}  & - \partial^2_{\eta y} p_{\e} \\
\end{matrix}
\right)  \label{Lflow}
\end{eqnarray}
applying $\partial_{\e}$ to the motion equations yields
$$   \dot{X}^t_{\e} = M_{\e}(t) X_{\e}^t + Y_{\e}^t  $$ with the notations $
X_{\e}^t = \partial_{\e} \phi^t_{\e}$, $ M_{\e}(t) = M_{\e}
(\phi^{\e}_t)$ and $ Y_{\e}^t = \left( \partial_{\e} {\mathbf
H}_{\e} \right) (\phi_t^{\e}) $.
 By $\reff{elliptic0}$ and $\reff{bound}$, we have
\begin{eqnarray}
 M_{\e} (t) & = & M + {\mathcal O} \left(  \exp (- 2 r_{\e}^t )
g_{\e} \left( r_{\e}^t  , y_{\e}^t , \eta_{\e}^t \right) \right),
\label{Initial} \\
 Y_{\e}^t & = &  {\mathcal O} \left(  \exp (- 3 r_{\e}^t )
g_{\e} \left( r_{\e}^t , y_{\e}^t , \eta_{\e}^t \right) \right)
\end{eqnarray}
where the matrix $M$ is given by
 $$ M  = \left(
\begin{matrix}
0 & 0 & 2 & 0 \\
0 & 0 &  0 & 0 \\
0 & 0 & 0  & 0 \\
0 & 0 & 0  & 0 \\
\end{matrix}
\right) .  $$ In particular, it is very easy to check that $M^2=0$
and this implies that
$$ \exp (tM) = 1 + t M , \qquad t \in \Ra . $$ We do a  change of unknown
function by considering $\widetilde{X}_{\e}^t = \exp (- t M)
X_{\e}^t$. It satisfies the equation
\begin{eqnarray}
\dot{\widetilde{X}_{\e}^t} & = & \widetilde{M}_{\e}(t)
\widetilde{X}_{\e}^t + \widetilde{Y}_{\e}^t \label{Equationtilde}
\end{eqnarray}
where the matrix $\widetilde{M}_{\e} (t)$ and the vector
$\widetilde{Y}_{\e}^t$ are given by
\begin{eqnarray}
\widetilde{M}_{\e}(t) & = & \exp (- t M )(M_{\e}(t)-M) \exp ( t M ) \nonumber \\
\widetilde{Y}_{\e}^t & = & \exp(- t M ) Y_{\e}^t . \nonumber
\end{eqnarray}
Lemma $ \refe{Lemme2} $ and proposition $ \refe{aprioritraj} $
imply easily the following key estimates
\begin{eqnarray}
\widetilde{M}_{\e}(t) & = & {\mathcal O} \left( \scal{t}^2
e^{-2t-r} |\eta |
  \right),  \label{Estimetilde} \\
\widetilde{Y}_{\e}^t & = & {\mathcal O} \left( \scal{t}  \exp( - 2
t - 2 r ) |\eta| \right). \label{Estimetilde2}
\end{eqnarray}
We shall deduce the estimates on $\widetilde{X}^t_{\e}$ from
$\reff{Equationtilde}$, $\reff{Estimetilde}$,$\reff{Estimetilde2}$
and the well known Gronwall's lemma which we quote under the
following form.
\begin{lemm}[Gronwall's lemma] Let $u(t) \geq 0$ be a continuous function
 on
$[t_0,+\infty)$ for some $t_0 \geq 0$. Assume that
$$  u (t) \leq a \int_{t_0}^t u (s) \ \emph{d}s + b, \qquad t \geq t_0 $$
for some $a \geq 0$ and $b \geq 0$. Then we have
$$  u (t) \leq b e^{a(t-t_0)}, \qquad t \geq t_0.  $$
\end{lemm}
Since $\widetilde{X}_{\e}^0  = X_{\e}^0=0$, we can turn
$\reff{Equationtilde}$ into the following integral equation
$$  \widetilde{X}_{\e}^t = \int_0^t \widetilde{M}_{\e}(s)
\widetilde{X}_{\e}^s \ \mbox{d}s + \int_{0}^t \widetilde{Y}_{\e}^s
 \ \mbox{d}s. $$
Then Gronwall's lemma combined with the estimates
$\reff{Estimetilde2}$ shows that for any $T \geq 0$, there exists
$C_T$ independent of the initial conditions in $ \Upsilon^{+}
(R,\sigma_{+},w,\Omega) $ and of $\e$ such that
\begin{eqnarray}
|\widetilde{X}_{\e}^t| \leq C_T  e^{-2 r} |\eta| , \qquad 0 \leq t
\leq T . \label{infT}
\end{eqnarray}
We can also write $\reff{Equationtilde}$ as follows
$$ \widetilde{X}_{\e}^t = \int_T^t \widetilde{A}_{\e}(s)
\widetilde{X}_{\e}^s \ \mbox{d}s + \int_{0}^t
\widetilde{Y}^{\e}(s)
 \ \mbox{d}s + \widetilde{X}_{\e}^T , \qquad t \geq T  . $$
In particular, by choosing $T$ large enough we can show easily
that
$$ |\widetilde{X}_{\e}^t| \leq \int_T^t
|\widetilde{X}^{\e}(s)| \ \mbox{d}s + C_{T}^{\prime}  e^{- 2
r}|\eta|  , \qquad t \geq T $$ for some $C_T^{\prime}
> 0$. This is a simple consequence of the fact that
$\widetilde{M}_{\e}(t) \rightarrow 0$ as $t \rightarrow + \infty$.
Thus another application of Gronwall's lemma  shows that for some
$C_T^{\prime} \geq 0$
\begin{eqnarray}
 |\widetilde{X}_{\e}^t| \leq C^{\prime}_T  e^{- 2 r} |\eta|
e^{(t-T)}, \qquad t
 \geq T. \label{supT}
\end{eqnarray}
Putting $\reff{infT}$ and $\reff{supT}$ into
$\reff{Equationtilde}$ and using $\reff{Estimetilde}$, one gets
$$   \dot{\widetilde{X}_{\e}^t} = {\mathcal O} \left( \scal{t}^2 e^{-t}
 e^{- 2r}|\eta| \right) , \qquad t \geq 0. $$ We can now come
back to $X_{\e}^t$ since one checks easily that $\dot{X}_{\e}^t =
M X^t_{\e} + \exp(tM)\dot{\widetilde{X}_{\e}^t}$ which implies
$$
X_{\e}^t = \exp(tM) \int_0^t \dot{\widetilde{X}_{\e}^s} \
\mbox{d}s . $$ Using the explicit form of $\exp(tM)$ we conclude
that
\begin{eqnarray}
\partial^{\e} r^t_{\e}  =  {\mathcal O} \left( \scal{t} e^{- 2 r } |\eta| \right), \\
|\partial^{\e} y^t_{\e}| +
 |\partial^{\e} \rho^t_{\e} |
+ | \partial^{\e} \eta^t_{\e}|   =  {\mathcal O} \left(  e^{- 2 r}
|\eta| \right) .
\end{eqnarray}
This completes the proof of theorem $\refe{Precised}$ if
$\partial^{\gamma} =
\partial_{\e}$. The other estimates for $|\gamma| = 0$ can be
proved easily using lemma $\refe{Lemme}$ and the motion equations.
For $|\gamma| \ne 0$, we proceed by induction by applying
$\partial^{\gamma}$ to the motion equations. \finpreuve

\medskip

\noindent {\bf Remark.} Using the results of theorem
$\refe{Precised}$, we can improve the estimates on $y^t_{\e}$.
Since
$$ y^t_{\e} = y + \int_0^t e^{-2 r^s_{\e}} \partial_{\eta} g_{\e}(r^s_{\e},y^s_{\e},\eta^s_{\e}) \ \mbox{d}s $$
where $e^{-r^s_{\e}}\partial_{\eta}
g_{\e}(r^s_{\e},y^s_{\e},\eta^s_{\e}) $ is bounded, we obtain
easily
\begin{eqnarray}
|\partial^{\gamma}(y^t_{\e}-y)| \leq C_{\gamma} e^{-(j_+ +
1)r}|\eta|, \qquad j_+ = \max (j,1). \label{yimproved}
\end{eqnarray}
\subsection{Proofs  propositions ${\bf \refe{eiko1}}$ and ${\bf \refe{eiko2}}$}
We start with a technical lemma.
\begin{lemm} \label{connected} Assume that $\Omega$ is connected. Then $\Gamma^{+}
(R,\varepsilon,w,\Omega)$ is connected.
\end{lemm}

\noindent {\it Proof.} Consider $(r_0,y_0,\rho_0,\eta_0)$ and
$(r_1,y_1,\rho_1,\eta_1)$ in $\Gamma^{+}
(R,\varepsilon,w,\Omega)$. We prove the result under the
conditions $\eta_0 \ne 0$ and $ \eta_1 \ne 0 $ (the other cases
are simpler). We first connect $(r_0,y_0,\rho_0,\eta_0)$ to
$(r_1,y_1,\rho_0, \tau \eta_0)$ for some $\tau > 0$ by considering
$$ ((1-s)r_0 + s r_1  , y (s) , \rho_0 , \tau (s) \eta_0) , \qquad s \in [0,1] $$
with $y(s)$ a path joining $y_0$ to $y_1$ in $\Omega$ and
$\tau(s)$ the unique positive continuous function such that
$$  \exp \left( - 2 ((1-s)r_0 + s r_1) \right) g (y(s),\tau(s)\eta_0) =  \exp \left( - 2 r_0 \right) g (y_0,\eta_0) . $$
Then we connect $(r_1,y_1,\rho_0,\tau \eta_0)$ to
$(r_1,y_1,\rho_1, \eta_1)$ using a path $(\rho (s), \eta (s))$
such that
$$ \eta (0) = \tau \eta_0 , \qquad \eta (1) = \eta_1, \qquad {\mathcal Q}(\rho(s), \eta (s)) \in [E_0,E_1] , \qquad
E_j = p_0 (r_j,y_j,\rho_j,\eta_j), \ \ j =1,2,  $$ where
${\mathcal Q}$ is the quadratic form on $\Ra^n$ defined by
${\mathcal Q }(\rho,\eta) = \rho^2 + e^{-2r_1}g(y_1,\eta)$. This
is possible since the fact that $n \geq 2$ implies that the region
defined by ${\mathcal Q}(\rho,\eta) \in [E_0,E_1]$ and $\rho > 0$
is connected in $\Ra^{n}$.
 \finpreuve

\medskip

We split the proof of proposition $\refe{eiko1}$ into three
lemmas.

\begin{lemm} \label{methodediffeo} 
Let $0<w<1/2$ and $\Omega \subset \Ra^{n-1}$ be an open subset.
There exists $\varepsilon > 0$ small enough and $R$ large enough
such that $\Phi_{\e}^t$ is diffeomorphism from $ \Gamma^{+}
(R,\varepsilon,w,\Omega)$ onto its range, for all $\e \in [0,1]$
and $t \geq 0$.
\end{lemm}

\noindent {\it Proof.} The estimates $\reff{rhot2}$ and
$\reff{etat2}$ show that $\nabla \Phi_{\e}^t$ is as close to the
identity as we want by choosing $ \varepsilon$ small enough. Under
such a condition, $\Phi_{\e}^t$ is a local diffeomorphism onto its
range, thus it is global iff it is injective. If $\Phi_{\e}^t
(r,y,\rho,\eta) = \Phi_{\e}^t
(r^{\prime},y^{\prime},\rho^{\prime},\eta^{\prime})$ then $r =
r^{\prime}$ and $y = y^{\prime}$. Moreover $\reff{rhot2}$ and
$\reff{etat2}$ implies that
\begin{eqnarray}
 \left(
\rho_{\e}^t,\eta_{\e}^t  \right) (r,y,\rho,\eta) - \left(
\rho_{\e}^t, \eta_{\e}^t \right)
 (r,y,\rho^{\prime},\eta^{\prime}) = (\rho - \rho^{\prime} , \eta -
 \eta^{\prime}) + {\mathcal O}(\varepsilon)
 \label{Estimeinjectif}
\end{eqnarray}
thus if the left hand side vanishes we have $\rho - \rho^{\prime}
= {\mathcal O}(\varepsilon)$ and $\eta - \eta^{\prime} = {\mathcal
O}(\varepsilon)$. On the other hand, Taylor's formula to the
second order combined with $\reff{rhot2}$ and $\reff{etat2}$ show
that
\begin{eqnarray*}
 \left( \rho_{\e}^t,\eta_{\e}^t  \right) (r,y,\rho,\eta)
- \left( \rho_{\e}^t, \eta_{\e}^t \right)
 (r,y,\rho^{\prime},\eta^{\prime})  =  (1 - {\mathcal
 O}(\varepsilon))(\rho - \rho^{\prime} , \eta -
 \eta^{\prime}) + {\mathcal O} (\varepsilon) |(\rho - \rho^{\prime} , \eta -
 \eta^{\prime} )|^2
 \end{eqnarray*}
This implies that
 $$ | \left( \rho_{\e}^t,\eta_{\e}^t  \right) (r,y,\rho,\eta)
- \left( \rho_{\e}^t, \eta_{\e}^t \right)
 (r,y,\rho^{\prime},\eta^{\prime}) |  =  (1 - {\mathcal
 O}(\varepsilon))|(\rho - \rho^{\prime} , \eta -
 \eta^{\prime})| $$
 if the left hand side of $\reff{Estimeinjectif}$ vanishes
 and shows that $ \rho = \rho^{\prime} $ and $\eta = \eta^{\prime}$, if $\varepsilon$ is small enough.  \finpreuve

\begin{lemm} Let $\Omega_1,\Omega_2$ be connected open subsets of $\Ra^{n-1}$ such that
$\Omega_1 \csubset \Omega_2$. Then for all $\varepsilon_2$ small
enough and $R_2$ large enough there exists exists $w_1,w_2,
\varepsilon_1 > 0 $ small enough and $R_1 > 0$ large enough so
that $\reff{Part0}$ holds for all $t \geq 0$ and all $\e \in
[0,1]$.
\end{lemm}

\noindent {\it Proof.} The proof follows a rather standard scheme
(see \cite{Robe1}) . It is enough to find parameters such that
\begin{eqnarray}
 \Gamma^{+} (R_1,\varepsilon_1,w_1,\Omega_1) \cap \Phi_{\e}^t \left(
\partial \Gamma^{+} (R_2,\varepsilon_2,w_2,\Omega_2) \right) = \emptyset
\label{Intevide} \\
 \Gamma^{+} (R_1,\varepsilon_1,w_1,\Omega_1) \cap \Phi_{\e}^t \left(
 \Gamma^{+} (R_2,\varepsilon_2,w_2,\Omega_2) \right) \ne \emptyset .
\label{Part}
\end{eqnarray}
Actually $\reff{Part}$ is always satisfied since any element of
the form $(r,y,1,0)$, with $r$ large, is invariant by
$\Phi_{\e}^t$. Thus we look for conditions on the parameters
$R_j,\varepsilon_j,w_j, \Omega_j$ ensuring the fact that
$\reff{Intevide}$ holds for all $t \geq 0 $ and $\e \in [0,1]$.
Any element in the intersection of $\reff{Intevide}$ can be
written
$$ ( r,y,\rho,\eta ) = ( r , y , \rho_{\e}^t (r,y,\rho_0,\eta_0) ,
\eta_{\e}^t (r,y,\rho_0,\eta_0) ) $$ with $( r,y,\rho_0,\eta_0 )
\in \partial \Gamma^{+} (R_2,\varepsilon_2,w_2,\Omega_2) $. Since
$\Omega_1 \csubset \Omega_2$,  choosing $R_1 > R_2$ implies that
one of the following two conditions is satisfied
\begin{eqnarray}
\rho_0^2 + e^{-2r}g (y,\eta_0) & = & 1 \pm w_2 , \label{Condition01} \\
e^{-2r} g (y,\eta_0) & = & \varepsilon_2 . \label{Condition02}
\end{eqnarray}
Using $\reff{rhot2}$ and $\reff{etat2}$ we see that
\begin{eqnarray}
 \rho^2 + e^{-2r} g (y,\eta) & = &  \rho_0^2 + e^{-2r}g (y,\eta_0) +
{\mathcal O} \left( \varepsilon_2^{1/2}  \right)
\nonumber \\
& = & 1 \pm w_2 \left( 1 + {\mathcal O} \left( \varepsilon_2^{1/2}
/ w_2  \right)   \right) \nonumber .
\end{eqnarray}
Thus $\reff{Condition01}$ can certainly not happen if the
following conditions are satisfied simultaneously
\begin{eqnarray}
w_1 \ll w_2, \qquad \varepsilon_2 \ll w_2^2 ,
 \label{Contradiction1}
\end{eqnarray}
since $\rho^2 + e^{-2r}g(y,\eta) = 1 + {\mathcal O} \left( w_1
\right)$. Similarly, we deduce easily from theorem
$\refe{Precised}$ that
$$ e^{-2r} g (y,\eta) = e^{-2r} g (y,\eta_0) + {\mathcal O} \left( e^{-r} \right) =
\varepsilon_2 \left( 1 +  {\mathcal O} \left( e^{-R_1} /
\varepsilon_2 \right)  \right) . $$ Thus the conditions
$e^{-2r}g(y,\eta) < \varepsilon_1  $ and $\reff{Condition02}$
cannot hold simultaneously  if
\begin{eqnarray}
\varepsilon_1 \ll \varepsilon_2 \qquad \mbox{and} \qquad e^{-R_1}
\ll \varepsilon_2 . \label{Contradiction2}
\end{eqnarray}
The existence of parameters satisfying $\reff{Contradiction1}$ and
$\reff{Contradiction2}$ is clear. We choose can choose
 for instance
$$
 w_1 = w_2^2 , \qquad \varepsilon_2 = w_2^3 \qquad \varepsilon_1 = w_2^4 , \qquad
e^{-R_1} = w_2^4 , \qquad R_1 > R_2 . $$ Thus $\reff{Intevide}$
and $\reff{Part}$ hold. Since $\Gamma^{+}(R_1,
\varepsilon_1,w_1,\Omega_1)$ is connected and $\Phi_{\e}^t$ is an
homeomorphism on a neighborhood of $\overline{\Gamma^{+}(R_2,
\varepsilon_2,w_2,\Omega_2)}$ onto its range (by possibly
decreasing $w_2,\varepsilon_2,\Omega_2$ and increasing $R_2$),
this implies $\reff{Part0}$ . \finpreuve

\begin{lemm} If $\varepsilon_1 > 0$ is small enough, the estimates
$\reff{Diffe}$ hold on $\Gamma^+ (R_1,\varepsilon_1,w_1,\Omega_1)$
for $t \geq 0$ and $\e \in [0,1]$.
\end{lemm}

\noindent {\it Proof.} By $\reff{Part0}$ we know that
$(r,y,\tilde{\rho}_{\e}^t,\tilde{\eta}_{\e}^t) \in
\Gamma^{+}(R_2,\varepsilon_2,w_2,\Omega_2)$ if $(r,y,\rho,\eta)
\in \Gamma^{+}(R_1,\varepsilon_1,w_1,\Omega_1)$. Moreover if $w_2$
is small enough, then $\Gamma^{+}(R_2,\varepsilon_2,w_2,\Omega_2)
\subset \Upsilon^+ (R_2, - 1/2 ,w_2,\Omega_2)$. Thus we can use
the estimates $\reff{rhot2}$ and $\reff{etat2}$ with the initial
conditions $ (r,y,\tilde{\rho}_{\e}^t (r,y,\rho,\eta),
\tilde{\eta}_{\e}^t (r,y,\rho,\eta))$. This shows that
$$ | \tilde{\eta}_{\e}^t - \eta
| \leq C e^{-r}|\tilde{\eta}_{\e}^t| $$ and by choosing $r$ large
enough so that $C e^{-r}<1/2$ we get
$$ |
\tilde{\eta}_{\e}^t - \eta | \leq 2 C e^{-r}|\eta| .
$$
 Using this estimate and theorem $\refe{Precised}$, we obtain
similarly
$$ | \tilde{\rho}_{\e}^t - \rho | \leq C e^{-r} |\tilde{\eta}_{\e}^t| \leq \tilde{C}
e^{-r}|\eta| . $$ This shows that $\reff{Diffe}$ holds with
$\gamma=0$. For $|\gamma| \geq 1$ we proceed by induction, by
differentiating the equality $\Phi_{\e}^t \circ \left( \Phi_{\e}^t
\right)^{-1} = \mbox{id} $. For instance, if $\partial^{\gamma} =
\partial_{\e} $, we have
$$  ( \nabla \Phi_{\e}^t )_{|\left( \Phi_{\e}^t
\right)^{-1} } . \partial_{\e} \left( \Phi_{\e}^t \right)^{-1} = -
( \partial_{\e} \Phi_{\e}^t )_{|\left( \Phi_{\e}^t \right)^{-1} }
$$
where the right hand side is ${\mathcal O}(e^{-r})$ and, on the
left hand side, $( \nabla \Phi_{\e}^t )-1$ is small enough hence
$$
\partial_{\e} \left( \Phi_{\e}^t \right)^{-1} = - ( \nabla \Phi_{\e}^t
)^{-1} . \partial_{\e} \! \left( \Phi_{\e}^t \right)^{-1} =
{\mathcal O}(e^{-r}) . $$ We don't go any further into details.
\finpreuve

\medskip

This lemma completes the proof of proposition $\refe{eiko1}$.

\medskip

\noindent

{\it Proof of proposition $\refe{eiko2}$.}
The fact that $S_{\e}^{+}$, defined as the composition of the right hand
side of $ \reff{eikonalexplicit} $ with the inverse of
$\Phi^t_{\e}$, solves $ \reff{eikonal} $ is a standard result. See
for instance \cite{DeGe} or \cite{Robe1}. Thus we focus on the
proof of $\reff{approxeiko}$. We first remark that
$\reff{eikonalexplicit}$ can be rewritten as follows
\begin{eqnarray}
r_{\e}^t \rho^t_{\e} + y^t_{\e}.\eta^t_{\e} - 2 \int_{0}^t \left(
\rho^s_{\e} \right)^2 \mbox{d} s + t p_{\e} - \int_0^t
e^{-r^s_{\e}}
\partial_{\eta} g_{\e}
(r_{\e}^s,y_{\e}^s,\eta^s_{\e}).e^{-r^s_{\e}}\eta^s_{\e} \
\mbox{d}s . \label{rewrite}
\end{eqnarray}
This is easily obtained by integration by part using the motion
equations.
Furthermore
$$ \int_{0}^t \left(
\rho^s_{\e} \right)^2 \mbox{d}s = t p_{\e} - \int^t_0 e^{-2
r^s_{\e}} g_{\e}(r_{\e}^s, y_{\e}^s, \eta^s_{\e}) \ \mbox{d}s  $$
and this implies that
\begin{eqnarray}
 \reff{eikonalexplicit} = r^t_{\e}\rho^t_{\e} +
y^t_{\e}.\eta^t_{\e} - t p_{\e} - \! \! \int_0^t \! e^{-r^s_{\e}}
\partial_{\eta} g_{\e}
(r_{\e}^s,y_{\e}^s,\eta^s_{\e}).e^{-r^s_{\e}}\eta^s_{\e} - 2 e^{-2
r^s_{\e}} g_{\e}(r_{\e}^s, y_{\e}^s, \eta^s_{\e})   \mbox{d}s .
\label{eikonalbase}
\end{eqnarray}
where one must notice that the integral is convergent and
${\mathcal O}(e^{-r}|\eta|)$. We rewrite $\reff{eikonalbase}$
using the fact that $p_{\e} = p_{\e} \circ \phi^t_{\e}$. In
particular, if we note $\tilde{r}^s_{\e} =
r^s_{\e}(r,y,\tilde{\rho}^t_{\e},\tilde{\eta}^t_{\e})$ and
$\tilde{y}^s_{\e} = y^s_{\e}
(r,y,\tilde{\rho}^t_{\e},\tilde{\eta}^t_{\e}) $ we obtain
\begin{eqnarray*}
 S_+ (t,r,y,\rho,\eta,\e) & = & \tilde{r}^t_{\e}\rho +
\tilde{y}^t_{\e}.\eta - t \left( \rho^2 +
e^{-2\tilde{r}^t_{\e}}g_{\e}(\tilde{r}^t_{\e},\tilde{y}^t_{\e},\eta)
\right) + {\mathcal O} \left( e^{-r}|\eta| \right) \\
& = & \tilde{r}^t_{\e}\rho + y.\eta - t  \rho^2  + {\mathcal O}
\left( e^{-r}|\eta| \right)
\end{eqnarray*}
using $\reff{yimproved}$ to estimate $\tilde{y}^t_{\e}-y$ and the
exponential decay w.r.t. $t$ of $
e^{-2\tilde{r}^t_{\e}}g_{\e}(\tilde{r}^t_{\e},\tilde{y}^t_{\e},\eta)
$. In order to estimate $\tilde{r}^t_{\e}$, we use the motion
equations to write
$$ \tilde{r}^t_{\e} = r + 2 \int_0^t \rho^s_{\e} \left(r,y,\tilde{\rho}^t_{\e},\tilde{\eta}^t_{\e} \right) \ \mbox{d}s
, \qquad  \rho^s_{\e}
\left(r,y,\tilde{\rho}^t_{\e},\tilde{\eta}^t_{\e} \right) = \rho +
\int_s^t (\partial_r p_{\e}) \left(
\tilde{r}^u_{\e},\tilde{y}^u_{\e},
\eta^u_{\e}(r,y,\tilde{\rho}^t_{\e},\tilde{\eta}^t_{\e}) \right) \
\mbox{d} u . $$ By $\reff{rt}$ and theorem $\refe{Precised}$, we
see that the second integral is $ {\mathcal O}(e^{-s}
e^{-r}|\eta|) $ and this implies that
$$  \tilde{r}^t_{\e} = r + 2 t \rho + {\mathcal O}( e^{-r}|\eta|) . $$
and finally we obtain
$$ S_{\e}^{+}
(t,r,y,\rho,\eta,\e) = r \rho + y. \eta + t \rho^2 + {\mathcal O}(e^{-r}|\eta|) . $$
 Note that the remainder ${\mathcal
O}(e^{-r}|\eta|)$ is explicit and theorem $ \refe{Precised} $
combined with proposition $ \refe{eiko1} $ show that
$\reff{approxeiko}$ holds. This completes the proof of proposition
$\refe{eiko2}$. \finpreuve

\section{Pseudo-differential operators} \label{pseudo}

\setcounter{equation}{0}
\subsection{Local theory}
In this part, we consider pseudo-differential operators on $\Ra^n$. We apply the
results to operators on $X$ in the next subsection.

\begin{defi} \label{symbol} For any $m,m^{\prime} \in \Ra$, the class
$S^{m,m^{\prime}}$ is the set of smooth
functions $a$ such that, for all $R > 0$, $ \Omega \Subset
\Ra^{n-1} $ and $k,l,\alpha,\beta$ there exists $C$ such that
$$ \left| \partial_r^k \partial_{\rho}^l \partial_y^{\alpha} \partial_{\eta}^{\beta} a (r,y,\rho,\eta) \right|
\leq C  \scal{\rho}^{m} \scal{e^{-r}\eta}^{m^{\prime}} ,
$$
for all $ r \geq R , \ y \in \Omega , \ \rho \in \Ra , \ \eta \in
\Ra^{n-1}$. We set $ S^{ - \infty} =
\cap_{m,m^{\prime}} S^{m,m^{\prime}}
 $.
\end{defi}
As usual, the best constants $C$ are semi-norms which define the
topology of $ S^{m,m^{\prime}} $. We also mention that we shall mainly consider
 cases where $m , m^{\prime} \in \Ra^-$.

We give two examples of special interest for us.
 If $b$ belongs to
$ S^{m,m^{\prime}}_{0,1} $, that is if
$$ \left| \partial_r^k \partial_{\rho}^l \partial_y^{\alpha} \partial_{\eta}^{\beta} b (r,y,\rho,\eta) \right|
\leq C  \scal{\rho}^{m} \scal{\eta}^{m^{\prime}-  |\beta|}
$$
 for $r > 0$, $y \in \Omega$, $\rho \in \Ra$ and $\eta \in
 \Ra^{n-1}$, then one checks easily that the function $ a $ defined by
\begin{eqnarray}
   a (r,y,\rho,\eta) = b \left( r , y , \rho , e^{-r} \eta \right)
\label{formprincipal}
\end{eqnarray}
 is an element of $S^{m,m^{\prime}}$. In particular, if $f \in {\mathcal
 S}$  then
 $ f \left( \rho^2 + e^{-2r} g (y,\eta) \right) \in S^{ - \infty} . $
The second example is the following.
If  $ \Omega_0 \Subset \Omega $ and $  w^{\prime} > w  $, there
exists $C > 0$ such that for all $R$ large enough and all
 $ a \in S^{ - \infty} $
 suppported in $  \Upsilon^{+} (R,\sigma_{+},w,\Omega_0) $
(see $ \reff{presquesortante} $) we have
\begin{eqnarray}
 a \circ \phi_{\e}^{-t} \in S^{ - \infty}
  \qquad \mbox{and} \qquad
 \mbox{supp} \ a \circ \phi_{\e}^{-t} \subset
\Upsilon^{+} (R + t - C,\sigma_{+}, w^{\prime} ,\Omega)
 , \label{transportcutoff}
\end{eqnarray}
for all $ t \geq 0 $ and $ \e \in [0,1] $.
This follows from proposition $ \refe{aprioritraj} $ and theorem
$ \refe{Precised} $.

The pseudo-differential operators that we will use
are of the form $ O \! p_h (a) $, with
\begin{eqnarray}
 O \! p_h (a) u (r,y) = (2 \pi)^{-n} \int \! \! \int e^{ir\rho +
i y.\eta} a(r,y,h\rho,h\eta) \widehat{u}(\rho,\eta) \ \mbox{d}
\rho \mbox{d} \eta , \qquad h \in (0,1], \label{opha}
\end{eqnarray}
 where
$\hat{u}(\rho,\eta) = \int \! \! \int
e^{-ir\rho-iy.\eta}u(r,\rho)\ \mbox{d}r \mbox{d}y $ is the
Fourier transform of $u \in {\mathcal S}$ . They depend on the parameter
$h \in (0,1]$ and it is natural to consider symbols depending on
$h$ as well. Following the standard definitions of
\cite{Robe1}, we  say that $a = a(h)$ is an admissible
symbol in $ S^{m,m^{\prime}} $ and note
$$ a \sim a_0 + h a_1 + h^2 a_2 + \cdots . $$
to mean that for all $N$, $ a = a_0 + h a_1 + \cdots + h^{N-1}
a_{N-1} + h^N r_N (h) $ with $a_j \in S^{m,m^{\prime}} $,
independent of
$h$,  and $r_N (h ) $
bounded family of $ S^{m,m^{\prime}} $.




\medskip

We now give  estimates on $ O \! p_h (a) $ in Schatten classes.

\begin{prop} \label{Schattenlocal} Let $\chi$ be a bounded function supported in $ \Ra^+ \times
\Omega$, with $ \Omega \Subset \Ra^{n-1} $ and $\nu > 0$, $ q \geq
1$ be positive real numbers such that $ \nu > (n-1)/ q $. 
Then for any $ a \in S^{-\varepsilon - 1 / q , - \varepsilon -
(n-1)/q} $, with $ \varepsilon > 0 $, the operator $ e^{- \nu r}
\chi O \! p_h (a) $ belongs to $ \Schn{q} $ and
\begin{eqnarray}
 \left| \left| e^{- \nu r} \chi O \! p_h (a)  \right|
\right|_{\bf q} \leq C h^{-n/q} , \qquad h \in (0,1] ,
\label{Schatten0}
\end{eqnarray}
 where $ C
$ depends on finitely many  semi-norms of $a$. In particular,
if $\nu
> (n-1)$, then $ e^{- \nu r} \chi O \! p_h (a) $ is trace class
and
\begin{eqnarray}
 \emph{tr} \left(  e^{- \nu r} \chi O \! p_h (a)  \right) = (2
\pi h)^{-n} \int \! \! \! \int \! \! \! \int \! \! \! \int e^{-
\nu r} \chi (r,y)  a (r,y,\rho, \eta) \emph{d} \rho \emph{d} \eta
\emph{d} r \emph{d}y . \label{trace0}
\end{eqnarray}
\end{prop}

\medskip

\noindent {\it Proof.} With no loss of generality we can assume
that $ \varepsilon $ is a small as we want and in particular that
\begin{eqnarray}
 \nu > 2 \varepsilon + \frac{n-1}{q} . \label{nu}
\end{eqnarray}
  We can also assume that
$\chi(r,y) = \phi (r) \psi (y)$ with $ \psi \in
C_0^{\infty}(\Omega)$ and $ \phi \equiv 1$ near infinity. The
estimate $ \reff{Schatten0} $ will follow from the fact that we
can write $ e^{- \nu r} \chi O \! p_h (a) = A_h B_h $ with
\begin{eqnarray}
\ || A_h ||_{\bf q} \leq C h^{-n / q} , \qquad || B_h ||_{\infty}
\leq C , \qquad h \in (0,1] . \label{reducSchatten}
\end{eqnarray}
In order to construct $A_h$, we pick $ \widetilde{\psi} \in
C_0^{\infty} $ such that $\widetilde{\psi} \psi = \psi$, and we
set $A_h = A_h^1 \otimes A_h^2$ with
$$ A_h^1 = e^{- \varepsilon r} \phi (r) \scal{h D_r}^{-\varepsilon - 1 /q} , \qquad
A_h^2 = \widetilde{\psi}(y) \scal{h D_y}^{- \varepsilon - (n-1) /
q  } . $$ By standard estimates, we know
that $ A_h^1 = {\mathcal O} (h^{-1/q}) $ in $
\Schn{q}(L^2(\Ra,\mbox{d}r)) $ and that $ A_h^2 = {\mathcal O}
(h^{-(n-1)/q}) $ in $ \Schn{q}(L^2(\Ra^{n-1},\mbox{d}y)) $, thus
the first estimate of $ \reff{reducSchatten} $ holds. We now have
to consider $B_h$ which we define, using $\widetilde{\phi}$ such
that $ \widetilde{\phi} \phi = \phi$, by
$$ B_h = \scal{h D_r}^{\varepsilon + 1 / q} \scal{h D_y}^{\varepsilon + (n-1) q}
 e^{-(\nu - \varepsilon)r} \widetilde{\phi} (r)  \psi (y) O \! p_h
 (a) . $$
 In order to show
 the second estimate of $ \reff{reducSchatten} $, it is
 enough to show that $B_h = O \! p_h (b)$ with $b  = b (h)$
  bounded in $S^{0,0}$. We first consider $  \scal{h D_y}^{\varepsilon + (n-1) q}
 e^{-(\nu - \varepsilon)r} \widetilde{\phi}(r) \psi (y) O \! p_h
 (a)  = O \! p_h (b_1)  $ with
\begin{eqnarray}
  b_1 (r,y,h \rho,h\eta,h)=  e^{-(\nu-\varepsilon)r} \widetilde{\phi}(r) (2
 \pi)^{1-n} \!
\int e^{i y.\xi} \scal{h \eta + h \xi}^{\varepsilon + (n-1)/q}
\widehat{a}_{\psi} (r,\xi,h \rho , h  \eta) \mbox{d} \xi
\label{integrand}
\end{eqnarray}
 where $
\widehat{a}_{\psi} $ is the Fourier transform of
$\psi(y)a(r,y,\rho,\eta)$ with respect to $y$. By Peetre's
inequality we have $ \scal{h \eta + h \xi}^{\varepsilon + (n-1)/q}
\leq C \scal{h \eta }^{\varepsilon + (n-1)/q} \scal{ h
\xi}^{\varepsilon + (n-1)/q} $ and the integrand of $
\reff{integrand} $ is dominated, for any $M$,  by
$$ C \scal{\xi}^{-M} e^{- (\nu - \varepsilon )r} \scal{h \eta }^{\varepsilon + (n-1)/q}
 \scal{h \rho}^{-\varepsilon - 1 / q} \scal{h e^{-r} \eta}^{-\varepsilon - (n-1) / q} \ \leq \ \widetilde{C}
 \scal{\xi}^{-M} \scal{h \rho}^{-\varepsilon - 1 / q} , $$
 where we used $ \reff{nu} $.
The same holds for the derivatives of $b_1$ and shows this $ b_1
\in S^{-\varepsilon - 1 /q , 0} $. Furthermore it depends
continuously on $b$. Thus $  \scal{h D_r}^{\varepsilon + 1 / q } O
\! p_h (b_1) = O \! p_h (b) $ with $b$ bounded and the second
estimate of $ \reff{reducSchatten} $  holds.
The proof of $ \reff{trace0} $ is well known.
\finpreuve

\medskip

This proposition is the key of  theorem $ \refe{theo1} $ . For
 theorem $ \refe{theo2} $ we shall need the slightly stronger estimate
 $ \reff{Schatten1} $ below,
 in order to use $ \reff{propagationestimate} $.
By standard results on pseudo-differential calculus
we know that, for any $ M \in \Ra$ and $ b \in S^{ 0 , 0} $
\begin{eqnarray}
 \scal{r}^M O \! p_h (b) \scal{r}^{-M} = O \! p_h ( \widetilde{b}
) , \qquad \widetilde{b} \in S^{0,0} . \label{polynomial}
\end{eqnarray}
 Furthermore, for all $ M $,
$ \scal{r}^M  e^{- \varepsilon r} \phi (r) \scal{h
D_r}^{-\varepsilon - 1 /q}  \scal{r}^M $ still belongs to $
\Schn{q} $ with norm $ {\mathcal O}(h^{-n/q}) $. Using this remark and
the proof of proposition $ \refe{Schattenlocal} $, we see easily
that
\begin{eqnarray}
 \left| \left| e^{- \nu r} \scal{r}^M \chi O \! p_h (a) \scal{r}^M  \right|
\right|_{\bf q} \leq C h^{-n/q} , \qquad h \in (0,1] ,
\label{Schatten1}
\end{eqnarray}
for any $M \in \Ra$, $ a \in S^{- \varepsilon - 1 / q , -
\varepsilon - (n-1)/q} $ and $ \nu > (n-1)/q $. Furthermore, if $a^{(\kappa)}$
is a family of symbols which is bounded in $ S^{- \varepsilon - 1 / q , -
\varepsilon - (n-1)/q} $ such that $ a^{(\kappa)} \rightarrow a^{(0)} $
in $C^{\infty}(\Ra^{2n})$ (or ${\mathcal D}^{\prime}(\Ra^{2n})$), as
$ \kappa \downarrow 0 $ and such that $\rho^2 + e^{-2r}|\eta|^2$ is
  bounded, independently of $ \kappa $, on their support then
  \begin{eqnarray}
  e^{- \nu r} \scal{r}^M \chi O \! p_h (a^{(\kappa)}) \scal{r}^M \rightarrow
  e^{- \nu r} \scal{r}^M \chi O \! p_h (a^{(0)}) \scal{r}^M
   \qquad \mbox{in} \ \Schn{q}, \qquad \kappa \downarrow 0 .
\label{Schatten1parametre}
\end{eqnarray}
This follows simply from the fact that $e^{- \delta r} a^{(\kappa)}
\rightarrow e^{- \delta r} a^{(0)} $ in  $ S^{- \varepsilon - 1 / q , -
\varepsilon - (n-1)/q} $ for any $\delta > 0$, and from the fact that we can
choose $ \delta $ small enough such that $ \nu - \delta > (n-1)/q $.

\bigskip

Now we are going to study symbols and operator depending smoothly
on $\e \in [0,1]$.
\begin{defi} The class $S^{\nu,m,m^{\prime}}_{\e}$ is the set of
functions $a_{\e}(r,y,\rho,\eta)$ which are smooth w.r.t.
$r,y,\rho,\eta$ and $\e$ such that for all $R > 0$, $ \Omega
\Subset \Ra^{n-1} $ and $j,k,l,\alpha,\beta$ there exists $C$
satisfying
$$ \left| \partial_{\e}^j
\partial_r^k \partial_{\rho}^l \partial_y^{\alpha} \partial_{\eta}^{\beta} a_{\e} (r,y,\rho,\eta) \right|
\leq C  e^{-(j+\nu)r} \scal{\rho}^{m} \scal{e^{-r}\eta}^{m^{\prime}} ,
$$
for all $ r \geq R , \ y \in \Omega , \ \rho \in \Ra , \ \eta \in
\Ra^{n-1}$. We set $ S^{\nu , - \infty}_{\e} = \cap_{ m,m^{\prime}}
S^{\nu,m,m^{\prime}}_{\e} $.
\end{defi}

 The typical example of
such symbols is given by $ a_{\e}(r,y,\rho,\eta) = f (\rho^2 +
e^{-2r} g_{\e}(r,y,\eta)) $, which belongs to
$a_{\e} \in S^{0,-\infty}_{\e}$ if $f \in {\mathcal S}$. If we replace $a$
by $ a_{\e} \in  S^{0,-\infty}_{\e} $ then $\reff{transportcutoff} $
holds with $  S^{0,-\infty}_{\e} $ instead of
$ S^{-\infty} $.
This class is particularly natural and convenient since proposition
$ \refe{Schattenlocal} $ shows that for all
$a_{\e} \in S^{0, - \varepsilon - j / q ,- \varepsilon - (n-1) j / q }_{\e}$
\begin{eqnarray}
 \left| \left| \partial_{\e}^j \chi O \! p_h (a_{\e}) \right|
\right|_{\bf q/j}
\leq C h^{-n j / q}, \qquad \mbox{provided} \ \ \ q > n -1, \ \
 q \geq j \geq 1 .
\end{eqnarray}


\medskip

The main drawback of these classes is the following. If $a_{\e} \in
S^{\nu,m,m^{\prime}}_{\e}$ then we don't have in general $ O \! p_h
(a_{\e})^* = O \! p_h (\widetilde{a}_{\e})$ with $
\widetilde{a}_{\e} \in S^{\nu,m,m^{\prime}}_{\e}$, which is due to the
fact that pseudo-differential operators do not preserve
exponential decay.
This problem can be overcome by considering properly supported
operators. This is the purpose of what follows.

 Recall that the Schwartz kernel of $ O \! p_h (a_{\e})
$ is given by the following oscillatory integral 
$$ {\mathcal K}_{h,\e} (r,r^{\prime},y,y^{\prime}) = (2 \pi h)^{-n}
\int \! \! \! \int
e^{i(r-r^{\prime})\rho/h + i (y-y^{\prime}).\eta/h} a_{\e}(r,y,
\rho ,   \eta) \ \mbox{d} \rho \ \mbox{d} \eta $$ 
which we can write, using $ \theta \in C_0^{\infty}(\Ra)$ and $
\theta \equiv 1 $ on $(-\delta,\delta) $, as
$ {\mathcal K}_{h,\e} = {\mathcal K}_{h,\e}^{diag} +
{\mathcal K}_{h,\e}^{o \! f \! f} $ with
\begin{eqnarray}
 {\mathcal K}_{h,\e}^{diag}(r,r^{\prime},y,y^{\prime}) =   
{\mathcal K}_{h,\e}
(r,r^{\prime},y,y^{\prime}) \theta (r-r^{\prime}) \theta
(|y-y^{\prime}|) . \label{cutoffdiagonale} 
\end{eqnarray}
\begin{defi} If $ a_{\e} \in S^{\nu,m,m^{\prime}}_{\e} $, we define
$O \! p_h^{pr}(a_{\e})$ as the operator with Schwartz kernel  $
{\mathcal K}_{h,\e}^{diag} $.
\end{defi}

The following proposition explains in which sense
$  O \! p_h^{pr}(a_{\e}) - O \!
p_h(a_{\e}) $ is small.
\begin{prop} \label{restepropre} Assume that $\nu \geq 0$, 
$m \leq 0$ and $ m^{\prime} \leq 0$. Let
$\chi$ be a bounded function supported in $ \Ra^+ \times \Omega $,
with $ \Omega \Subset \Ra^{n-1} $. Let $K_{h,\e}$ the
operator with kernel $\chi (r,y)  {\mathcal K}_{h,\e}^{o \! f \!
f}(r,r^{\prime},y,y^{\prime})$. Then for all $s \in \Ra$, $j \geq
0 $ and $N \geq 0$, there exists $C > 0$ such that
\begin{eqnarray}
 \left| \left|
e^{(j+\nu)r} \scal{r}^{-s} \frac{d^j}{d \e^j}   K_{h,\e}
\scal{r}^{s} \right| \right|_{\infty}  \leq C h^{N}, \qquad h \in
(0,1], \ \ \e \in [0,1] . \label{off}
\end{eqnarray}
\end{prop}

\medskip

\noindent {\it Proof.} It is standard. On the support of $
{\mathcal K}_{h,\e}^{o \! f \! f}  $, we have either $|r-r^{\prime}| \geq
1 $ or $|y-y^{\prime}| \geq  1 $ and we can do as many
integrations by parts as we want with $|r^{\prime}-r|^{-1} h
\partial_{\rho}$ or $|y^{\prime}-y|^{-2} h^2 \Delta_{\eta}$. This shows that $
\partial_{\e}^j {\mathcal K}_{h,\e}^{o \! f \! f}
(r,r^{\prime},y,y^{\prime}) $ is a linear combination of
\begin{eqnarray}
  h^N (2 \pi h)^{-n}  e^{-jr} \int \! \! \! \int
e^{i(r-r^{\prime})\rho / h + i (y-y^{\prime}). \eta / h}
 b_{N,\e} (r,r^{\prime},y,y^{\prime},h\rho,h\eta) \ \mbox{d} \rho
 \mbox{d}\eta \label{noyau}
\end{eqnarray}
 with $b_{N,\e}$ a product of derivatives of $\partial_{\e}^j a_{\e}$, $\chi$
 and $|r-r^{\prime}|^{-1}$, $|y-y^{\prime}|$.   The result follows from $
 \reff{CalderonVaillancourt} $ and $ \reff{polynomialstable} $. \finpreuve

\medskip

Note that the adjoint $  K_{h,\e}^* $ satisfies $
\reff{off} $ with $e^{(j+\nu)r}$ on the right. This leads to the
following definition.
\begin{defi} \label{negligeable}
 A family of bounded operators $A_{\e,h}$ is said $ h^N $
negligible  if for all $j \geq 0$ and $M \in \Ra$, we
can write $ \partial_{\e}^j A_{\e,h} = A^{0}_{\e,h} + A^{1}_{\e,h} + \cdots +
 A^{j}_{\e,h} $ where, for $ 0 \leq k \leq j $,
$$  \left| \left|
e^{\nu_{j-k}r} \scal{r}^{-M}   A^{k}_{\e,h} e^{\nu_k r} \scal{r}^{M}
\right| \right|_{\infty}  \leq C_{N,M,j} h^N , \qquad h \in (0,1] , \ \ \e
\in [0,1] , $$
for some pair $ \nu_{k}, \nu_{j-k} \geq 0  $ such that 
$ \nu_j + \nu_{j-k} = \nu + j $.
It is called negligible if it is $ h^N $ negligible for all $ N $.
\end{defi}
Proposition $ \reff{restepropre} $ shows that $ \chi O \!
p_h^{pr}(a_{\e}) -  \chi O \!p_h (a_{\e})  $ is negligible,
provided $m \leq 0$ and $m^{\prime} \leq 0$. If $A_{\e,h}$ is
negligible, then $A_{\e,h}^*$ is clearly negligible as well, hence all this shows
that, by using standard methods for the calculus of the adjoint
of a pseudo-differential, we obtain easily the following result.

\begin{prop} \label{adjoint} Let  $a_{\e} \in S^{\nu,m,m^{\prime}}_{\e} $ be supported in $(0,\infty)
\times \Omega$ with $\Omega$ bounded and $m,m^{\prime} \leq 0$.
Then $ O \! p_h (a_{\e})^*$ is the sum of $ O \! p_h (
\widetilde{a}_{\e}) $ and of a negligible operator, with  
$\widetilde{a}_{\e} \in S^{\nu,m,m^{\prime}}_{\e} $ such that
$$ \widetilde{a}_{\e,j} \sim \sum_j h^j  \sum_{|\alpha|+k =j} \frac{1}{k!} 
\frac{1}{\alpha!} \partial_r^k D_{\rho}^k
\partial_y^{\alpha} D_{\eta}^{\alpha} \overline{a_{\e}} . $$
\end{prop}

We end this subsection with a simple remark. In practice, we will use
the
negligibility as follows: if $A_{\e,h}$ is $h^N$
negligible, we can write for all $j \geq 1$, $M \geq 0$ and $ \varepsilon > 0 $
as
\begin{eqnarray}
\partial_{\e}^j A_{\e,h} = h^N \sum_{k=0}^j e^{- \nu_k r} \scal{r}^{-M}
B_{\e,h}^k \scal{r}^{-M} e^{- \nu_{j - k}} \label{astuceformelle}
\end{eqnarray}
where the operators $ B^{k}_{\e,h}$ are bounded in operator norm  and 
 $ \nu_k + \nu_{j-k} = \nu + j -\varepsilon $,  
with $ \nu_k , \nu_{j-k} \geq 0 $ for all $ k $.

\subsection{Global theory and functional calculus}
The pseudo-differential operators that we are going to consider on
$X$ will be of the form
\begin{eqnarray}
A = \sum_{k \in I} \widetilde{\chi}_k O \! p_h (a^k) \chi_k,
\qquad I \ \mbox{finite}  \label{formepseudo}
\end{eqnarray}
For each $k$, $\chi_k$ and $\widetilde{\chi}_k$ are supported in
the same chart which is either relatively compact  or a
chart at infinity, and $a^{k}$ is a symbol expressed in the
coordinates associated to the chart.
If we work in a
chart at infinity, we will always use the radial variable $r$ and
variables $y_1, \cdots , y_{n-1} $ associated to the manifold at
infinity $Y$. 
 

Let us recall some (standard) abuse of notations which are
convenient. We use the same notation for $ \widetilde{\chi}_k O \!
p_h ( a_{k} ) \chi_{k} $ as an operator acting on $L^2(\Ra^n)$ or
on $L^2(X)$. Furthermore if, for each $k \in I$,
$\widetilde{\chi}_k O \! p_h ( a^{k} )
\chi_{k} $ is bounded on $L^2(\Ra^n)$ with the Lebesgue measure
associated to the corresponding coordinates, then $A$ is bounded
on $L^2(X)$ and its norm can be estimated by the sum  of norms of $
\widetilde{\chi}_k O \! p_h ( a^{k} ) \chi_{k}  $ on $L^2(\Ra^n)$.
This is due to our choice of the density $ \mbox{d}vol
$.
The same remark holds for estimates in $ \Schn{q} (L^2(X)) $ and
we will therefore use the notations $||.||_{\infty}$ and
$||.||_{\bf q}$, initially used on $ L^2 ( \Ra^n  )$,
 for the respective norms of bounded operators and
Schatten classes relative to $L^2(X)$.

\bigskip

We now apply the results of the previous subsection to the analysis of functions of
$H_{\e}$. We will only consider Schwartz functions $f$ and
 use the
 Helffer-Sj\"ostrand formula,
\begin{eqnarray}
  f (H_{\e}) = \frac{1}{2 \pi} \int \! \! \! \int_{\Ra^2}
\partial_{\bar{z}} \widetilde{f} (s+it)
 R_{\e}(s+i t) \ \mbox{d}s \ \mbox{d} t \label{HelfferSjostrand}
\end{eqnarray}
where $ \partial_{\bar{z}} = \partial_s + i \partial_t $
and a $ \widetilde{f} $ is a quasi-analytic extension of $f$. i.e. a
$ C^{\infty} $ function on $ \Ca $, supported in the strip $|t| \leq 1$ and
such that
$$   \widetilde{f}_{|{\mathbb R}} = f , \qquad
\left| \partial_{\bar{z}} \widetilde{f} (s+it) \right| \leq C_{M}
|t|^M \scal{s}^{-M}  $$
for all $s,t \in \Ra$ and $M \geq 0$.
We do not insist on the construction of $ \widetilde{f} $ which can be chosen
depending continuously on  $ f $  and refer for instance to \cite{DiSj1} for the 
details.

The formula $ \reff{HelfferSjostrand} $ shows that we only have to study the
pseudo-differential expansion of $R_{\e}(z)$.      We recall that
the pseudo-differential analysis of the resolvent is
well known for operators on compact
manifolds (or for elliptic
operators on $\Ra^n$)  and thus we will only focus on the calculations
 in charts at infinity.
We look for a parametrix $ Q_{\e}(z) $ of $R_{\e}(z)$ of the form
$\reff{formepseudo}$ and more precisely
$$  Q_{\e}(z) = \sum_{k \in I} \widetilde{\chi}_{k} O \! p_h (q^k_{\e}(z)) \chi_k $$
where $ \sum_{k \in I} \chi_k = 1 $ is an admissible partition of unit,
 $   \widetilde{\chi}_k = 1 $ near the support of $ \chi_k $, and
$q^k_{\e}(z)$ is an admissible symbol for each $k$. We explain the construction
of these symbols in a single chart at infinity and drop the index $ k $ for
convenience. Since we want to get $ (H_{\e}-z) Q_{\e}(z) = 1 + {\mathcal
O}(h^{\infty}) $,
 we seek $  q_{\e} (z) \sim \sum_j h^j q_{j,\e} (z)  $
 satisfying, in the chart that we consider,
\begin{eqnarray}
 \sum_{j,l,m}  h^{j+l+m} \left( p^z_{l,\e} \# q_{j,\e} (z)  \right)_m \sim 1 .
 \label{identification}
\end{eqnarray}
Here  $p^z_{l,\e}$ are the symbols of $ H_{\e} - z $, i.e.
$ p_{0,\e}^z = p_{\e}
-z $ and $ p_{l,\e} = p_{\e}^{(l)} $ for $ l = 1,2 $,
and the notation $ ( a \#  b)_m$ stands for the $ m $-th symbol of the (finite)
expansion of the product $ O \! p_h (a) O \! p_h (b)  $ if $ O \! p_h (a) $
is a differential operator. The condition $ \reff{identification} $ yields
\begin{eqnarray}
   q_{0,\e} (z) & = &  (\rho^2 +
g_{\e}(r,y,e^{-r}\eta) - z )^{-1} \\
q_{j,\e} (z) & = & - q_{0,\e} \sum_{j_0 + j_1 + j_2 = j , \
j_1 < j} \left( p^{(j_0)}_{\e} \# q_{j_1,\e} \right)_{j_2} , \qquad j
\geq 1 .
\end{eqnarray}
This procedure is the standard one used by Seeley \cite{Seel1} and 
Helffer-Robert \cite{HeRo1} but
the point that we want to make here is the following: since the symbols
of $H_{\e}$ are of the form $ \reff{formprincipal} $, it follows
clearly  that $q_{j,\e} ( z ) $ defined as above is of
the form $ \reff{formprincipal} $ as well. By an easy induction,
we get, for $j \geq 1$,
\begin{eqnarray}
  q_{j,\e} = \sum_{l}^{2j-1} d_{l,j,\e} (p_{\e}-z)^{-1-l},
  \label{formsymbol0}
\end{eqnarray}
where $d_{l,k,\e}$ is polynomial w.r.t. the variables $\rho$ and $ e^{-r}\eta$
which is
independent of $z$ and a linear combination of elements of
$S_{\e}^{0,m,m^{\prime}}$ with $m + m^{\prime} \leq 2 l - j$. In
particular we can write $q_{j,\e}(r,y,\rho,\eta)  =
\widetilde{q}_{j,\e}(r,y,\rho,e^{-r}\eta) $ for some $
\widetilde{q}_{j,\e} \in S_{0,1}^{-1 - j/2,-1-j/2}$ and we have
$$ q_{j,\e} \in S_{\e}^{0,m,m^{\prime}} , \qquad \mbox{for all }
 \ m , m^{\prime} \leq 0   \ \ \mbox{such that}
\ \ m + m^{\prime} = - 2 - j ,  .  $$
These remarks, combined with the fact that $ |p (p-z)^{-1}| \leq C
\scal{\re z }/ |\im z| $, for $p \in \Ra$ and $z \notin \Ra$ show that the 
semi-norms of $ q_{j,\e} $ in
$S^{0,m,m^{\prime}}_{\e}$ ($m+m^{\prime} = - 2 -j$) are dominated by
\begin{eqnarray}
 C  \frac{\scal{\re z}^M}{|\im z|^M} , \label{estimeseminorme}
\end{eqnarray}
  for some $C,M$
depending on the semi-norm.

We can now construct the global parametrix. We fix $N \geq
0 $ and define it as
$$ Q_{N,\e}(z) = \sum_{k \in I} \widetilde{\chi}_k  \left( \sum_{j \leq N}
 h^j O \! p_h (q_{j,\e}^k(z)) \right) \chi_k , $$
 with the $ q_{j,\e}^k(z) $ defined in each chart by the preceeding procedure.
Then we have
\begin{eqnarray}
 ( H_{\e} - z ) Q_{N,\e}(z) & = &  \sum_{k \in I} \widetilde{\chi}_k
(H_{\e}-z)  \left( \sum_{j \leq N} h^j O \! p_h (q_{j,\e}^k(z))
\right) \chi_k \nonumber \\
& &
 +  \sum_{k \in I} \left[ H_{\e} , \widetilde{\chi}_k
\right]  \left( \sum_{j \leq N} h^j O \! p_h (q_{j,\e}^k(z))
\right) \chi_k . \label{globalization}
\end{eqnarray}
\begin{lemm} \label{resteSchwartz} For each $k$ and $j$, the Schwartz kernel
of  $ [ H_{\e} , \widetilde{\chi}_k ]   O \! p_h (q_{j,\e}^k(z))
\chi_k $ is $ {\mathcal O}(h^{\infty}) $ in the Schwartz space
${\mathcal S} (\Ra^{2 n})$. More precisely, if we note it
$ {\mathcal K} (r,r^{\prime},y,y^{\prime},\e,z,h)$ then we have
$$  \left| \partial^{\gamma}_{r,r^{\prime},y,y^{\prime},\e}
{\mathcal K}  (r,r^{\prime},y,y^{\prime},\e,z,h) \right| \leq C h^M e^{-M r}
\scal{r^{\prime}}^{-M}
\scal{y}^{-M}  \scal{y^{\prime}}^{-M}   \frac{ \scal{
\emph{Re} \ z }^M}{| \emph{Im} \ z|^M}  $$ for all
 multiindex $ \gamma $, $M \geq 0$, $r,r^{\prime} \geq 0 $, $ y, y^{\prime} \in
 \Ra^{n-1} $, $ \e \in [0,1] $, $ h \in [0,1] $ and $ z \notin \Ra $.
\end{lemm}

\medskip

\noindent {\it Proof.} It follows easily by integrations by parts
similar to those of lemma $ \refe{restepropre} $. The exponential
decay is due to the following fact.  We can choose
$\widetilde{\chi}_k = \widetilde{\phi}_k (r) \widetilde{\psi}_k
(y)$, with  $\widetilde{\phi}_k \equiv 1$ near infinity. Then
$[H_{\e},\widetilde{\chi}_k]$ is either compactly supported w.r.t.
$r$, or $|y-y^{\prime}|
\ne 0$ in which case we integrate by part using $\Delta_{\eta}$
and get as many powers of $e^{-r}$ as we want. We omit the
other details. \finpreuve

\medskip

This lemma shows  in particular that that the second term of the right hand side of $
\reff{globalization} $ is negligible, in the sense of definition $ \refe{negligeable} $.
On the other hand, by construction, the first sum of the right hand side of
$ \reff{globalization} $ is
$$ \sum_{k \in I}  \widetilde{\chi}_k   \left( 1 +  h^N O \! p_h (\varrho_{N,\e}^k(z)) \right)
\chi_k  = 1 + h^N  \sum_{k \in I}  \widetilde{\chi}_k   O \! p_h (\varrho_{N,\e}^k(z))
\chi_k
 $$ where  $ \varrho_{N,\e}^k(z) =
\sum_{j+l+ m \geq  N} h^{j+l + m }  ( p^z_{l,\e} \# q_{j,\e}^k )_m $.
 It is  not hard to check that it belongs to
$S^{0,-N/2,-N/2}_{\e}$ (for instance) using the form of $ q_{j,\e}^k $.
 We can summarize our result as
follows.

\begin{prop} \label{domaine} The operators $H_{\e}$ are
essentially self-adjoint from $C_0^{\infty}(X)$ and the domain of their
self-adjoint realizations is independent on $ \e \in [0,1] $.
We have for all $N$
\begin{eqnarray}
R_{\e}(z) = Q_{N,\e}(z) + h^N R_{\e}(z) {\mathcal R}_{N,\e}(z) .
\label{basereso}
\end{eqnarray}
Here $h^N {\mathcal R}_{N,\e}(z) $ is the sum of $ \sum_{k \in I}
h^N \widetilde{\chi}_k O \! p_h ( \varrho_{N,\e}^k(z) ) \chi_k  $
and of the operators with kernel studied in lemma $
\refe{resteSchwartz}$.
\end{prop}

\medskip

\noindent {\it Proof.} By the standard trick (see \cite{Robe1}), we
see that $(H_{\e} \pm i)^*$ is injective for $h$ small enough
since $ (H_{\e} \pm i) Q_{N,\e}(\pm i) = 1 + {\mathcal O}(h^N) $
 in operator norm on $L^2(X)$. This
implies the existence of a unique self-adjoint realization.
Hence the resolvent $R_{\e}(z)$ is well defined and we can apply it to the left
of $ \reff{globalization} $ which yields $ \reff{basereso} $. In
particular, for $h$ small enough and fixed $z$ we have
$$ R_{\e}(z) = Q_{N,\e}(z) \left(1 - h^N {\mathcal R}_{N,\e}(z)  \right)^{-1}  $$
and since $H_{\e_1} Q_{N,\e_2}(z)$ is bounded for all $\e_1,\e_2$,
the same holds for $H_{\e_1} R_{\e_2}(z)$. This proves the
independence of the domain w.r.t. $\e$. \finpreuve

\medskip

Let us give a first application of this proposition.

\begin{lemm} \label{aprioriSchatten0} Let $A$ be a differential operator on
$X$ of order
$m$ of the following form in any chart at infinity
$$ A = \sum_{l+|\alpha| \leq m} a_{l,\alpha}(r,y)  (e^{-r}D_y)^{\alpha} D_r^l $$
with $a_{l,\alpha}$ bounded. Then for all $  \nu,q,k $ such that $  \nu > (n-1)/q
$ and  $2 k  - m > n /q$ , there exists
$C$ such that
$$ \left| \left| e^{-\nu r} h^m \! A R_{\e}(i)^{1+k} \right| \right|_{\bf q} \leq C h^{- n / q} , \qquad
h \in ( 0 ,1 ] , \ \e \in [0,1] . $$ If $\nu = 0$, by convention
the above norm is the operator norm $||.||_{\infty}$ and $ h^{n/q}
= 1 $.
\end{lemm}

Note that this lemma implies the following estimate, for $ \im z
\ne 0$,
$$ \left| \left| e^{-\nu r} h^m \! A R_{\e}(z)^{1+k} \right| \right|_{\bf q}
 \leq C h^{- n / q}   \frac{ \scal{\mbox{Re} \ z }^{1+k}}{| \mbox{Im} \ z|^{1+k}}
 , \qquad
h \in ( 0 ,1 ] , \ \e \in [0,1] . $$


\noindent {\it Proof.} We use the fact that $ k ! R_{\e}(z)^{k+1}
=
\partial_z^k R_{\e}(z) $ and Leibnitz rule into $ \reff{basereso}
$ to obtain
\begin{eqnarray}
R_{\e}(i)^{k+1}  = \frac{1}{k!} \left( \partial^{k}_z Q_{N,\e} \right) (i)  \left(1 - h^N
 \sum_{\mu \leq k+1} ( H_{\e} - i)^{\mu} \partial_{z}^{\mu}
{\mathcal R}_{N,\e}( i )
 \frac{1}{\mu ! (k - \mu) ! } \right)^{-1}  \label{Neumann}
\end{eqnarray}
for $ h  \leq h_0$ small enough. Note that $ ( H_{\e} - i)^{\mu}
\partial_{z}^{\mu} {\mathcal R}_{N,\e}( i ) $ is bounded for all $ \mu
$. Using the form of the symbols of $ Q_{N,\e} $ given by $
\reff{formsymbol0} $  and by proposition $
\refe{Schattenlocal} $ applied to $ e^{-\nu r} h^m A
\partial_z^k Q_{N,\e}(i) $, we obtain the result for $h \leq h_0$.
Note that we use the fact that $ k ! \partial_z^{k} q_{0,\e}(z) = q_{0,\e}(z)^{k+1} $
belongs to $  S_{\e}^{ 0, - j - \varepsilon - 1 / q , - |\alpha |-\varepsilon -(n-1) / q }$
for all $ j , |\alpha| $ such that $ j + |\alpha|=m $ and $ \varepsilon > 0 $
small enough.
The result for $ h \leq 1$ follows by writing $R_{\e}(i) =
(\bar{h}^2 P_{\e}+i)^{-1}(\bar{h}^2 P_{\e}+i) R_{\e}(i)$ with $
\bar{h} = \min (h_0 , h) $. \finpreuve

\medskip

This proposition will be used in order to estimate  $
\partial_{\e}^j ( R_{\e}(z) {\mathcal R}_{N,\e}(z) ) $ in Schatten
classes. We shall also need estimates similar to $
\reff{Schatten1} $. To that end, we first recall that for
any $ 0 \leq M_0 \leq 1 $
$$ \left[ \scal{r}^{M_0} , R_{\e}(z) \right] = - R_{\e}(z) \left[ \scal{r}^{M_0} ,
H_{\e} \right] R_{\e}(z) $$
where $ \left[ \scal{r}^{M_0} , H_{\e} \right] $ is a differential
operator of order $1$ which is $H_{\e}$ bounded since $
[\partial_r , \scal{r}^{M_0}] $ is bounded. This implies in
particular that
$$ \scal{r}^{M_0} R_{\e}(z) \scal{r}^{-M_0} =
R_{\e} (z) - R_{\e} (z) \left[ \scal{r}^{M_0} , H_{\e} \right]
R_{\e}(z) \scal{r}^{- M_0} . $$ Since any $M \geq 0 $ can be
written $l M_0$ with $l \in \Na $, this formula can be iterated to
show that
$$
 \left| \left| (H_{\e}+i) \scal{r}^M R_{\e}(z) \scal{r}^{-M} \right|
\right|_{\infty} \leq C  \frac{ \scal{\re z}^{M^{\prime}}}{| \im z |^{M^{\prime}}}
 , $$  for some $ M^{\prime} \geq 0 $,
 uniformly w.r.t. $h \in (0,1]$ and $ \e \in [0,1]$.
More generally, we can obtain rather easily for any $k \in \Na$
and $M \geq 0$
\begin{eqnarray}
\left| \left| (H_{\e}+i)^{1+k} \scal{r}^M R_{\e}(z)^{1+k}
\scal{r}^{-M} \right| \right|_{\infty} \leq C  \frac{
\scal{\re z}^{M^{\prime}}}{| \im z |^{M^{\prime}}}  . \label{raccourci}
\end{eqnarray}
for some $ C $ and $ M^{\prime} $ independent of $ \e  $ and $ h $.
On the other hand,
$ \partial_{\e}^j R_{\e} (z) = (-1)^j j  ! \ R_{\e}(z) ( V R_{\e}(z) )^j $
 can be written for any $ M \in \Ra$ as
$$ \partial_{\e}^j R_{\e} (z) =  ( - 1 )^j j ! \  \scal{r}^{-M} \left(
\scal{r}^M R_{\e}(z) \scal{r}^{-M} \right) \left( \widetilde{V}
\scal{r}^{-M} R_{\e}(z) \scal{r}^M \right)^j  \scal{r}^{-M}  $$
where the differential operator $ \widetilde{V}  = \scal{r}^M V \scal{r}^{M} $
is of the form $e^{- \nu r} h^2 A$ for any $\nu < 1$ with the
notations of lemma $\refe{aprioriSchatten0}$.
 Using these remarks combined with
lemma $ \refe{aprioriSchatten0} $, we will obtain
\begin{lemm} \label{restemoitie} For any $ q > n - 1 $, $M \geq 0$ real and $ j \geq 1$ integer,
there exists $C,M^{\prime} $ such that
$$  \left| \left| \scal{r}^M \left( \partial_{\e}^j R_{\e}(z) \right) R_{\e}(i)^{j-1}
 \scal{r}^M \right| \right|_{\bf q/j} \leq C h^{- n j / q}  \frac{
 \scal{ \emph{Re} \ z }^{M^{\prime}} }{| \emph{Im} \ z |^{M^{\prime}}}  $$
 for all $ h \in (0,1]$ and $ \e \in [0,1] $.
\end{lemm}

\medskip

\noindent {\it Proof.} We proceed by induction on $j$. If $ j = 1
$, we write
$$ \scal{r}^M  \partial_{\e} R_{\e}(z)
 \scal{r}^M = \left( \scal{r}^M R_{\e} (z) \scal{r}^{-M} e^{-\nu r} \right)
 \left( h^2 A  \scal{r}^{-M}  R_{\e}(z) \scal{r}^M \right) $$
with $ (n-1)/q < \nu < 1 $ and $ e^{- \nu r} h^2 A = \scal{r}^M V \scal{r}^M
$. The first factor belongs to $\Schn{q}$ whereas the second one
is bounded, by lemma $ \refe{aprioriSchatten0} $. If $j \geq 2$,
we have $ \partial_{\e}^j R_{\e} (z) = - j (\partial_{\e}^{j-1}
R_{\e}(z)) V R_{\e}(z) $ and then the operator that we want to
estimate can be written
$$ - j \left( \scal{r}^M \left( \partial_{\e}^{j-1} R_{\e}(z) \right) R_{\e}(i)^{j-2}
 \scal{r}^M \right) \left(  \scal{r}^{-M} (H_{\e}+i)^{j-2} V R_{\e}(z) R_{\e}(i)
 \scal{r}^M \right) . $$
We use the the induction assumption for the first factor, which
belongs to $\Schn{q/(j-1)}$. We estimate the second factor using
lemma $ \refe{aprioriSchatten0} $ and $ \reff{raccourci} $ and
conclude using $ \reff{HolderSchatten} $. \finpreuve

\medskip

The main consequence of this lemma is the following.
\begin{prop} \label{reste4} For each $N \geq N_0$ large enough and $ M \geq 0 $,
 $ q > n-1 $ as above and $ j\geq  1 $ there exists $ C $ and $ M^{\prime} $ such that
\begin{eqnarray}
\left| \left| \scal{r}^M  \partial_{\e}^j  \left(  R_{\e}(z) -
Q_{N,\e}(z) \right) \scal{r}^M \right| \right|_{\bf q / j} \leq C
h^{N - n j / q}  \frac{ \scal{ \emph{Re} \ z}^{M^{\prime}}}{| \emph{Im}
\ z |^{M^{\prime}}}
\end{eqnarray}
for all $ \e \in [0,1] $ and $ h \in (0,1] $.
\end{prop}

\medskip

\noindent {\it Proof.} Since $h^N \partial_{\e}^j \left( R_{\e}(z)
{\mathcal R}_{N,\e}(z) \right) $ which is a linear combination of
$ h^N \partial_{\e}^{j_1}  R_{\e}(z)
\partial_{\e}^{j_2} {\mathcal R}_{N,\e}(z) $ with $j_1 + j_2 = j$,
we have to estimates operators of the form
$$ h^N \left( \scal{r}^M  (\partial_{\e}^{j_1} R_{\e}(z)) R_{\e}(i)^{j_1 -1}
\scal{r}^{M} \right) \left( \scal{r}^{-M}
(H_{\e} + i)^{j_1-1} \partial_{\e}^{j_2} {\mathcal R}_{N,\e}(z)
 \scal{r}^{M} \right) . $$
The first factor can be estimated by proposition $
\refe{restemoitie} $ and the second one by $ \reff{Schatten1} $
 and lemma $ \refe{resteSchwartz} $.
\finpreuve

\medskip

Using formula $ \reff{HelfferSjostrand} $, we get directly the following result.

\begin{theo} \label{fonctionne}     
  Let $ q > n - 1 $, $f \in {\mathcal S}(\Ra)$ and $N \geq 1$. There
 exists a pseudo-differential operator $Q_{N,\e}^f = \sum_{ l < N
} h^l A_l$ with $ A_l $ of the form $\reff{formepseudo}$, with
symbols in $S_{\e}^{0, - \infty}$ such that for all $ M
\geq 0$
$$ \left| \left| \partial_{\e}^j \scal{r}^M \left(
f(H_{\e}) - Q_{N,\e}^f \right) \scal{r}^M \right| \right|_{\bf q /j}
\leq C h^{N-nj/ q}, \qquad h \in (0,1], \ \e \in [0,1].
$$
Here $C$ depends on a finite number of semi-norms of $ f $. In
each chart, the symbols of $A_l$ are linear combinations of
$d_{l^{\prime} \! ,l,\e} f^{(l)}(p_{\e})$ and in particular the
principal symbol is $f(p_{\e})$.
\end{theo}

\subsection{Proofs of theorem $\bf \refe{theo1}$ and lemma 
$\bf \refe{lemmeformexi} $}  \label{formeoperateur}      
Theorem $ \refe{theo1} $ is a direct consequence of theorem $
\refe{fonctionne} $.  We obtain 
$ \reff{heatexpansion} $ by considering $ h = t^{1/2} $ and 
$ f \in {\mathcal S} (\Ra) $ such that $ f (\lambda) = e^{-\lambda} $ near the 
spectra of the operators. A priori, $ \reff{heatexpansion} $  involves
 all powers of the form $ t^{ (k-n)/2} $, but a standard argument
shows that the coefficients corresponding to odd $k$ vanish, since they
correspond to integrals of odd functions on the sphere. \finpreuve

\smallskip

For the proof of lemma $ \refe{lemmeformexi} $, we have to show that
$ \mbox{tr} \{f(H_{\e}^{(\kappa)}) V^{(\kappa)} \}_{\bf q} \rightarrow 
\mbox{tr} \{f(H_{\e}) V \}_{\bf q}  $ as $ \kappa \downarrow 0 $, 
for all $ f \in {\mathcal S}(\Ra)$. 
Using the explicit expressions of the symbols of $ Q_{N,\e}^{(\kappa),f} $
associated to $H_{\e}^{(\kappa)}$ (with the notation of theorem
$\refe{fonctionne}$), it is easy to check that 
$$  \left\{ V^{(\kappa)} Q_{N,\e}^{(\kappa),f}  \right\}_{\bf q} \rightarrow
\left\{ V Q_{N,\e}^f \right\}_{\bf q} \qquad
\mbox{in the trace class} $$
with a trace norm uniformly bounded by $ C \scal{\re z}^{M}/ |\im z|^{M} $ for 
some $M$ and $ C $.
 Thus we are left to the study of the remainder which, via Helffer-Sj\"ostrand
 formula, reduces to the study of the remainder given by $ \reff{basereso} $.
By the resolvent identity, we have
$$ R_{\e}^{(\kappa)} (z) - R_{\e} (z) = - R_{\e}^{(\kappa)}(z) 
\left( V^{(\kappa)} -V \right) R_{\e}(z)  , \qquad R_{\e}^{(\kappa)}(z) = 
(H_{\e}^{(\kappa)}-z)^{-1} $$
with $ ( V^{(\kappa)} -V ) R_{\e}(z) \rightarrow 0 $ in operator norm. This
shows that $ R_{\e}^{(\kappa)} (z) \rightarrow R_{\e}(z) $ in operator norm. 
Using $
\reff{Neumann} $ with $ k = 0 $ for $  H_{\e}^{(\kappa)} $ it is easy to check
that $ H_{\e} R_{\e}^{(\kappa)} (i) $ is bounded in operator norm uniformly
w.r.t. $ \kappa $. The resolvent idendity then shows that 
$ H_{\e} R_{\e}^{(\kappa)} (i) \rightarrow H_{\e} R_{\e}(i) $ in operator
norm and thus $ H_{\e} R_{\e}^{(\kappa)}(z) $ converges as well. By induction
$ H_{\e}^{k+1} R_{\e}^{(\kappa)}(z)^{k+1} \rightarrow H_{\e}^{k+1} R_{\e}(z)^{k+1} $ 
in operator norm with a uniform bound of the form $ C \scal{\re z}^M / 
|\im z|^M $ (apply $\partial_z^k$ to the resolvent identity for instance).
Then, we see that $ (\partial_{\e}^j R_{\e}^{(\kappa)}(z))
R_{\e}^{(\kappa)}(z)^{j-1} $ converges in $ \Schn{\bf q/j} $ with the same kind
of bound as in lemma $ \refe{restemoitie} $ and we can repeat the proof of proposition
$ \refe{reste4} $. The expected convergence follows from Helffer-Sj\"ostrand 
formula and dominated convergence. \finpreuve

\section{The proof of theorem $\bf \refe{theo2} $} \label{main}
\setcounter{equation}{0}
This section is entirely devoted to the proof of theorem $ \refe{theo2} $ and more particularly
to the proof of the asymptotic expansion. The continuity of $ \xi_{\bf q} $ on the absolutely
continuous spectrum is a consequence of the method. More precisely, we shall explicitely show
that $\xi_{\bf q}(\mu,h)$ is continuous in a neighborhood of $ \mu = 1 $. The proof of the
continuity of $ \xi_{\bf q}(\lambda) $, which we omit, would follow from the same method using 
$ \reff{propagationestimate} $ for fixed $ h $.
\subsection{Isozaki-Kitada's method in the asymptotically hyperbolic
case} \label{main1} Let us consider a covering of $Y$ by a finite
number of open sets $U_{n-1}$, each one of them being a relatively compact
susbet of a coordinate patch $ \widetilde{U}_{n-1} $ and consider their
respective images on $ \Ra^{n-1} $, i.e. $ \Omega  $ and 
$\widetilde{\Omega}$, under the coordinates maps. For
each such $\Omega$, we choose $\Omega_k$ for $ k = 1,\cdots,5$
such that
$$ \Omega = \Omega_5 \Subset \Omega_4 \Subset \Omega_3 \Subset
\Omega_2 \Subset \Omega_1 \Subset \widetilde{\Omega} . $$
Furthermore we can
assume that $ \Omega_k $ is open and convex (this will be usefull only for
proposition  $ \refe{convexfactor} $) for all
$k$. Then we consider, the outgoing areas $ \Gamma_1^+
\supset \Gamma_2^+ \supset \cdots \supset \Gamma_5^+ $ defined by
$$ \Gamma^+_{k} = \Gamma^+ (R^{k} , \varepsilon^{k} , w^{k} , \Omega_k) , \qquad k = 1 , \cdots , 5 . $$
We assume that $R > 1$ and $ \varepsilon, w \in (0,1) $ thus it is
clear that $ \overline{\Gamma^+_{k} } \subset \Gamma^+_{k-1} $ for
$ k = 2 , \cdots, 5 $. More precisely  one can always find a
cutoff function supported in $\Gamma^+_{k-1}$ which is $ \equiv 1
$ on $ \Gamma^+_k $. We can choose for instance
\begin{eqnarray}
\chi_{R^k} (r) \chi_{\Omega_k} (y) \chi_+ (\rho) \chi_{
\varepsilon^k} (e^{-2r}g(y,\eta)) \chi_{w^k} (\rho^2 +
e^{-2r}g(y,\eta) ) \label{cutoffsortant}
\end{eqnarray}
with smooth functions $ \chi_{R^k} $, $ \chi_{\Omega_k} $, $
\chi_{ \varepsilon^k} $ and $ \chi_{w^k} $ such that
\begin{eqnarray*}
\mbox{supp} \ \chi_{R^k} \subset (R^{k-1} , \infty) \  \ &
\mbox{and} & \ \  \chi_{R^k} \equiv 1 \ \ \mbox{on} \ \ (R^k ,
\infty), \\ \mbox{supp} \ \chi_{\Omega_k} \subset \Omega_{k-1}
 \ \ & \mbox{and} &  \ \
\chi_{\Omega_k} \equiv 1 \ \ \mbox{on} \ \ \Omega_k \\
\mbox{supp} \ \chi_{\varepsilon^k} \subset (- \infty ,
\varepsilon^{k-1} ) \ \ & \mbox{and} & \ \ \chi_{\varepsilon^k}
\equiv
1 \ \ \mbox{on} \ \ ( - \infty , \varepsilon^k), \qquad \\
\mbox{supp} \ \chi_{w^k} \subset (1 - w^{k-1} , 1 + w^{k-1}) \ \
 & \mbox{and} & \ \ \chi_{w^k} \equiv 1 \ \ \mbox{on} \ \ (1-w^k , 1 +
w^k) .
\end{eqnarray*}
 Since we will choose $w$ and $ \varepsilon $ small, it is enough
to consider $ \chi_+$ supported into $(0, \infty)$ and such that
$\chi_+ \equiv 1$ on $ (1/2 , \infty) $. Note that on  $
\Gamma^{+}_k $, we always have
\begin{eqnarray}
  \left( 1 - w^{k} - \varepsilon^{k} \right)^{1/2}   \leq
\rho  \leq \left( 1 + w^{k}  \right)^{1/2} . \label{rhononzero}
\end{eqnarray}
We are now ready to construct the functions needed for
Isozaki-Kitada's method in the hyperbolic case. We start with the
following proposition.
\begin{prop} \label{propquatre1} For $R$ large enough, and $ \varepsilon , w $ small
enough, there exists a smooth function $ \varphi_{\e}^+
(r,y,\rho,\eta) $ defined, for any $\e \in [0,1]$, on $ \Gamma_1^+
$ such that
$$  (\partial_{r} \varphi_{\e}^+ )^2 + e^{-2r} g_{\e}(r,y,\partial_y \varphi_{\e}^+) = \rho^2 . $$
This function satisfies the following estimates for $ \e \in [0,1]
$ and  $ (r,y,\rho,\eta) \in \Gamma_1^+ $
\begin{eqnarray}
  \left| \partial_{\e}^j \partial_r^k \partial_{\rho}^l
\partial_y^{\alpha} \partial_{\eta}^{\beta}
 \left(  \varphi_{\e}^+ - r \rho - y.\eta \right) \right| \leq C_{j,k,l,\alpha,\beta} e^{-(j+1)r}|\eta| .
 \label{phasesortante}
\end{eqnarray}
\end{prop}
This proposition solves the equation $ \reff{HamiltonJacobi} $ in
the case $p_{\e}(x,\xi) = \rho^2 + e^{-2r} g_{\e}(r,y,\eta)$ and
explains how differentiation w.r.t. $\e$ provides exponential
decay.

\medskip

\noindent {\it Proof.} We follow the principle explained in subsection
$ \refe{IsozakiKitadaformel} $, using the function
$ S_{\e}^{+} (t,r,y,\rho,\eta) $ given by proposition $ \refe{eiko2} $.
Since $S_{\e}^+$ solves $ \reff{eikonal} $ and is a generating function of the
flow (see $\reff{generatingflow}$) we have
\begin{eqnarray*}
 \partial_t S_{\e}^+ =
p_{\e} (r,y,\partial_{r} S_{\e}^+ , \partial_y S_{\e}^+) & =
 & p_{\e} (\partial_{\rho} S_{\e}^+ , \partial_{\eta} S_{\e}^+,\rho,\eta) \\
 & = & \rho^2 + e^{- 2 \partial_{\rho} S_{\e}^+}
 g_{\e} \left( \partial_{\rho} S_{\e}^+ , \partial_{\eta} S_{\e}^+  ,\eta \right)
\end{eqnarray*}
where one remark that the last term of the second line
 is ${\mathcal O}(e^{-r - 2t\rho})$
by $\reff{approxeiko}$ hence is integrable.
Thus we can use formula $ \refe{definitionphase} $ with
 $\widetilde{S}_{\e}^{+}(t,\rho,\eta) = t \rho^2 $ to define $\varphi_{\e}^+$,
and then $ \reff{phasesortante} $ is a direct consequence of
 proposition $ \refe{eiko2} $. \finpreuve

\medskip

The next lemma is a preparation lemma for the resolution of the
transport equations. We do not prove it since it can be obtained very similarly
to the estimates on the geodesics of section $ \refe{classical} $.
\begin{lemm} \label{lemmequatre1} There exist $R$ large enough and $w,\varepsilon$ small
enough, the following result holds: for all $(r,y,\rho,\eta) \in
\Gamma_2^+$ and $ \e \in [0,1] $, the solution $ (
\check{r}_{\e}^t , \check{y}_{\e}^t )(r,y,\rho,\eta)$ of
$$ \frac{d}{dt} \left( \begin{matrix}  \check{r}_{\e}^t \\ \check{y}_{\e}^t \end{matrix} \right) =
\left( \begin{matrix} 2 \partial_r \varphi_{\e}^+ (
\check{r}_{\e}^t , \check{y}_{\e}^t , \rho , \eta ) \\ e^{- 2
\check{r}_{\e}^t} (\partial_{\eta} g_{\e} ) ( \check{r}_{\e}^t ,
\check{y}_{\e}^t , \partial_y \varphi_{\e}^+ (\check{r}_{\e}^t ,
\check{y}_{\e}^t , \rho, \eta ) )
\end{matrix} \right) , \qquad
 \left( \begin{matrix}  \check{r}_{\e}^0 \\
\check{y}_{\e}^0 \end{matrix} \right) = \left(
\begin{matrix}  r \\ y  \end{matrix}
\right)
$$
is defined for all $t \geq 0$ and satisfies $ (  \check{r}_{\e}^t
, \check{y}_{\e}^t , \rho , \eta ) \in \Gamma_1^+ $. Furthermore,
we have the estimates
\begin{eqnarray}
\left| \partial_{\e}^j \partial_r^k \partial_{\rho}^l
\partial_y^{\alpha} \partial_{\eta}^{\beta}
 \left(  \check{r}_{\e}^t - r  - 2 \rho t \right) \right| + \left| \partial_{\e}^j \partial_r^k \partial_{\rho}^l
\partial_y^{\alpha} \partial_{\eta}^{\beta}
 \left(  \check{y}_{\e}^t  - y \right) \right| \leq C_{j,k,l,\alpha,\beta}
 e^{-(j+1)r}|\eta| \label{flowtransport}
\end{eqnarray}
for all $ \e \in [0,1] $ and $ (r,y,\rho,\eta) \in \Gamma^+_2 $.
\end{lemm}
This lemma, which gives in particular a precise sense to $
\reff{spacetime} $ in the current context, allows to define
functions $ a_{\e}^{(0)} , \cdots , a_{\e}^{(N)} $ of
$r,y,\rho,\eta \in \Gamma_2^+$ according to the formulas $
\reff{abrzerob} $, $ \reff{abrjb} $ where one has of course to
replace $ \check{x}_{\e}^t  $ by $ \check{r}_{\e}^t,
\check{y}_{\e}^t $ and use the following explicit expression
$$ c_{\e} = \partial_r^2 \varphi_{\e}^+ + e^{-2r} g_{\e}(r,y,\partial_y)
\varphi_{\e}^+ + e^{-r} 
\sum_{l+|\alpha|=1} \tilde{v}_{\e}^{l,\alpha}(e^{-r},y) (e^{-r} \partial_y)^{\alpha}
\partial_r^l
 \varphi_{\e}^+ ,
$$
where $ \tilde{v}_{\e}^{l,\alpha}  = v_0^{l,\alpha} + \e (v_1^{l,\alpha}-
v_0^{l,\alpha}) $, with the notations of $ \reff{ordreundiff} $. 
By $ \reff{phasesortante} $ we see easily  that  $
c_{\e}(r,y,\rho,\eta) = {\mathcal O}(e^{-r}\scal{\eta}) $ (note the factor 
 $ \scal{\eta} $ instead of $|\eta|$ which is caused by the term $ \partial^l_r 
\varphi^+_{\e}$ with $l=1$ in the expression of $c_{\e}$).
 This implies
that the integral in $ \reff{abrzerob} $ is convergent since $
\reff{rhononzero} $ and $ \reff{flowtransport} $ show that, for
some $c > 0 $,
$$ c_{\e}(\check{r}_{\e}^t , \check{y}_{\e}^t , \rho , \eta) = 
{\mathcal O}(e^{-ct-r} \scal{\eta}) , \qquad t \geq 0 . $$
Thus $ a_{\e}^{(0)} $ is well defined on $ \Gamma^+_2 $. By
induction, one checks that, for $ m \geq 1 $, we have
$$ P_{\e}(r,y,D_r, D_y)
a_{\e}^{(m-1)} (\check{r}_{\e}^t, \check{y}_{\e}^t , \rho,\eta) =
{\mathcal O} (e^{-ct-r}\scal{\eta}) , \qquad t \geq 0,
$$
on $\Gamma_2^+$, for all $ \e \in [0,1] $ and thus $a_{\e}^{(m)}$
is well defined on $ \Gamma^+_2 $ for all $ m $. More generally,
by mean of proposition $\refe{propquatre1}$, lemma 
$\refe{lemmequatre1}$ and $\reff{rhononzero}$, we obtain easily
the following result.
\begin{prop} \label{propquatre2} The functions $a_{\e}^{(0)}, \cdots , a_{\e}^{(N)}$
defined in $\Gamma_2^+$ by $ \reff{abrzerob} $ and $ \reff{abrjb}
$ satisfy
\begin{eqnarray}
 \left| \partial_{\e}^j \partial_r^k \partial_{\rho}^l
\partial_y^{\alpha} \partial_{\eta}^{\beta}
 \left(  a_{\e}^{(0)}(r,y,\rho,\eta)  - 1  \right) \right| & \leq & C_{j,k,l,\alpha,\beta}
 e^{-(j+1)r}\scal{\eta} , \nonumber \\
 \left| \partial_{\e}^j \partial_r^k \partial_{\rho}^l
\partial_y^{\alpha} \partial_{\eta}^{\beta}
  a_{\e}^{(m)}(r,y,\rho,\eta)   \right| & \leq & C_{j,k,l,\alpha,\beta}
 e^{-(j+1)r}\scal{\eta} , \qquad 1 \leq m \leq N \nonumber
 \end{eqnarray}
 for all $ \e \in [0,1] $ and $ (r,y,\rho,\eta) \in \Gamma_2^+ $.
\end{prop}
The details of the proof are left to the reader. They follow from
the explicit expressions of the functions $ a_{\e}^{(m)} $.

In order to consider functions defined on $ \Ra^{2n} $, we choose
a cutoff $\chi_{2,3}(r,y,\rho,\eta)$ of the form $
\reff{cutoffsortant} $ which is supported in $ \Gamma_2^+ $ and $
\equiv 1 $ on $ \Gamma_3^+ $. Then we multiply the functions
$a_{\e}^{(0)}, \cdots , a_{\e}^{(N)}$ by $\chi_{2,3}$ and define
$$ a_{\e} = \chi_{2,3} \left( a_{\e}^{(0)} + h a_{\e}^{(1)} + \cdots + h^N a_{\e}^{(N)} \right) . $$
Then, following Isozaki-Kitada's method as explained in subsection
$ \refe{IsozakiKitadaformel} $, we consider
$$ H_{\e} J (\varphi_{\e}^+ , a_{\e}) - J (\varphi_{\e}^+ , a_{\e}) P , \qquad \mbox{with} \qquad P = h^2 D_r^2 . $$
Of course, the choice of $ P $ follows from $\reff{HJapprox}$
since $ \lim_{r \rightarrow \infty} p_{\e}(r,y,\rho,\eta) = \rho^2
= : p ( \rho) $. The operator $ H_{\e} J (\varphi_{\e}^+ , a_{\e})
- J (\varphi_{\e}^+ , a_{\e}) P $ has the following expression
\begin{eqnarray}
  J  \left( \varphi_{\e}^+  , a_{\e}^{\prime}
\right) +
 h^{N+2} J \left( \varphi_{\e}^+ , \chi_{2,3} P_{\e}(r,y,D_r,D_y)
a_{\e}^{(N)} \right) \label{negligeable0}
\end{eqnarray}
where $ a_{\e}^{\prime} $ is a linear combination of products of {
derivatives (of order $\geq 1$) of $\chi_{2,3}$ and of derivatives
of $a_{\e}^{(0)}, \cdots , a_{\e}^{(N)}$. The first term of $
\reff{negligeable0} $ is produced by all the derivatives due to
$H_{\e}$ which may fall on $\chi_{2,3}$. The amplitude of the
second term is nothing but $ \chi_{2,3} (\widetilde{a}_{\e}^{(0)}
+ \cdots + h^{N+2} \widetilde{a}_{\e}^{(N+2)} ) $ with the
notations of subsection $ \refe{IsozakiKitadaformel} $. This
follows from the construction of $\varphi_{\e}^+$ and
$a_{\e}^{(0)}, \cdots , a_{\e}^{(N)}$ which satisfy respectively
Hamilton-Jacobi's equation and the transport equations on the
support of $ \chi_{2,3} $. The negligibility of  $
\reff{negligeable0} $, or more precisely of $ \reff{Isozakitrick}
$, will be a consequence of the next lemma.
\begin{lemm} \label{clefnegligeable} For $R$ large enough, and $ \varepsilon , w $ small
enough we have the following property: for any symbol $b
(r,y,\rho,\eta)$ supported in $ \Gamma_4^+ $ and such that, for
all $M$,
$$ \left|  \partial_r^k \partial_{\rho}^l
\partial_y^{\alpha} \partial_{\eta}^{\beta} b(r,y,\rho,\eta)  \right|  \leq
C_{k,l,\alpha,\beta,M} \scal{r}^{-M} $$ we have the following
estimates for $ s \geq 0 $
\begin{eqnarray}
\left| \left| \scal{r}^M  J  \left( \varphi_{\e}^+  ,
a_{\e}^{\prime} \right) U (s) J (\varphi_{\e}^+,b)^* \scal{r}^{M}
\right| \right|_{\infty}  \leq  C_{M} h^{M} \scal{s}^{-M} ,
\label{supernegligeable}
  \\
\left| \left| \scal{r}^M  J  \left( \varphi_{\e}^+  ,  \chi_{2,3}
P_{\e}(r,y,D_r,D_y) a_{\e}^{(N)} \right) U (s) J
(\varphi_{\e}^+,b)^* \scal{r}^{M} \right| \right|_{\infty}  \leq
C_{M} h^{-n_0} \scal{s}^{-M} \label{justenegligeable}
\end{eqnarray}
for all $M$ and some universal constant $n_0$. Here $U(s) = e^{-i
\frac{s}{h} P}$ is the propagator of $P$.
\end{lemm}
Note that the power $h^{-n_0}$ on the right hand side of $
\reff{justenegligeable} $ is harmless since we have a power
$h^{N+2}$ in $ \reff{negligeable0} $, with $N$ arbitrarily large.
The proof of this lemma is not very hard and follows from suitable
integrations by parts on the Schwartz kernels of the operators which are
explicitly given by oscillatory integrals. However its proof is a
bit long and we have postponed it to appendix $ \refe{rappelmicrolocaux} $.

The last step of the construction is the factorization of
pseudo-differential operators. This is the purpose of the
following proposition.
\begin{prop} \label{convexfactor}  There exists $R$ large enough and $ \varepsilon , w $
small enough such that the following property holds: for any
symbol $ c_{\e} \in S$, supported in $\Gamma^+_5$, one can find
symbols $b_{\e}^{(0)} , \cdots , b_{\e}^{(N)} \in S $ supported in
$ \Gamma^+_4 $ such that
$$  J ( \varphi_{\e}^+ , a_{\e} ) J ( \varphi_{\e}^+, b_{\e} )^* - O \! p_h (c_{\e}) \ \ \mbox{is} \ h^N \
\mbox{negligible} $$ with $b_{\e} = b_{\e}^{(0)} + h b_{\e}^{(1)}
+ \cdots + h^N b_{\e}^{(N)} $.
\end{prop}

\noindent {\it Proof.} As explained in subsection
$\refe{IsozakiKitadaformel}$, we need to study the map
\begin{eqnarray}
(\rho,\eta) \mapsto \Phi_{\e}^{+}(r,y,r^{\prime},y^{\prime},\rho,\eta) = 
\int_0^1
\partial_{r,y} \varphi_{\e}^+
(r^{\prime}+t(r-r^{\prime}),y^{\prime}+t(y-y^{\prime}),\rho,\eta)
\  \mbox{d} t  \label{mapdiffe}
\end{eqnarray}
It is then clear that if $(r,y,\rho,\eta) $ and 
$ (r^{\prime},y^{\prime},\rho,\eta)  $  belong to $  \Gamma^+
(R,\varepsilon,w,\Omega)  $ with $\Omega$ convex, then 
$ (r^{\prime}+t(r-r^{\prime}),y^{\prime}+t(y-y^{\prime}),\rho,\eta) \in 
\Gamma^+ (R, C_0 \varepsilon , w^{\prime} ,\Omega) $ for all $ t \in
[0,1]$, and any $ w^{\prime} > w $ which can be chosen as close to $ w $ as we
want by decreasing $ \varepsilon $. Here $ C_0 $ is given by 
$ \reff{elliptic0} $. In particular,
$$  \left| \partial^{\gamma} 
 \left( \Phi^+_{\e} - (\rho,\eta) \right) \right| \leq  C_{\gamma} 
 e^{-(j+1)  \min (r,r^{\prime}) } |\eta| $$
with $ \partial^{\gamma} = \partial_{\e}^j \partial_r^k 
\partial_{r^{\prime}}^{k^{\prime}} \partial_{\rho}^l
\partial_y^{\alpha} \partial_{y^{\prime}}^{\alpha^{\prime}}
\partial_{\eta}^{\beta} $,
for all $(r,y,\rho,\eta) $ and 
$ (r^{\prime},y^{\prime},\rho,\eta)  $  in $  \Gamma^+
(R,\varepsilon,w,\Omega)  $. Using the formula $ \reff{formebzero} $
(with $\theta_{\e} = \Phi_{\e}^+$), we can define $ b^{(0)}_{\e}$ and by
iterations $ b_{\e}^{(1)}, \cdots $, from $ c_{\e} $. In particular, if
$ c_{\e} $ is supported in $ \Gamma^+_5 $, 
$ b^{(0)}_{\e}, b_{\e}^{(1)} , \cdots $ will be clearly supported in $
\Gamma^+_4 $ if $R$ is large enough and $ \varepsilon , w $ small enough.
In order to check that the difference
$J ( \varphi_{\e}^+ , a_{\e} ) J ( \varphi_{\e}^+, b_{\e} )^* - O \! p_h
(c_{\e})$ is $ h^N $ negligible, we look at its Schwartz kernel
which we split into three terms using the partition of unit
$$ 1 = \theta_{0}(r-r^{\prime}) + \theta_+ (r-r^{\prime}) + \theta_+ 
(r^{\prime} -r ) $$
with $ \theta_0 = 1 $ close to $ 0 $  and $\theta_+$ supported in
$[1,\infty)$. The term corresponding to $ \theta_0 $
 behaves nicely by construction. 
The off diagonal part is ${\mathcal O}(h^{\infty})$ and is the sum
of two oscillatory integrals which are either ${\mathcal O}(e^{-(j+\nu)r})$
or $ {\mathcal O}(e^{- (j + \nu ) r^{\prime}}) $ after application 
$\partial_{\e}^j$, according to the sign of $ r - r^{\prime} $.
 We don't go any further into details. \finpreuve

\subsection{The dependence on ${\bf \kappa}$ and the remainders}

Our next task is explain how to deal with `remainders'. In subsection
$ \refe{representation} $ or in section $ \refe{pseudo}  $, we have seen why and how
a lot of expansions of operators, in powers of $h$, could be obtained. It
 is now necessary to show why the contributions of the
remainders of such expansions can indeed be neglected in the
expected expansion. We introduce the notation $ \equiv_N^{n_1} $,
which we will use extensively in the next subsection. Its meaning
is the following.
\begin{defi} \label{definitionmodulo} We write
$$ \emph{tr} \left(  T_{\e}^{(\kappa)}(t,h)  \right) \  \equiv_N^{n_1} \  0  $$
if $T_{\e}^{(\kappa)}(t,h)$ is a family of
 trace class operators
for $\kappa \in (0,1]$, whose trace is $C^{q-1}$ with respect to
$\e \in [0,1] $, measurable with respect to $t\in \Ra$, and if
there exists $C_{N,n_1}$, independent of $h$, such that
\begin{eqnarray}
\lim_{\kappa \downarrow 0} \left\{ \emph{tr} \left(
T_{\e}^{(\kappa)}(t,h) \right)
 \right\}_{\bf q} \ \mbox{exists in } \ {\mathcal
 S}^{\prime}(\Ra_t) \ \mbox{and belongs to} \ L^1 (\Ra,\mbox{d}t),  \label{limiteSchwartz1}  \\
 \int_{-\infty}^{+\infty} \left|  \lim_{\kappa \downarrow 0}
\left\{ \emph{tr} \left( T_{\e}^{(\kappa)}(t,h) \right)
 \right\}_{\bf q}  \right| \ \emph{d}t \  \leq \  C_{N,n_1} h^{N-n_1}
 . \label{limiteSchwartz2}
\end{eqnarray}
If the same result holds when integrating on $[0,+\infty)$, then
we will write $ \equiv_{N,+}^{n_1} $ instead of $\equiv_N^{n_1}$.

If $  S_{\e}^{(\kappa)}(t,h)     $ is another family of trace
class operators for $\kappa \in (0,1] $, then
$$  \emph{tr}  \left( T_{\e}^{(\kappa)}(t,h) \right) \equiv_N^{n_1}
 \emph{tr}  \left( S_{\e}^{(\kappa)}(t,h) \right) \qquad
\Longleftrightarrow \qquad \emph{tr}
 \left( T_{\e}^{(\kappa)}(t,h)  -
 S_{\e}^{(\kappa)}(t,h) \right) \ \equiv_N^{n_1} \  0 ,   $$
and similarly for $\equiv_{N,+}^{n_1}$.
\end{defi}
Operators (or rather traces) satisfying $ \reff{limiteSchwartz1} $
and $ \reff{limiteSchwartz2} $ are of interest since
\begin{eqnarray}
 \lim_{\kappa \downarrow 0} {\mathcal F}_h^{-1} \left\{ \mbox{tr}
\left( T_{\e}^{(\kappa)}(t,h) \right) \right\}_{\bf q} = {\mathcal
F}_h^{-1} \lim_{\kappa \downarrow 0} \left\{ \mbox{tr} \left(
T_{\e}^{(\kappa)}(t,h) \right) \right\}_{\bf q} = {\mathcal
O}(h^{N-n_1-1}) \ \ \mbox{in} \ C^0(\Ra).
\label{passagelimitefaible}
\end{eqnarray}

 In practice, 
  it is enough to show, for instance, that $ \lim_{\kappa \rightarrow 0} \{ T_{\e}^{(\kappa)}(t,h)
  \}_{\bf q}  =: \{ T_{\e}^{(0)}(t,h)
  \}_{\bf q}   $ exists, in the weak topology of bounded operator
and that
\begin{eqnarray}
 \left| \left| \{ T_{\e}^{(0)}(t,h)
  \}_{\bf q} \right| \right|_{\bf 1} & \leq  & C_N h^{N-n_1} \psi (t), \qquad  t \in \Ra, \ h \in (0,1],
  \label{suffisant0}
  \\
\lim_{\kappa \rightarrow 0} \  \mbox{tr} \{
T_{\e}^{(\kappa)}(t,h)
  \}_{\bf q} & = & \mbox{tr} \{ T_{\e}^{(0)}(t,h)
  \}_{\bf q} \qquad t \in \Ra, \ h \in (0,1], \label{suffisant1}
  \\
 \left| \left| \{ T_{\e}^{(\kappa)}(t,h)
  \}_{\bf q} \right| \right|_{\bf 1} & \leq & C h^{-n_2} \scal{t}^{n_2},
  \qquad t \in \Ra, \ h \in (0,1] , \ \kappa \in [0,1] \label{suffisant2}
\end{eqnarray}
for some $L^1$ function $\psi$ in $\reff{suffisant0}$, and some
$n_2 \geq 0$ in $ \reff{suffisant2} $. This last condition could
be weakened since we could clearly allow $C$  to depend
arbitrarily on $h$. Usually, $ \reff{suffisant1} $ and $
\reff{suffisant2} $ are easy to check and the non trivial part of
the job is to show $ \reff{suffisant0} $.

Before giving explicit examples, let us explain how we will use
definition $ \refe{definitionmodulo} $. Recall that we are
studying $ \reff{formezero} $ which is nothing but
\begin{eqnarray}
\frac{1}{2 \pi h} \int_{- \infty}^{+ \infty} e^{\frac{i}{h} \mu t
} \ \mbox{tr} \left( f (H_{\e}^{(\kappa)}) U_{\e}^{(\kappa)}(t)
V^{(\kappa)} \tilde{f}(H_{\e}^{(\kappa)}) \right) \ \mbox{d} t
\label{coreFourierinverse}
\end{eqnarray}
 with $ \tilde{f}f = f $ which we choose
to be supported close to $1$.  The method is the following: if for
any $N$ large enough, we can find $T_{\e,N}^{(\kappa)}(t,h)$ such
that
$$  \mbox{tr} \left( f (H_{\e}^{(\kappa)}) U_{\e}^{(\kappa)}(t)
V^{(\kappa)} \tilde{f}(H_{\e}^{(\kappa)}) \right) \ \equiv_N^{n_1}
\ \mbox{tr} \left(   T_{\e,N}^{(\kappa)}(t,h) \right)  $$ for some
$n_1$ independent of $N$,  and moreover such that $\{ \mbox{tr} (
T_{\e,N}^{(\kappa)}(t,h) ) \}_{\bf q}$ converges in ${\mathcal
S}^{\prime}(\Ra_t)$ as $\kappa \downarrow 0$ with a limit in
$L^1(\Ra,\mbox{d}t)$ satisfying
$$  \frac{1}{2 \pi h} \int_{- \infty}^{+ \infty} e^{\frac{i}{h} \mu t
} \  \lim_{\kappa \downarrow 0} \left\{ \mbox{tr} \left(
T_{\e,N}^{(\kappa)}(t,h) \right) \right\}_{\bf q} \ \mbox{dt} =
h^{-n} \sum_{k=0}^N h^k \alpha_k (\mu) + {\mathcal O} (h^{N-n})
$$ with $\alpha_k \in C^0(I)$ and  ${\mathcal O}(h^{N-n})$
understood in the topology of $C^0(I)$, then we get the existence
of the expansion $ \reff{expansionphase}$. This is due to $
\reff{passagelimitefaible} $ since it implies  that
$$ \xi_{\bf q}(\mu,h) -  \frac{1}{2 \pi h} \int_{- \infty}^{+ \infty} e^{\frac{i}{h} \mu t
} \  \lim_{\kappa \downarrow 0} \left\{ \mbox{tr} \left(
T_{\e,N}^{(\kappa)}(t,h) \right) \right\}_{\bf q} \ \mbox{dt} =
{\mathcal O}(h^{N-n_1-1}) \ \ \mbox{in} C^0 $$
 with $N-n_1 -1$ arbitrarily large if $N$ is large.

\medskip

We shall now give explicit and useful examples of operators
satisfying $ \reff{limiteSchwartz1} $ and $ \reff{limiteSchwartz2}
$. To that end, we will use extensively the fact that for all $f
\in C_0^{\infty}(\Ra)$, in the strong sense,
$$ U_{\e}^{(\kappa)}(t) f (H_{\e}^{(\kappa)}) \rightarrow U_{\e}(t) f(H_{\e}), 
\qquad \kappa \downarrow 0 $$
 with a norm bounded by $ \sup |f| $. We omit the proof of this
 easy
 fact since it follows by  section $ \refe{pseudo}
 $ or more directly, by Helffer-Sj\"ostrand formula for instance.
 We will also use the following well known lemma, which is valid
 if $B^{(\kappa)}$ is a family of bounded operators and
 $T^{(\kappa)}$ a family of $ \Schn{q}$.
\begin{lemm} \label{superflu} If $B^{(\kappa)} \rightarrow B^{(0)}$ strongly
 and $T^{(\kappa)} \rightarrow T^{(0)}$ in $ \Schn{q} $
 then $  B^{(\kappa)} T^{(\kappa)}  \rightarrow  B^{(0)} T^{(0) } $ in 
 $ \Schn{q} $.
\end{lemm}


Let ${\mathcal R}_{\e,N}^{(\kappa)} (h) = V^{(\kappa)} \left( \tilde{f} (H_{\e}^{\kappa}) - 
Q^{\tilde{f}, (\kappa)}_{N,\e} \right)   $, with the notations of theorem
$ \refe{fonctionne} $ where we include the depence on $ \kappa $.

\begin{prop} Let $ I $ be an intervalle such that $ \reff{propagationestimate} $
 holds. Then  all $f,\tilde{f} \in C_{0}^{\infty}(I)$, with $ f \tilde{f} = f $,
 we have
$$ \emph{tr} \left( f(H_{\e}^{(\kappa)}) U_{\e}^{(\kappa)}(t)
{\mathcal R}_{N,\e}^{(\kappa)} (h) \right) \equiv_N^{n_1} 0  , 
\qquad n_1 = n + q    . $$
\end{prop}

\noindent {\it Proof.} We only show $ \reff{suffisant0} $. We first rewrite the trace as
$$  \mbox{tr}
\left( \scal{r}^{-M} f(H_{\e}^{(\kappa)}) U_{\e}^{(\kappa)}(t) \scal{r}^{-M} \
\scal{r}^M {\mathcal R}_{\e,N}^{(\kappa)} (h) \scal{r}^M \right) $$
and then apply $\partial_{\e}^{q-1}$. By Leibnitz rule and lemma
$ \refe{Leibnitz0} $ we obtain an explicit expression of this derivative.
We note that, for any $ 0 \leq k \leq q - 1  $,
\begin{eqnarray}
\scal{r}^M \partial_{\e}^{q-1-k}  {\mathcal R}_{\e,N}^{(\kappa)} (h)
\scal{r}^{M} = {\mathcal O}(h^{N - (q-k) n / q}) \qquad \mbox{in} \
\Schn{q/(q-k)} \label{casparticulierpseudo}
\end{eqnarray}
and that this estimate holds for $ \kappa \geq 0$, continuously w.r.t $\kappa$
and $\e$, which will allow
to let $\kappa \downarrow 0$.  This can been seen with the same argument as the
one used for the proof of lemma $ \refe{lemmeformexi} $ in subsection $
\refe{formeoperateur} $.
On the other hand, using the notations of lemma
$ \refe{Leibnitz0} $ in which we now take $\kappa$ into account,
we can write $ \scal{r}^{-M} \reff{integrale}
\scal{r}^{-M} $ as the integral over $F_t^j$ of
\begin{eqnarray}
h^{-j} \left( V^{(\kappa)} \scal{r}^{-M}\partial_{\e}^{l_0} S_{\e}^{(\kappa),\tilde{f}_0 }
\scal{r}^M \right)
\left( \scal{r}^{-M} U_{\e}^{(\kappa),f_0}(t_0) \scal{r}^{-M} \right) \left(
\scal{r}^M V^{(\kappa)}
  \partial_{\e}^{l_1} S_{\e}^{(\kappa),\tilde{f}_1 }  \scal{r}^M  \right)
  \cdots \nonumber \\ \cdots  \left( \scal{r}^{-M} U_{\e}^{(\kappa),f_j} (t)
  \scal{r}^{-M} \right)
 \left(  \scal{r}^{M}\partial_{\e}^{l_{j+1}} S_{\e}^{(\kappa),\tilde{f}_{ j + 1 },
  (\kappa)} \scal{r}^{-M} \right) . \qquad \label{astuceecriture}
\end{eqnarray}
Note that we used the notation $ \refe{notationpropagation} $. Then, as before,
it is not hard to check that
$$ \scal{r}^M (V^{(\kappa)})^j \partial_{\e}^l S^{(\kappa),\tilde{f}}_{\e}
\scal{r}^M
\rightarrow
\scal{r}^M V^{\tau} \partial_{\e}^l S^{(0),\tilde{f}}_{\e} 
\scal{r}^M \qquad \in \Schn{q/(q-l- {\bf \tau})}
$$
for $ \tau = 0 $ or $ 1 $. By lemma $ \refe{superflu} $ we can let $ \kappa
\downarrow 0 $ in $ \reff{astuceecriture} $, and using 
$ \reff{HolderSchatten} $
we get for $ \kappa = 0 $,
$$ \left| \left| \reff{astuceecriture} 
\scal{r}^M \partial_{\e}^{q-1-k}  {\mathcal R}_{\e,N}^{(0)} (h)
\scal{r}^{M} \right| \right|_{\bf 1} \leq C h^{N-n-j} 
\prod_{l=0}^{j} \left| \left| \scal{r}^{-M}
U_{\e}^{f_l}(t_l) \scal{r}^{-M}  \right| \right|_{\infty}  $$
Then the result follows from $ \reff{Convolution} $. \finpreuve

\smallskip

The same method can be applied to other kind of operators. We only explain
which modifications to do. 
For instance, if $ \partial_{\e}^j {\mathcal R}_{\e,N}^{(\kappa)} (h) $ 
can be written in any chart as $ \reff{astuceformelle} \times e^{-r} $
then everything works almost as well. This typically waht happens
when one considers remainders of the pseudo-differential expansion of the
adjoint of
 $ v^{(\kappa)} \tilde{f}(H_{\e}^{(\kappa)})  $ or with negligible operator
 involved in proposition $ \refe{convexfactor} $. 
Actually, in this case we don't have
$ \reff{casparticulierpseudo} $ anymore, but the exponential weights
in $ \reff{astuceformelle} $ can be put on both sides of $
\reff{astuceecriture}$ which then becomes trace class and we can conclude
similarly. This situation is also of interest when one has to consider 
the remainder involved in $ \reff{evolutiontrick} $. For the
latter, the exponential weights are provided by the differentiation w.r.t. 
$ \e $ in view of the explicit form of the symbols in term of the flow and
 theorem $ \refe{Precised} $.

We can also study the contribution of the right hand side of
$ \reff{Isozakitrick} $ using lemma $ \refe{clefnegligeable} $. Here there is 
one more
integral over $ s $ but we clearly have $ L^1 $ estimates by lemma $
\refe{clefnegligeable} $ and the convolution trick works again. 
Note finally that in this case we will choose $ n_1 = n+ q + n_0 $.

\subsection{The core of the proof}
Below, we will use extensively the notation $ \equiv_N^{n_1} $ of the previous
subsection with $ n_1 = n + q + n_0 $, which is independent of $ N $.

We start from the semi-classical Fourier transform of $ \reff{formezero} $ 
with $ f $ supported in $ [1-w^5/4,1+w^5/4] $ and $ f =
1 $ close to $ 1 $. We furthermore assume that the coefficients of 
$  V $ (hence of $V^{(\kappa)} $) are all supported in $ \{ r > R \} $.
We shall indicate how to modify the proof in the general case, however recall
that it is the case for the principal symbol. 

 Using theorem $ \refe{fonctionne} $, we get a pseudo-differential expansion
of $ V^{(\kappa)} f (H_{\e}^{(\kappa)}) $  whose symbols can be splited in two
parts  using a partition of unit
associated to $ \reff{split} $ so that get, for any
$N$,
\begin{eqnarray}
 \mbox{tr} \left( f (H_{\e}^{(\kappa)}) U_{\e}^{(\kappa)}(t)
V^{(\kappa)} \tilde{f}(H_{\e}^{(\kappa)}) \right) \equiv_N^{n_1} \sum_{k
\leq N , \ \Omega } h^k \mbox{tr} \left( f (H_{\e}^{(\kappa)})
U_{\e}^{(\kappa)}(t) O \! p_h ( \gamma_{\e,k}^{(\kappa)}) \right)
\label{dev0}
\end{eqnarray}
  where $
\gamma_{\e,k}^{(\kappa)} \in S_{\e}^{1,-\infty} $ and  either $
\mbox{supp} \gamma_{\e,k}^{(\kappa)} \subset \Upsilon^+
(R, w^5/2,\sigma_+,\Omega)$  or $ \mbox{supp}
\gamma_{\e,k}^{(\kappa)} \in \Upsilon^- (R , w^5/2,\sigma_-,\Omega) $
 for
some $ \Omega $ as described in the beginning of subsection $
\refe{main1} $, some $ \sigma_{\pm} > 0 $ and $ R  $ large enough.
Let us consider only one term of the right hand side of $
\reff{dev0} $ corresponding to a symbol supported in $ \Upsilon^+
(R,w^5 / 2,\sigma_+,\Omega) $  (the case $-$ is similar). Then, using
the trick $ \reff{shifttrick} $, lemma $ \refe{Lemme} $ and
Egorov's theorem we see that
\begin{eqnarray}
 \mbox{tr} \left( f
(H_{\e}^{(\kappa)}) U_{\e}^{(\kappa)}(t) O \! p_h (
\gamma_{\e,k}^{(\kappa)}) \right) \equiv_N^{n_1} \sum_{l \leq N} h^l
\mbox{tr} \left( f (H_{\e}^{(\kappa)}) U_{\e}^{(\kappa)}(t) O \!
p_h ( \tilde{\gamma}_{\e,l}^{(\kappa)}) \right) \label{dev1}
\end{eqnarray}
 with $
\tilde{\gamma}_{\e,l}^{(\kappa)} \in S_{\e}^{1,-\infty} $ supported in
$\Gamma_5^+$. Here again, we study only one term of the right hand
side of $ \reff{dev1} $.  We split the
integral corresponding to the inverse Fourier transform
$ \reff{coreFourierinverse} $  into two parts and consider first
\begin{eqnarray}
 \frac{1}{2 \pi h} \int_{0}^{+ \infty} e^{\frac{i}{h} \mu t
} \ \mbox{tr} \left( f (H_{\e}^{(\kappa)}) U_{\e}^{(\kappa)}(t)  O
\! p_h ( \tilde{\gamma}_{\e,l}^{(\kappa)})
 \right) \ \mbox{d} t \label{future} .
\end{eqnarray}
The interest of considering positive times and a symbol supported
in $ \Gamma^{+}_5 $ is that we can use $ \reff{Isozakitrick} $
with the results of subsection $ \refe{main1}$. Using the
notations of this subsection, we have
\begin{eqnarray}
\mbox{tr} \left( f (H_{\e}^{(\kappa)}) U_{\e}^{(\kappa)}(t)  O \!
p_h ( \tilde{\gamma}_{\e,l}^{(\kappa)})
 \right) \equiv_{N,+}^{n_1}
\mbox{tr} \left( f (H_{\e}^{(\kappa)}) J (\varphi^{+,(\kappa)}_{\e},
a_{\e}^{(\kappa)}) U (t) J (\varphi^{+,(\kappa)}_{\e},b_{\e}^{(\kappa)})^*
 \right) \label{dev2}
\end{eqnarray}
and, by centrality, we can rewrite the right hand side as
\begin{eqnarray}
\mbox{tr} \left( J (\varphi^{+,(\kappa)}_{\e},b_{\e}^{(\kappa)})^*
f (H_{\e}^{(\kappa)}) J
(\varphi^{+,(\kappa)}_{\e},a_{\e}^{(\kappa)})
 U (t)
 \right) \equiv_N^{n_1} \sum_{m \leq N} h^m \mbox{tr}  
 \left( O \! p_h ( c_{\e,m}^{(\kappa)}) U (t)
 \right)\label{dev3}
\end{eqnarray}
with $ c_{\e,m}^{(\kappa)} \in S_{\e}^{1,\infty} $. Now the contribution of
each term of the right hand side of $ \reff{dev3} $ to the
expected asymptotic is easy to obtain since
\begin{eqnarray}
 \lim_{\kappa \downarrow 0}  \frac{1}{2 \pi h} \! \int_{0}^{+ \infty} \! \! e^{\frac{i}{h} \mu t
} \ \left\{ \mbox{tr}  \left( O \! p_h ( c_{\e,m}^{(\kappa)}) U (t)
\right) \right\}_{\bf q}  \mbox{d}t = \frac{1}{2 \pi h} \!
\int_{0}^{+ \infty} \! \! \! \! \int \! \! \! \int \! \! \! \int
\! \! \! \int \! c_{m}(r,y,\rho,\eta) e^{\frac{i}{h}(\mu -
\rho^2)t} \mbox{d}r \mbox{d}y \mbox{d}\rho \mbox{d}\eta \
\mbox{d}t \nonumber
\end{eqnarray}
thanks to the easy and crucial remark that
\begin{eqnarray}
 c_m = \lim_{\kappa \rightarrow 0} \left\{  c_{\e,m}^{(\kappa)}
 \right\}_{\bf q} \in S_{\e}^{ q , - \infty} .
\end{eqnarray}
 The stationary phase theorem (in the variables
$t,\rho$) yields easily the asymptotic of the last integral in
integer powers of $h$ with coefficients which are smooth functions
of $\mu$ in the neighborhood of $ 1 $. Note that we use the fact
that $\rho$ is close to $1$ on the support of $c_m$ (see for instance
$\reff{rhononzero}$)
thus the
Hessian of the phase $ (\mu - \rho^2)t $ is non degenerate w.r.t.
$t,\rho$ at the stationary point $t=0, \rho = \mu^{1/2}$.

We still have to explain how to deal with the contribution of
negative times. Here we use the following trick due to Robert
\begin{eqnarray}
 \frac{1}{2 \pi h} \! \int_{-\infty}^{0} \! \! \! \! e^{\frac{i}{h} \mu t
} \ \mbox{tr} \left( f (H_{\e}^{(\kappa)}) U_{\e}^{(\kappa)}(t)  O
\! p_h ( \tilde{\gamma}_{\e,l}^{(\kappa)})
 \right) \ \mbox{d} t =
 \frac{1}{2 \pi h} \! \int_{0}^{+ \infty} \! \! \! \! e^{- \frac{i}{h} \mu t
} \ \overline{\mbox{tr} \left( \bar{f} (H_{\e}^{(\kappa)})
U_{\e}^{(\kappa)}(t) O \! p_h (
\tilde{\gamma}_{\e,l}^{(\kappa)})^*
 \right)} \ \mbox{d} t  \nonumber
\end{eqnarray}
which is a simple consequence of the fact that $ \mbox{tr}(A^*) =
\overline{\mbox{tr}(A)} $ combined with the centrality of the
trace. Then we can use the same method since
$$  \mbox{tr} \left( \bar{f} (H_{\e}^{(\kappa)})
U_{\e}^{(\kappa)}(t) O \! p_h (
\tilde{\gamma}_{\e,l}^{(\kappa)})^*
 \right) \equiv_N^{n_1} \sum_{m \leq N} h^m \mbox{tr} \left( \bar{f} (H_{\e}^{(\kappa)})
U_{\e}^{(\kappa)}(t) O \! p_h ( \check{\gamma}_{\e,m}^{(\kappa)})
 \right) $$
where the symbols $\check{\gamma}_{\e,m}^{(\kappa)}$ are
derivatives of the complex
conjugate of $ \tilde{\gamma}_{\e,l}^{(\kappa)} $. Since they are
supported in $ \Gamma^+_5 $ as well, we can repeat the same method
and the conclusion follows.  

This completes the proof under the condition that the coefficients of $ V $
vanish outside $ \{ r > R \} $. Otherwise, we proceed as follows. We write 
$ V = V_0 + V_{\infty} $ with $ V_{\infty} $ supported in $ \{ r > R \} $
and $ V_0 $ compactly supported outside $ \{ r > R + 1 \} $. The contribution of
$ V_{\infty} $ (or more precisely $ V_{\infty}^{(\kappa)}$) is treated as
before and thus we only have to consider 
$$ \mbox{tr}  \ \left( f(H_{\e}^{(\kappa)}) U_{\e}^{(\kappa)} (t) V_0
\tilde{f}(H_{\e}^{(\kappa)}) \right) . $$
Let $ K $  be the compact subset of $ T^* X $ defined by the conditions
 $ r \leq R + 1 $
and $ | p_{\e} - 1 | \leq w_5 / 4  $. 
By the non trapping condition, there exists $ T > 0 $ and $ R^{\prime} $ 
such that
$ \Pi_X \phi^{\pm T}_1 (K) \subset \{  R < r < R^{\prime} \} $, if $ \Pi_X $ is the
projection onto the base. Then we choose $ r_0 $ such that $ \Pi_X \phi^t_1(K) $
is disjoint from $ \{ r > r_0 \} $ for $ |t| \leq T $ and this ensure the fact that
$ \phi_{\e}^t = \phi_1^t $ for $ |t| \leq T $. Thus we can use the trick
$ \reff{shifttrick} $ with $ K_{\e} = V_0 \tilde{f}(H_{\e}^{(\kappa)}) $
combined with Egorov's theorem, and then repeat the same method. \finpreuve

\appendix

\section{Proof of lemma $ {\bf \refe{clefnegligeable}} $} \label{rappelmicrolocaux}
\setcounter{equation}{0} Let us first recall that for any smooth
and bounded function $a(x,x^{\prime},\xi)$ defined on $\Ra^{3n}$
the operator $A_h$, with Schwartz kernel $ {\mathcal A}_h
(x,x^{\prime}) $ defined as
$$ {\mathcal A}_h
(x,x^{\prime}) = (2 \pi h)^{-n} \int
e^{\frac{i}{h}(x-x^{\prime}).\xi} a (x,x^{\prime},\xi) \ \mbox{d}
\xi $$ is bounded on $ L^2 (\Ra^n) $ and we have the following
bound
\begin{eqnarray}
\left| \left| A_h \right| \right|_{\infty} \leq C
\max_{|\alpha+\alpha^{\prime}+\beta| \leq n_0} \sup_{\Ra^{3n}} \
\left|
\partial^{\alpha}_x
\partial^{\alpha^{\prime}}_{x^{\prime}} \partial^{\beta}_{\xi} a (x,x^{\prime},\xi)
\right|, \qquad h \in (0,1] \label{CalderonVaillancourt}
\end{eqnarray}
for some constants $C$ and $n_0$ independent of $a$ and $h$. This
is the usual Calderon-Vaillancourt's theorem.
This kind of operators preserve polynomial decay in the sense
that, for any $M \in \Ra$ we have
\begin{eqnarray}
 \left| \left| \scal{x}^{M} A_h \scal{x}^{-M} \right|
\right|_{\infty} \leq C_M \max_{|\alpha+\alpha^{\prime}+\beta|
\leq n_M} \sup_{\Ra^{3n}} \ \left|
\partial^{\alpha}_x
\partial^{\alpha^{\prime}}_{x^{\prime}} \partial^{\beta}_{\xi} a (x,x^{\prime},\xi)
\right|, \qquad h \in (0,1] \label{polynomialstable}
\end{eqnarray}
 for some $C_M$ and $ n_M $ depending
only on $M$. This is an easy consequence of $
\reff{CalderonVaillancourt} $.

Below, we shall use $ \reff{CalderonVaillancourt} $ as follows.
Assume that we have an operator $K_h^s$, with Schwartz kernel
${\mathcal K}_h^s$ of the form
$${\mathcal K}_h^s (r,y,r^{\prime},y^{\prime}) = (2 \pi h)^{-n} \int \! \! \int e^{\frac{i}{h}(r-r^{\prime})\rho
+ \frac{i}{h}(y-y^{\prime}).\eta} e^{\frac{i}{h} \psi^s
(r,y,r^{\prime},r,\rho,\eta)} a_h^s
(r,y,r^{\prime},y^{\prime},\rho,\eta) \ \mbox{d} \rho \mbox{d}
\eta
$$
with a smooth function $ \psi^s $ real valued and such that, on
the support of $a_h^s$,
$$ \left| \partial^{\gamma} \psi^s (r,y,r^{\prime},y^{\prime},\rho,\eta) \right| \leq C_{\gamma} \scal{s} , $$
for all $\gamma$. If we know moreover that for some $N \geq 0$ and
$M \geq 0$, we have
$$
 \left| \partial^{\gamma} a^s_h
(r,y,r^{\prime},y^{\prime},\rho,\eta) \right| \leq C_{\gamma} h^N
\scal{s}^{-M} \scal{r}^{-M} \scal{r^{\prime}}^{-M}
$$
then we obtain the estimate
\begin{eqnarray}
 \left| \left| \scal{r}^{M} K_h^s \scal{r}^{M}
\right| \right|_{\infty} \leq C h^{N - n_0} \scal{s}^{-M+n_{0}} .
\label{normeFourier}
\end{eqnarray}
This follows simply from $ \reff{CalderonVaillancourt} $ by
considering $ a = e^{i \psi^s / h} a_h^s $.

\medskip

For any $a(r,y,\rho,\eta)$, the kernel of $ J (\varphi_{\e}^+ , a)
U (s) J (\varphi_{\e}^+,b)^* $ is
\begin{eqnarray}
(2 \pi h)^{-n} \int \! \! \int e^{\frac{i}{h} \Phi_{\e}^s
(r,y,r^{\prime},y^{\prime},\rho,\eta) } a (r,y,\rho,\eta)
\overline{b (r^{\prime},y^{\prime},\rho,\eta)} \ \mbox{d}\rho
\mbox{d} \eta \label{noyauFourier}
\end{eqnarray}
where the phase function $ \Phi_{\e}^s $ is real valued and given
by
\begin{eqnarray}
 \Phi_{\e}^s (r,y,r^{\prime},y^{\prime},\rho,\eta) & = &
 \varphi_{\e}^+ (r,y,\rho,\eta) - s \rho^2 - \varphi_{\e}^+
 (r^{\prime},y^{\prime},\rho,\eta) . \nonumber
\end{eqnarray}
We shall use extensively the fact that for any $\gamma$,
$\partial^{\gamma} \left( \varphi_{\e}(r,y,\rho,\eta) - r \rho -
y.\eta \right) = {\mathcal O} (\varepsilon^k) $ on $\Gamma_k^+$.
This is a direct consequence of $\reff{phasesortante}$ and show
that, if $ \mbox{supp} \ a \subset \Gamma_2^+ $, $ \mbox{supp} \ b
\subset \Gamma_4^+ $, then
\begin{eqnarray}
\partial_{\rho} \Phi_{\e}^s (r,y,r^{\prime},y^{\prime},\rho,\eta)  & =
&  r - r^{\prime} - 2 s \rho + {\mathcal O} (\varepsilon^2)  \label{rhophase}  \\
\partial_{\eta} \Phi_{\e}^s (r,y,r^{\prime},y^{\prime},\rho,\eta)
& = & y - y^{\prime} + {\mathcal O} (\varepsilon^2)
\label{integrationeneta}
 \end{eqnarray}
 Recall moreover that, in view of $ \reff{rhononzero} $, there
 exists $C_1,c_1 > 0$ such that
\begin{eqnarray}
   c_1 < \rho < C_1 ,  \label{minmaxspeed}
\end{eqnarray}
on the support of $ a (r,y,\rho,\eta) \overline{b
(r^{\prime},y^{\prime},\rho,\eta)} $, if $ \mbox{supp} \ a \subset
\Gamma_2^+ $ and $ \mbox{supp} \ b \subset \Gamma_4^+ $.

Let us start the proof of $ \reff{supernegligeable} $. Here we
shall use the fact that $ a^{\prime}_{\e} $ in $ \Gamma_2^+
\setminus \Gamma_3^+ $. More precisely, we need to study symbols
$a (r,y,\rho,\eta) $ of the form
$$ \left( \partial_r^k \chi_{2,3}(r,y,\rho,\eta) \right) \partial^{\gamma} a_{\e} (r,y,\rho,\eta), \qquad
 \left( e^{- |\alpha| r} \partial_y^{\alpha} \chi_{2,3}(r,y,\rho,\eta) \right) \partial^{\gamma} a_{\e}(r,y,\rho,\eta) $$
 with $ k \geq 1$, $ |\alpha| \geq 1 $. We first consider the
 symbols $ a$ involving $\partial_r^k \chi_{2,3}$. \newline
 $ 1-$ {\it If $\partial_r$ falls on $ \chi_{R^3} $.}

 In this case, we
 have $ R^2 \leq r \leq R^3  $, hence we get the fast decay w.r.t.
 to $r$. Furthermore $r^{\prime} \geq R^4$ on the support of $b$,
 thus we have $ r - r^{\prime} \leq - R^4 + R^3  \ll 0 $ for $R$
 large, and we obtain
\begin{eqnarray}
 \partial_{\rho} \Phi_{\e}^s (r,y,r^{\prime},y^{\prime},\rho,\eta) \leq - 1 - 2 \rho s \leq - 1 - 2 c_1 s
 . \label{minorationtime}
\end{eqnarray}
 Thus we can integrate by part with $h \left( \partial_{\rho} \Phi_{\e}^s  \right)^{-1} D_{\rho}
 $ in $ \reff{noyauFourier} $, and we get as many powers of $h
 \scal{s}^{-1}$ as we want. \newline
$ 2-$ {\it If $\partial_r$ falls on $ \chi_{\varepsilon^3} $.}

Then, we have $ \varepsilon^3 \leq e^{-2r}g(y,\eta) \leq
\varepsilon^2 $, thus we get
$$   e^{-2 r^{\prime}} g
(y^{\prime},\eta) \leq
 \varepsilon^4 \ll \varepsilon^3 \leq  e^{-2r}g(y,\eta) .  $$
Since $  e^{- 2 r^{\prime}} |\eta|^2 \leq C_0 e^{-2 r^{\prime}} g
(y^{\prime},\eta) $ and $ e^{-2r}g(y,\eta) \leq C_0
e^{-2r}|\eta|^2 $, this implies easily that
$$  e^{2 r - 2 r^{\prime}} \leq C_0^2 \varepsilon . $$
Thus $ r - r^{\prime} = {\mathcal O} (\log \varepsilon) \ll 0 $.
In particular $ r \leq r^{\prime} $ and since we have as much
powers of $\scal{r^{\prime}}^{-1}$ as we want, we use the simple
fact that $ \scal{r} \scal{r^{\prime}}^{-1} $ is bounded to get as
many powers of $ \scal{r}^{-1} $ as we want. Furthermore, $
\reff{minorationtime} $ still holds, for $ \varepsilon $ small
enough, thus we can integrate by part as before and get as many
powers of $ h \scal{s}^{-1} $ as we wish. \newline $ 3-$ {\it If
$\partial_r$ falls on $ \chi_{w^3} $.}

Here we have $ | \rho^2 + e^{-2r}g(y,\eta) - 1 | \geq w^3 $. On
the other hand we also have
\begin{eqnarray}
\rho^2 + e^{-2r}g(y,\eta) & = & \rho^2 +
e^{-2r^{\prime}}g(y^{\prime},\eta) + \left( e^{-2r}g(y,\eta) -
e^{-2r^{\prime}}g(y^{\prime},\eta) \right) \nonumber \\ & = &
\rho^2 + e^{-2r^{\prime}}g(y^{\prime},\eta) + {\mathcal O}
(\varepsilon^{2})
\end{eqnarray}
with $ | \rho^2 + e^{-2r^{\prime}}g(y^{\prime},\eta) - 1 | \leq
w^4 $. Thus, on the support of $ a (r,y,\rho,\eta) \overline{ b
(r^{\prime},y^{\prime},\rho,\eta) } $, we obtain
\begin{eqnarray}
 | \rho^2 + e^{-2r}g(y,\eta) - 1 | \geq w^3 \qquad \mbox{and}
\qquad | \rho^2 + e^{-2r}g(y,\eta) - 1 | \leq w^4 + {\mathcal
O}(\varepsilon^2) . \label{contradictionsupport}
\end{eqnarray}
If we choose $ \varepsilon \leq w^2 $ small enough, then $
\reff{contradictionsupport} $ can not hold and then $ a
(r,y,\rho,\eta) \overline{ b (r^{\prime},y^{\prime},\rho,\eta) }
\equiv 0 $.

\medskip

We continue our analysis by considering symbols involving $
e^{-|\alpha|r} \partial_y^{\alpha} \chi_{2,3} $.

\medskip

\noindent $ 4- $ {\it If $ e^{-r} \partial_y $ falls on $
\chi_{\Omega_3} $.}

Then $y \in \Omega_2 \setminus \Omega_3$ and $ y^{\prime} \in
\Omega_4 \Subset \Omega_3 $, thus $|y-y^{\prime}| \geq c > 0$.
Using $\reff{integrationeneta}$, with $ \varepsilon  $ small
enough, we have $ |\partial_{\eta} \Phi_{\e}^s| \geq c / 2 $ and
we can integrate by part using $ h^2 |\partial_{\eta}
\Phi_{\e}^s|^{-2} \Delta_{\eta} $. This provides as many powers of
$h$ as we want. Furthermore the factor $e^{-r}$ yields the fast
decay w.r.t. $r$. We still need to explain how to get fast decay
w.r.t.  $s$. To that end, we introduce the following partition of
unit
\begin{eqnarray}
 1 = \theta_- \left( \frac{ r-r^{\prime} }{1+s} \right)
 + \theta_0  \left( \frac{ r - r^{\prime} }{ 1+s } \right)
 + \theta_+  \left( \frac{ r-r^{\prime}}{ 1+s } \right)
\label{partitiontimespace}
\end{eqnarray}
 Here $\theta_-$ is supported
in $ (- \infty , c_1) $, $\theta_+ $ in $ (3 C_1 , + \infty) $ and
we can assume that, on the support of $\theta_0$, we have
$$  c_1 / 2 \leq \frac{r-r^{\prime}}{1+s} \leq 4 C_1  . $$
On this support, we can write $ 1 = \scal{r-r^{\prime}}^M
\scal{r-r^{\prime}}^{-M} = \scal{r-r^{\prime}}^M {\mathcal O}
(\scal{s}^{-M}) $ , for any  $M$. Since we already have  fast
decay w.r.t. $ r $ and $ r^{\prime} $, the fast decay w.r.t. $s$
follows. On the support of $\theta_- ( (r - r^{\prime})  / (1 + s)
) $ (resp. $\theta_+$), using $ \reff{rhophase} $ and
$\reff{minmaxspeed}$ , we see that,
\begin{eqnarray}
 \partial_{\rho} \Phi_{\e}^s
(r,y,r^{\prime},y^{\prime},\rho,\eta) \leq -  c_1 s + c_1
 + {\mathcal O}( \varepsilon^2 )
 \qquad (\mbox{resp.} \  \geq C_1 s + 3 C_1 + {\mathcal O} ( \varepsilon^2 )
 ) \label{lowerupper}
\end{eqnarray}
 thus, for $s$ large enough, we can integrate by part as in $1-$
 and get as many powers of $ \scal{s}^{-1} $ as we want.
\newline
 $ 5- $ {\it If $ e^{-r} \partial_y $ falls on $
\chi_{\varepsilon^3} $.} We proceed as in $2-$ \newline $ 6- $
{\it If $ e^{-r} \partial_y $ falls on $ \chi_{w^3} $.} We proceed
as in $3- $

\medskip

This completes the proof of $ \reff{supernegligeable} $ using
$\reff{normeFourier}$.

\medskip
We now turn to the proof of $ \reff{justenegligeable} $. Here we
just need to get as many powers of $ \scal{r}^{-1} \scal{s}^{-1} $
as we want. The estimates will rely upon the fact that
\begin{eqnarray}
\left| \partial_r^k \partial_{\rho}^l \partial_y^{\alpha}
\partial_{\eta}^{\beta} P_{\e}(r,y,D_r,D_y) a_{\e}^{(N)} \right| \leq C_{k,l,\alpha,\beta} e^{-r}
\scal{\eta} \label{astucetime}
\end{eqnarray}
on $\Gamma^{+}_2$, for all $k,l,\alpha,\beta$, by proposition $
\refe{propquatre2} $. We proceed as follows. We introduce the
partition of unit $ \reff{partitiontimespace} $ and consider first
what happens on the support of $ \theta_0 ((r-r^{\prime})/(1+s))
$. Here we remark that
$$ e^{-r}\scal{\eta} \leq e^{- c_1 s / 2 - r^{\prime}} \scal{\eta} $$
where $ e^{-r^{\prime}}|\eta| $ is bounded on the support of $b$.
This yields an exponential decay in $ s $. Furthermore, we have $
r \leq r^{\prime} + 4 C_1 s $, thus we can write for any $M$
\begin{eqnarray}
 1 = \scal{r}^{-M} \scal{r}^{M} = \scal{r}^{-M} {\mathcal O}
\left( (\scal{r^{\prime}} + \scal{s})^M \right) . \label{Peetretime}
\end{eqnarray}
 The positive
powers of $\scal{r^{\prime}} $ and $ \scal{s} $ are respectively
controlled by the fast decay w.r.t. $ r^{\prime} $ of $b$ and the
exponential decay w.r.t. $s$, we obtain thus the expected decay.

On the support of $ \theta_- ((r-r^{\prime})/(1+s))  $, we have
similarly $ r - r^{\prime} \leq c_1 (1+s) $ and we can use $
\reff{Peetretime} $ again. This yields $ \scal{r}^{-M} $ for any
$M$ but we have to control $ \scal{s}^{M} $. By the same method as
in $4-$, using the upper bound given by $ \reff{lowerupper} $, we
get as many negative powers of $\scal{s}$ as we want.

The last step is to study what happens on the support of $
\theta_+ ((r-r^{\prime})/(1+s))  $. We have  $ r \geq r^{\prime} +
3 C_1 (s+1) $ and we get exponential decay in time, since
$$ e^{-r} \scal{ \eta }  \leq e^{-3 C_1 s - r^{\prime}} \scal{\eta } . $$
Furthermore, by choosing $ \varepsilon $ small enough, the lower
bound in $ \reff{lowerupper} $ shows that
$$ \partial_{\rho} \Phi_{\e}^s (r,y,r^{\prime},y^{\prime},\rho,\eta) \geq C_1 s + 2 C_1 > 0 . $$
Thus we can integrate by part  and get as many negative powers of
$ \partial_{\rho} \Phi_{\e}^{s} $ as we wish. Then we note that
$$ \left| \partial_{\rho} \Phi_{\e}^s (r,y,r^{\prime},y^{\prime},\rho,\eta) \right|^{-M} \leq
C \scal{r-r^{\prime}}^{-M} \scal{s}^M \leq C^{\prime}
\scal{r}^{-M} \scal{r^{\prime}}^M \scal{s}^{M}
$$
and the fast decay with respect to $r$ follows as before. This
completes the proof. \finpreuve

\section{Propagation estimates} \label{appendixuf}
\setcounter{equation}{0} In this appendix, we give sufficient
conditions leading to $\reff{propagationestimate}$.

\begin{prop} \label{A1} Assume that there are positive numbers $w,h_0,M,C$ such that
\begin{eqnarray}
 \sup_{\delta > 0} \ \left| \left| \scal{r}^{-M} R_{\e}(\mu \pm i
\delta ) \scal{r}^{-M} \right| \right|_{\infty} \leq C h^{-1},
\label{absorption}
\end{eqnarray}
 for all $\mu \in I $, $ h \in (0,h_0]$ and $\e \in [0,1]$. Then
 $\reff{propagationestimate}$ holds.
\end{prop}

\begin{prop} \label{A2} Assume that the manifold $(X,G)$ is non trapping,
and that the principal symbols of $ P_0 $ and $ P_1 $ coincide outside
$ \{ r > r_0 \}$. Let $ I $ be a neighborhood of $ 1 $. Then there
 exists $ r_0 $ large enough, $ h_0 $ small enough and 
 $C > 0$ such that $\reff{absorption}$ holds with $M=1$.
\end{prop}

The first proposition follows from Kato's theory of smooth
perturbations. We recall its simple proof for the sake of
completeness and to emphasize the uniformity with respect to the parameters.
The proof of the second one uses Mourre theory \cite{Mour1} and
more particularly a combination of the ideas of \cite{FrHi0} and
\cite{GeMa0,Robe1}.

\medskip


 \noindent {\it Proof of proposition $
\refe{A1}$.} Since we can always write $f=f_1 f_2$ with $f_1,f_2
\in C^{\infty}_0(I)$ and
$$ || \scal{r}^{-M} f(H_{\e}) U_{\e} (t)
\scal{r}^{-M} ||_{\infty} \leq || \scal{r}^{-M} f_1 (H_{\e}) U_{\e}(t)
||_{\infty} \ || ( \scal{r}^{-M} f_{2} (H_{\e}) U_{\e}(-t) )^* ||_{\infty}
$$ it is enough to estimate the norm of $ || \scal{r}^{-M} f (H_{\e})
U_{\e}(t) ||_{\infty} $ in $L^2(\Ra,\mbox{d}t)$.
Let $ A_{\e}=\scal{r}^{-M}f(H_{\e}) $, then by Parseval's
identity (see \cite{Kato0,RS4}), we have, for all $u \in L^2(X)$,
\begin{eqnarray}
 2 \pi h^{-1} \int_{\Ra} e^{-2
\delta t / h }||A_{\e} U_{\e}(t)u||^2 \ \mbox{d}t =  \int_{\Ra} ||A_{\e}
(R_{\e}(\mu - i \delta)-R_{\e}(\mu+i\delta))u||^2  \ \mbox{d}\mu .
\label{Parseval}
\end{eqnarray}
 On the other hand,  $  ( i /
2)(R_{\e}(z)-R_{\e}(\bar{z})) =  (\im z) \
R_{\e}(z)R_{\e}(\bar{z}) $
defines a nonnegative operator for $\im z > 0$ and we denote its
positive square root by $K_{\e}(z)$. We have trivially $ ||K_{\e}(z)A_{\e}^*
u||^2 = (2 i \pi)^{-1} \left( A_{\e}^* u , (R_{\e}(z) -
R_{\e}(\bar{z}))A_{\e}^* u \right) $ and this yields
\begin{eqnarray}
4  \int_{\Ra} ||A K_{\e}(\mu+ i \delta )^2 u ||^2 \ \mbox{d} \mu \leq
 \tilde{C} h^{-1} \int_{\Ra} ||K_{\e}(\mu+i\delta)u||^2 \
\mbox{d} \mu = \tilde{C} h^{-1} \pi ||u||^2 \label{root}
\end{eqnarray}
since $ ||A_{\e} K_{\e}(\mu+ i \delta )^2 u ||  \leq
|| K_{\e}(\mu+ i \delta ) A^*_{\e} ||
|| K_{\e}(\mu+ i \delta ) u || $. Here
 $\tilde{C}$ depends only on  $C$ in $ \reff{absorption} $
and $f$, and is uniform w.r.t. $\e,h$ and $\delta$. The left
hand side of $ \reff{root} $ is nothing but the right hand side of
$ \reff{Parseval} $ and the result follows easily. \finpreuve

\medskip

\noindent {\it Proof of proposition $ \refe{A2}$.} The
estimate $ \reff{absorption} $ follows directly from the method of
\cite{Mour1} provided the following Mourre estimate holds, with
$E_{\e}(I) $ the spectral projector of $H_{\e}$ on $I$,
\begin{eqnarray}
   E_{\e}(I) i \left[
H_{\e} , A_h \right] E_{\e} (I) \geq c h E_{\e} (I) \label{Mourre}
\end{eqnarray}
 for some constant $c > 0$ independent
of $h \in (0,h_0]$ and $\e \in [0,1] $. If $ \reff{Mourre} $ holds
then, the theory of Mourre  shows that $ || \
|A_h+i|^{-1}R_{\e}(\mu+i\delta) |A_h+i|^{-1} \ ||_{\infty} =
{\mathcal O}( (ch)^{-1} ) $ thus $ \reff{absorption} $ holds, provided
\begin{eqnarray}
 || \scal{r}^{-1} A_h ||_{\infty} \leq C_1 , \qquad h \in (0,h_0
] . \label{rborne}
\end{eqnarray}
 This reduces the proof to the construction of $A_h$. We now
 sketch its construction and refer to \cite{FrHi0,GeMa0,Robe1}
 for the details. The idea is to construct $A_h = A_h^0 + A_h^{\infty}$
with $A_h^{\infty}$ supported near infinity and $A_h^0$ compactly supported.
Following \cite{FrHi0}, we define $A^{\infty}_h$ as
\begin{eqnarray}
 A^{\infty}_h = \frac{1}{2} \tilde{f} (H_1)
\left( \theta^2 \omega_S^2 u h D_r + h D_r \theta^2 \omega_S^2 u  \right)
 \tilde{f} (H_1)   \label{ainfini}
\end{eqnarray}
  with $\tilde{f} \in C_0^{\infty}$, $\tilde{f}= 1$ close
  to $1$, and $\theta = \theta (r/R) $  smooth, bounded and supported
near infinity, say in $r \geq R/2$. We also have
$$ \omega_S = \omega \left( \frac{2 r - \log (h^2 \Delta_Y + 1) }{S} \right) ,
\qquad u = 2 r + S - \log (h^2 \Delta_Y + 1) , \qquad S > 0, $$
where the function $ \omega \in C^{\infty}(\Ra) $ is supported in
$[-1,\infty)$ and such that $\omega \equiv 1$ on $[-1/2,\infty)$.
Then using the calculations of \cite{FrHi0} and pseudo-differential calculus,
we see that
\begin{eqnarray}
 f (H_{\e})   i [H_{\e} , A^{\infty}_h] f (H_{\e}) \geq h
 \theta \omega_S f^2 (H_{\e}) \omega_S \theta + \iota(R,S)
{\mathcal O}(h ) + {\mathcal O}(h^2) \label{infini}
\end{eqnarray}
if $f$ is supported where $\tilde{f} = 1$. Here  $ \iota (R,S)
\downarrow 0 $ as $R,S \rightarrow \infty$ and the notation ${\mathcal O}(h^k)$
holds in operator norm, uniformly w.r.t $\e$. Note that the
spectral cutoff
$\tilde{f}(H_1)$ in $ \reff{ainfini} $ doesn't commute
with $H_{\e}$
and we use the fact, among other ones, that
$ ( \tilde{f}(H_1) - \tilde{f}(H_{\e})) \theta h D_r =
{\mathcal O} (R^{- \infty})  $.

 If there was no term
$\omega_S$ in the right hand side of $\reff{infini}$ , we would have done half
of our program.
How to neglect $ \theta (1 - \omega_S)f(H_{\e})$? We can give a
pseudo-differential expansion of $\omega_S$ similar to the one given in
section $ \refe{pseudo} $ and thus, up to an operator which is
${\mathcal O}(h)$, we can replace $ 1- \omega_S $ by a pseudo-differential
operator
with principal symbol $ 1- \omega ( (2 r - \log(g(y,\eta)+1))/S) $. On its
support,
we have  $e^{-2 r} g (y,\eta) + e^{-2r} \geq e^{S/2} $.
 On the other hand, on the  support of the principal symbol of
$ f(H_{\e})$ we have $ e^{-2 r} g_{\e}(r,y,\eta) < 3/2 $ (if $f$ is supported
close to $1$)
and thus $ e^{-2r} g (y,\eta) < 3 / 2 + {\mathcal O}(e^{-r}) $. All this shows
that the principal symbol of  $ \theta (1 - \omega_S)f(H_{\e}) $ is
identically $0$ for all $R$ and $S$ large enough, which implies that
\begin{eqnarray}
 \theta (1 - \omega_S) \tilde{f}(h^2 P_{\e}) = {\mathcal O}(h)
 \label{semiclassique}
\end{eqnarray}
where one should notice that ${\mathcal O}$ depends on $R$ and $S$. Then we get
\begin{eqnarray}
 f (H_{\e})   i [H_{\e} , A^{\infty}_h] f (H_{\e}) \geq h
 \theta  f^2 (H_{\e})  \theta + \iota(R,S)
{\mathcal O}(h ) + {\mathcal O}(h^2) \label{infini2}
\end{eqnarray}
where ${\mathcal O}(h^2)$ depends on $R,S$, but not ${\mathcal O}(h)$.

We will now construct $A_{h}^0$ by the method of \cite{GeMa0,Robe1}, such that
\begin{eqnarray}
 f (H_{\e}) i [H_{\e}, A^0_h] f (H_{\e}) =  h \widetilde{\theta} f(H_{\e})
\widetilde{\theta}
+ \tilde{\iota}( R , K ) {\mathcal O}( h) + {\mathcal O}(h^2) , \label{compactMourre}
\end{eqnarray}
with $\widetilde{\theta} \in C_0^{\infty} $  such that $ \theta^2 +
\widetilde{\theta}^2 = 1 $, and $\tilde{\iota}( R , K ) \downarrow 0 $
as $ R, K \uparrow \infty $. Here $ K $ is another large parameter introduced
below.
If this holds, we easily get $ \reff{Mourre} $ by summing
$ \reff{infini2} $ and $ \reff{compactMourre} $, choosing first
 $R$, $ K $ and $S$
large enough, then $h$ small enough and then by multiplying both side by $E_{\e}(I)$.
Note that we also use the fact that $[f(H_{\e}),\theta] = {\mathcal O}(h)$.

The idea is to define $A_h^0$ as a (bounded) pseudo-differential operator
whose principal symbol is the following function $ a $, which is invariantly
defined on $T^*X$
$$ a = \tilde{\chi} \int_0^{\infty} \widetilde{\theta}^2 \circ \phi^t_{1} \
\mbox{d} t
\ \tilde{f} \circ p_{1}  $$
where $ p_1 $ and $ \phi_1^t$ are the principal symbol
of $H_1$ and the associated flow, and $\tilde{\chi} = \tilde{\chi}(r/K)$ is a
 $ C_0^{\infty} $ function such that $\tilde{\chi} \widetilde{\theta}^2 =
\widetilde{\theta}^2$. Note that the integral is  convergent thanks to the non
trapping condition. Then the crucial remark  is that the Poisson bracket is
$$ \{ p_{\e} , a \}
 = \tilde{f} (p_{1} ) \widetilde{\theta}^2  + {\mathcal O}(K^{-1}) +
 {\mathcal O} (R^{- \infty})
$$
with ${\mathcal O}(K^{-1})$ uniform w.r.t $\e$ (but not w.r.t $R$)
in the topology of smooth and bounded functions.
 This follows first  from the fact that
$$ \{ p_{1} , \widetilde{\theta}^2 \circ \phi_1^t \} =
\{ p_{1} , \widetilde{\theta}^2  \} \circ \phi_1^t =
\frac{d}{d t} \left( \widetilde{\theta}^2 \circ \phi_1^t  \right)   $$
and from the fact that $p_{\e} = p_{1}$ for
$r \leq r_0 $ with $ r_0 $ which we can choose $ \geq R $. 
Then $ \reff{compactMourre} $ follows from semi-classical
pseudodifferential calculus. (see \cite{Robe2}). \finpreuve

\end{document}